\documentclass[leqno,10pt]{article}
\usepackage{amssymb,amsmath}
\usepackage[all]{xy}

\newenvironment{proof}{{\bf Proof. }}{\par}{\bigskip}
\newtheorem{theo}{Theorem}[section]
\newtheorem{defi}[theo]{Definition}

\newtheorem{con}[theo]{Convention}
\newtheorem{lem}[theo]{Lemma}

\newtheorem{lemm}[theo]{Lemma}
\newtheorem{prop}[theo]{Proposition}

\newtheorem{rem}[theo]{Remark}

\newtheorem{coro}[theo]{Corollary}

\newtheorem{exam}[theo]{Example}
\newtheorem{assu}[theo]{Assumption}

\newcommand{\kgot}{\ensuremath{\mathfrak{k}}}
\newcommand{\ggot}{\ensuremath{\mathfrak{g}}}
\newcommand{\sogot}{\ensuremath{\mathfrak{so}}}

\newcommand{\so}{\ensuremath{\mathfrak{so}}}
\newcommand{\ugot}{\ensuremath{\mathfrak{u}}}
\newcommand{\spingot}{\ensuremath{\mathfrak{spin}}}

\newcommand{\bc}{\ensuremath{{\mathbf c}}}

\newcommand{\Acal}{\ensuremath{\mathcal{A}}}
\newcommand{\Ccal}{\ensuremath{\mathcal{C}}}
\newcommand{\Dcal}{\ensuremath{\mathcal{D}}}
\newcommand{\Fcal}{\ensuremath{\mathcal{F}}}
\newcommand{\Ecal}{\ensuremath{\mathcal{E}}}

\newcommand{\Hcal}{\ensuremath{\mathcal{H}}}
\newcommand{\Kcal}{\ensuremath{\mathcal{K}}}
\newcommand{\Lcal}{\ensuremath{\mathcal{L}}}
\newcommand{\Pcal}{\ensuremath{\mathcal{P}}}
\newcommand{\Scal}{\ensuremath{\mathcal{S}}}
\newcommand{\Ucal}{\ensuremath{\mathcal{U}}}
\newcommand{\Vcal}{\ensuremath{\mathcal{V}}}

\newcommand{\Zcal}{\ensuremath{\mathcal{Z}}}

\newcommand{\Cbb}{\ensuremath{\mathbb{C}}}
\newcommand{\Rbb}{\ensuremath{\mathbb{R}}}
\newcommand{\Zbb}{\ensuremath{\mathbb{Z}}}
\newcommand{\Lbb}{\ensuremath{\mathbb{L}}}

\newcommand{\A}{\ensuremath{\mathbb{A}}}
\newcommand{\F}{\ensuremath{\mathbf{F}}}
\newcommand{\x}{\ensuremath{\mathbf{x}}}

\newcommand{\bas}{\ensuremath{\hbox{\scriptsize \rm bas}}}
\newcommand{\hor}{\ensuremath{\hbox{\scriptsize \rm hor}}}

\newcommand{\dr}{\ensuremath{\hbox{\scriptsize \rm dec-rap}}}
\newcommand{\KK}{\ensuremath{{\mathbf K}^0}}

\newcommand{\rel}{\ensuremath{\rm rel}}

\newcommand{\f}{\ensuremath{\mathcal{C}^{\infty}}}

\newcommand{\T}{\ensuremath{\hbox{\bf T}}}
\newcommand{\End}{\ensuremath{\hbox{\rm End}}}
\newcommand{\sign}{\operatorname{ sign}}
\newcommand{\Pf}{\operatorname{Pf}}
\newcommand{\jdemi}{\ensuremath{{\rm j}^{1/2}}}
\newcommand{\str}{\operatorname{Str}}
\newcommand{\herm}{\operatorname{Herm}}
\newcommand{\tr}{\operatorname{Tr}}
\newcommand{\ch}{\operatorname{Ch}}
\newcommand{\chs}{\operatorname{Ch_{\rm sup}}}
\newcommand{\chg}{\operatorname{Ch_{\rm sup}}}
\newcommand{\chc}{\operatorname{Ch_{c}}}
\newcommand{\chcf}{\operatorname{Ch_{{\ensuremath{\hbox{\scriptsize \rm fiber cpt}}}}}}
\newcommand{\chq}{\operatorname{Ch_{\rm Q}}}

\newcommand{\Eul}{\operatorname{Eul}}

\newcommand{\cf}{\ensuremath{\hbox{\scriptsize \rm fiber cpt}}}

\newcommand{\chr}{\operatorname{Ch_{\rm rel}}}
\newcommand{\supp}{\operatorname{Supp}}
\newcommand{\p}{\operatorname{p}}
\newcommand{\bere}{\operatorname{T}}
\newcommand{\pol}{\operatorname{pol}}
\newcommand{\Ex}{\operatorname{Ex}}
\newcommand{\e}{\operatorname{e}}

\newcommand{\res}{\operatorname{\bf r}}

\def \cf {\hbox{\scriptsize \rm fiber cpt}}
\def \xx {\mathbf{x}}
\def \clif {\mathbf{c}}
\def \Id {{\rm Id}}

\def \tuc {{\rm Th_c}}
\def \tucf {{\rm Th_{{\ensuremath{\hbox{\scriptsize \rm fiber cpt}}}}}}
\def \tur {{\rm Th_{{\ensuremath{\rm rel}}}}}

\def \tumq {{\rm Th_{{\ensuremath{\rm MQ}}}}}

\def \todd {{\rm Todd}}
\def \pf {{\rm Pf}}
\def \cst {{\rm cst}}
\def \sm  {{\rm m}}

\def \spin {{\rm Spin}}
\def \spinc {{\rm Spin}^{\rm c}}
\def \So {{\rm SO}}

\def \Det {{\rm Det}}

\def \clif {{\bf c}}

\title{Equivariant relative Thom forms and Chern characters}
\author{ Paul-Emile Paradan, Mich\`ele Vergne}

\date{December 2007}

\begin{document}

\maketitle

 {\small
 \tableofcontents}

\section{Introduction}

These notes are the first chapter of a monograph, dedicated to a
detailed proof of the equivariant index theorem for transversally
elliptic operators.

In this preliminary chapter, we  prove a certain number of natural
relations in equivariant cohomology.  These relations include the
Thom isomorphism in equivariant cohomology, the multiplicativity
of the relative Chern characters, and the Riemann-Roch  relation
between the relative Chern character of the Bott symbol and of the
relative Thom class. In the spirit of Mathai-Quillen, we  give
``explicit" representatives of a certain number of relative
classes. We believe that this
construction has its interest in the non equivariant case as well.
As remarked by Cartan (see \cite{gui-ste99}) and emphasized in
Mathai-Quillen \cite{Mathai-Quillen}, computations in the ordinary de Rham
cohomology of vector bundles are deduced easily from computations
in the equivariant cohomology of vector spaces. In particular, we
give here an explicit formula for the relative Thom form
$\tur(\Vcal)\in \Hcal^*(\Vcal,\Vcal\setminus M)$ of a Euclidean
vector bundle $p:\Vcal\to M$ provided with a Euclidean connection:
$\tur(\Vcal):=[p^*\Eul(\Vcal),\beta_{\Vcal}]$, where $\Eul(\Vcal)$
is the Euler class and  $\beta_{\Vcal}$  an explicit form
depending of the connection, defined outside the zero section of
$\Vcal$, such that $p^*\Eul(\Vcal)=d\beta_{\Vcal}$. We similarly give
an explicit formula for the relative Chern character
$\chr(\sigma_b)\in \Hcal^*(\Vcal,\Vcal\setminus M)$  of the Bott
morphism on the  vector bundle $\Vcal$, if $\Vcal$ is provided
with a complex structure. The Riemann-Roch relation
$$
p^*\left({\rm Todd}(\Vcal)\right)\chr(\sigma_b)=\tur(\Vcal)
$$
holds in relative cohomology, and follows from the formulae.

Our constructions in the de Rham model for equivariant cohomology
are strongly influenced by Quillen's construction of
characteristic classes via super-connections and super-traces. The
articles of Quillen \cite{Quillen85} and Mathai-Quillen
\cite{Mathai-Quillen} are our main background. However, Quillen
did not use relative cohomology while  our constructions are
systemically performed in {\bf relative cohomology}, therefore are
more precise. This relative construction was  certainly present in
the mind of Quillen, and we do not pretend to a great originality.
Indeed, if a morphism $\sigma:\Ecal^+\to \Ecal^-$ between two
vector bundles over a  manifold $M$ is invertible outside a closed
subset $F$, the construction of Quillen of the Chern character
$\chq(\sigma)=\str(e^{\A^2_{\sigma}})$ is defined using a super-connection
$\A_{\sigma}$ with zero degree term the odd endomorphism $i(\sigma\oplus
\sigma^*)$ and  this construction provides also a form $\beta$
defined outside $F$ such that the equality
$\ch(\Ecal^+)-\ch(\Ecal^-)=d\beta$ holds on $M\setminus F$. Thus the couple
$\left(\alpha,\beta\right)$ of differential forms, with
$\alpha:=\ch(\Ecal^+)-\ch(\Ecal^-)$, defines naturally a class $\chr(\sigma)$
in the de Rham relative group $\Hcal^*(M,M\setminus F)$ that we
call the Quillen's relative Chern character. Now, if
$a=\left(\alpha,\beta\right)$ is a closed element in relative cohomology, e.g.
$\alpha$ is a closed form on $M$ and $\alpha=d\beta$ outside $F$, the couple $(\alpha,\beta)$ of
differential forms leads naturally to usual de Rham closed differential forms
on $M$ with support as close as we want from $F$. Indeed, using a
function $\chi$ identically equal to $1$ on a neighborhood of $F$,
the  closed differential form $\p(a):=\chi \alpha+d\chi \beta$ is
supported as close as we want of $F$. Thus Quillen's
super-connection construction gives us three representations of
the Chern character: the Quillen Chern character $\chq(\sigma)$,
the relative Chern character $\chr(\sigma)$ and the Chern
character $\chs(\sigma)=\p(\chr(\sigma))$ supported near $F$. We
study the relations between these classes and prove some basic
relations. Our previous article \cite{pep-vergne1} explained the
construction of the relative Chern character in ordinary
cohomology. Here these constructions are done in equivariant cohomology and are very similar.

As an important example, we consider $\sigma_b$  the Bott morphism
on a complex  vector bundle $\sigma_b:\wedge^+\Vcal\to
\wedge^-\Vcal$ over $\Vcal$, given by the exterior product by
$v\in \Vcal$. This morphism has support the zero section $M$ of
$\Vcal$. Quillen's Chern character $\chq(\sigma_b)$ is
particularly pleasant as it is represented by a differential form
with ``Gaussian look" on each fiber of $\Vcal$. However, for many
purposes, it is important to construct the Chern character of
$\sigma_b$, as a differential form  supported near the zero
section of $\Vcal$. More precisely, we here consider
systematically the relative class $\chr(\sigma_b)$ in
$\Hcal^*(\Vcal,\Vcal\setminus M)$ which contains all information.

A similar construction of the Thom form, using the Berezin
integral instead of a super-trace, leads to explicit formulae for
the relative Thom class. Again here, we have three representatives
of the Thom classes, the Mathai-Quillen Thom form which has a
``Gaussian look", the relative Thom form, and the Thom form with
support near the zero section. Our main result is Theorem
\ref{theo:thom-form} where these three representatives are given in the
equivariant cohomology of vector spaces.

 These explicit formulae allows us to
derive the well known relations between Thom classes in cohomology
and $K$-theory (Theorem \ref{prop:chern-V-spin}). Here again,
following Mathai-Quillen construction, we perform all calculations
on the equivariant cohomology groups of an Euclidean vector space,
and we apply Chern-Weil morphism to deduce relations in any vector
bundle.

We also include in this chapter proofs of Thom isomorphisms in
various equivariant cohomologies spaces,  using  Atiyah's
``rotation'' construction. For the case of relative cohomology, we
need to define the product in de Rham relative cohomology and some
of its properties. This is the topic of Section
\ref{section:coho-relative}.

In the second chapter of this monograph, we will generalize our
results to equivariant cohomology classes with
$C^{-\infty}$-coefficients. This will be an essential ingredient
of our new index formula for transversally elliptic operators (see
\cite{VergneICM}).

\bigskip

\section{Equivariant cohomologies }\label{sec:chg-sigma}

\subsection{Definitions and notations}

If $f$ is a map on a space $M$, the notation $f(x)$ means,
depending of the context, either the value of $f$ at the
point $x$ of $M$, or the function $x\mapsto f(x)$ where $x$ is a running
variable in $M$.

When a compact Lie group $K$ acts linearly on a vector space $E$, we denote
$E^K$ the sub-space of $K$-invariant vectors.

Let $N$ be a manifold, and let $\Acal(N)$ be the algebra of
differential forms on $N$. We denote by $\Acal_c(N)$ the
sub-algebra of compactly supported differential forms. We will
consider on $\Acal(N)$ and $\Acal_c(N)$ the $\Zbb$-grading by the
exterior degree. It induces a $\Zbb_2$-grading  on $\Acal(N)$ in
even or odd forms. We denote by $d$ the de Rham differential. If
$\alpha$ is a closed differential form, we sometimes denote also
by $\alpha$ its de Rham cohomology class.

If $S$ is a vector field on $N$, we denote by
$\Lcal(S):\Acal^k(N)\to \Acal^{k}(N)$ the Lie derivative, and by
$\iota(S):\Acal^k(N)\to \Acal^{k-1}(N)$  the contraction of a
differential form by the vector field $S$.

Let $K$ be a  compact Lie group with Lie algebra $\kgot$. We
denote  $\Ccal^{\pol}(\kgot)$
 the space of polynomial functions on $\kgot$ and $\Ccal^{\infty}(\kgot)$
 the space of $C^{\infty}$-functions on $\kgot$.
 The algebra   $\Ccal^{\pol}(\kgot)$ is
isomorphic to the symmetric algebra $S(\kgot^{*})$ of $\kgot^*$.

 We suppose that the manifold $N$ is provided with an
action of $K$. We denote $X\mapsto VX$ the corresponding morphism
from $\kgot$ into the Lie algebra of vector fields on $N$: for
$n\in N$,
$$
V_n X:=\frac{d}{d\epsilon} \exp(-\epsilon X)\cdot n|_{\epsilon=0},\quad X\in\kgot.
$$

%

Let $\Acal^{\pol}(\kgot, N)=(\Ccal^{\pol}(\kgot)\otimes
\Acal(N))^K$ be the $\Zbb$-graded algebra of equivariant
polynomial functions $\alpha: \kgot\to\Acal(N)$. Its
$\Zbb$-grading is the grading induced by the exterior degree and
where elements of $\kgot^*$ have degree two. Let $D= d-\iota(VX)$
be the equivariant differential:
$$(D\alpha)(X)=d(\alpha(X))-\iota(VX)\alpha(X).$$
Let $\Hcal^{\pol}(\kgot,N):= \mathrm{Ker} D/ \mathrm{Im} D$ be the
equivariant cohomology algebra with polynomial coefficients. It is
a module  over  $\Ccal^{\pol}(\kgot)^K$.

\begin{rem}If $K$ is not connected, $\Acal^{\pol}(\kgot, N)$
depends of $K$, and not only of the Lie algebra of $K$. However,
for notational simplicity, we do not include $K$ in the notation.
\end{rem}

If $g:M\to N$ is a $K$-equivariant map from the $K$-manifold $M$
to the $K$-manifold $N$, then we obtain a map
$g^*:\Acal^{\pol}(\kgot,N)\to \Acal^{\pol}(\kgot,M)$, which
induces a map $g^*$ in cohomology. When $U$ is an open invariant
subset of $N$, we denote by $\alpha\mapsto \alpha|_U$ the restriction
of $\alpha\in \Acal^{\pol}(\kgot,N)$ to $U$.

If $S$ is a $K$-invariant vector field on $N$, the operators $\Lcal(S)$ and $\iota(S)$ are
extended from $\Acal(N)$ to $\Ccal^{\pol}(\kgot)\otimes
\Acal(N)$: they commute with the $K$-action, so $\Lcal(S)$ and $\iota(S)$  act on
$\Acal^{\pol}(\kgot,N)$. Cartan's relation
holds:
\begin{equation}\label{Cartan}
\Lcal(S)=D\circ\iota(S)+\iota(S)\circ D.
\end{equation}

If $N$ is non-compact, we can also consider
 the space  $\Acal^{\rm pol}_c(\kgot,N):=(\Ccal^{\pol}(\kgot)\otimes \Acal_c(N))^K$
  of equivariant polynomial forms $\alpha(X)$
 which are  compactly supported on $N$. We denote by
$\Hcal^{\rm pol}_c(\kgot,N)$ the corresponding cohomology algebra.
If $N$ is an oriented manifold, integration over $N$ defines a map
$\Hcal^{\rm pol}_c(\kgot,N)\to \Ccal^{\pol}(\kgot)^K$. If  $\pi:
N\to B$ is a $K$-equivariant fibration  with oriented fibers, then the integral
over the fiber defines a map $\pi_*:\Hcal^{\pol}_{c}(\kgot,N)\to
\Hcal^{\pol}_c(\kgot,B)$.

Finally, we give more definitions in the case of a $K$-equivariant real vector
bundle $p:\Vcal\to M$. We may define two sub-algebras of $\Acal^{\pol}(\kgot, \Vcal)$ which are
stable under the derivative $D$.
The sub-algebra $\Acal^{\pol}_{\dr}(\kgot, \Vcal)$
consists of  polynomial equivariant forms on $\Vcal$ such that all
partial derivatives are rapidly decreasing along the fibers. We may also
consider the sub-algebra
$\Acal^{\pol}_{\cf}(\kgot,\Vcal)$
of $K$-equivariant forms on $\Vcal$ which have a compact support
in the fibers of $p:\Vcal\to M$. The inclusions
$\Acal^{\pol}_{\cf}(\kgot,\Vcal)\subset\Acal^{\pol}_{\dr}(\kgot, \Vcal)\subset \Acal^{\pol}(\kgot,\Vcal)$
give rise to the natural maps
$\Hcal^{\pol}_{\cf}(\kgot,\Vcal)\to\Hcal^{\pol}_{\dr}(\kgot, \Vcal)\to \Hcal^{\pol}(\kgot,\Vcal)$.
If the fibers of $\Vcal$ are oriented, integration over the fiber defines a map
$p_*:\Hcal^{\pol}_{\dr}(\kgot,\Vcal)\to \Hcal^{\pol}(\kgot,M)$.

\bigskip

Let $\Acal^{\infty}(\kgot, N)$ be the $\Zbb_2$-graded algebra of
equivariant {\em smooth} maps $\alpha: \kgot\to\Acal(N)$. Its
$\Zbb_2$-grading is the grading induced by the exterior degree.
The equivariant differential $D$ is well defined on
$\Acal^{\infty}(\kgot, N)$ and respects the $\Zbb_2$-grading. Let
$\Hcal^{\infty}(\kgot,N):= \mathrm{Ker} D/ \mathrm{Im} D$ be the
corresponding cohomology algebra with $C^{\infty}$-coefficients. We
denote by $\Acal^{\infty}_c(\kgot, N)$ the sub-algebra of
equivariant differential forms with compact support and by
$\Hcal_c^{\infty}(\kgot, N)$ the corresponding algebra
cohomology  algebra .  Then  $\Hcal^{\infty}(\kgot,N)$ and
$\Hcal^{\infty}_c(\kgot,N)$  are $\Zbb_2$-graded algebras. If $N$
is oriented, integration over $N$ defines a map
$\Hcal^{\infty}_c(\kgot,N)\to \Ccal^{\infty}(\kgot)^K$.

Let $\Vcal$ be a $K$-equivariant real vector bundle over a manifold $M$. We define
similarly $\Acal^{\infty}_{\dr}(\kgot, \Vcal)$ and
$\Hcal^{\infty}_{\dr}(\kgot, \Vcal)$ as well as
$\Acal^{\infty}_{\cf}(\kgot,\Vcal)$ and
$\Hcal^{\infty}_{\cf}(\kgot,\Vcal)$. There are natural maps
$\Hcal^{\infty}_{c}(\kgot, \Vcal)\to\Hcal^{\infty}_{\cf}(\kgot, \Vcal)
\to$ \break $\Hcal^{\infty}_{\dr}(\kgot,\Vcal)$ and an integration map
$\Hcal^{\infty}_{\dr}(\kgot, \Vcal) \to
\Hcal^{\infty}(\kgot,M)$ if the fibers of $\Vcal\to M$ are oriented.

After these lengthy definitions, we hope that at this point the
reader  is still with us.

\bigskip

If $K$ is the identity (we say the non-equivariant case), then the
operator $D$ is the usual de Rham differential $d$. We
systemically skip the letter $\kgot=\{0\}$ in the corresponding
notations. Thus the equivariant cohomology group
$\Hcal^{\pol}(\kgot,N)$ coincide with the usual de Rham cohomology
group  $\Hcal(N)$.  The compactly supported cohomology space is
denoted by  $\Hcal_c(N)$ and the rapidly decreasing cohomology
space by $\Hcal_{\dr}(\Vcal)$. (In this article, we will only work
with cohomology groups, so the notation $\Hcal$ refers always to
cohomology).

\subsection{Equivariant cohomology of vector bundles}\label{section:coho-vectorbundles}

It is well known that the cohomology of a vector bundle is the
cohomology of the basis. The same equivariant Poincar\'{e} lemma holds
in equivariant cohomology (see for example \cite{Duflo-Vergne2}).
We review the proof.

Let $p:\Vcal\to M$ be a $K$-equivariant (real) vector bundle. Let
$i:M\to \Vcal$ be the inclusion of the zero section.

\begin{theo}
$\,$

 $\bullet$ The map $i^*: \Hcal^{\pol}(\kgot, \Vcal)\to
\Hcal^{\pol}(\kgot,M)$ is an isomorphism with inverse $p^*$.

 $\bullet$ The map $i^*: \Hcal^{\infty}(\kgot, \Vcal)\to
\Hcal^{\infty}(\kgot,M)$ is an isomorphism with inverse $p^*$.

\end{theo}

\begin{proof}
We prove the statement in the $\Ccal^{\infty}$ case.
 As $p\circ i=\Id_M$, we have $ i^*\circ p^*=\Id_{\Hcal^{\infty}(\kgot,M)}$.

Let us prove $ p^*\circ i^*=\Id_{\Hcal^{\infty}(\kgot,\Vcal)}.$
 We denote by $(m,v)$ a
point of $\Vcal$ with $m\in M$ and $v\in \Vcal_m$. For $t\geq 0$,
let $h(t)(m,v)=(m,t v)$  be the homothety on the fiber. The
transformations $h(t)$ verify $h(t_1)h(t_2)=h(t_1t_2)$ and commute
with the action of $K$.

Let $\alpha\in \Acal^{\infty}(\kgot, \Vcal)$ be a closed element
and let $\alpha(t)=h(t)^*\alpha$. Thus
$\alpha(0)=p^*\circ i^*(\alpha)$, while $\alpha(1)=\alpha$. From
Formula (\ref{Cartan}), we obtain for $t>0$:
\begin{equation}\label{homo}
\frac{d}{dt}\alpha(t)=\frac{1}{t}\Lcal(S)\cdot(\alpha(t))=
D\Big(\frac{1}{t}\iota(S)(\alpha(t))\Big).
\end{equation}
Here $S$ is the Euler vector field on
$\Vcal$ : at each point $(m,v)$ of $\Vcal$, $S_{(m,v)}=v$.

It is easily checked that $\frac{1}{t}(\iota(S)\alpha(t))$ is
continuous at $t=0$. Indeed, locally if  $\alpha=\sum_{I,J}
\nu_{I,J}(X,m,v) dm_I dv_J$, $\alpha(t)=\sum_{I,J}
\nu_{I,J}(X,m,tv) dm_I t^{|J]}dv_J$, and $\iota(S)$ kills all
components with $|J|=0$. Integrating Equation (\ref{homo}) from $0$ to $1$, we obtain:
$$\alpha-p^*\circ i^*(\alpha)=D\Big(\int_0^1\frac{1}{t}\iota(S)\alpha(t)dt\Big).$$
Thus we obtain the relation
$p^*\circ i^*=\Id_{\Hcal^{\infty}(\kgot,\Vcal)}.$

\end{proof}\bigskip

\subsection{The Chern-Weil construction}

Let $\pi:P\to B$ be a principal bundle with structure group $G$.
For any $G$-manifold $Z$, we define the manifold $\Zcal=P\times_G Z$ which
is fibred over $B$ with typical fiber $Z$. Let $\Acal(P\times Z)_{\hor}\subset \Acal(P\times Z)$
be the sub-algebra formed by the differential forms on $P\times Z$ which are
{\em horizontal}: $\gamma\in\Acal(P\times Z)_{\hor}$ if  $\iota(VX)\alpha=0$ on $P\times Z$
for every $X\in\ggot$. The algebra $\Acal(\Zcal)$ admits a natural identification with the
basic subalgebra
$$
\Acal(P\times Z)_{\bas}:=(\Acal(P\times Z)_{\hor})^G.
$$

Let $\omega\in (\Acal^{1}(P)\otimes \ggot)^G$ be a  connection one
form on $P$, with curvature form
$\Omega=d\omega+\frac{1}{2}[\omega,\omega]\in (\Acal^{2}(P)_{\hor}\otimes \ggot)^G$.
The connection one form $\omega$ defines, for any $G$-manifold $Z$, a projection from
$\Acal(P\times Z)^G$ onto $\Acal(P\times Z)_{\bas}$.

The Chern-Weil homomorphism
\begin{equation}\label{eq:chern-weil}
\phi_\omega^Z:\Acal^{\pol}(\ggot,Z) \longrightarrow \Acal(\Zcal).
\end{equation}
is defined as follows (see \cite{B-G-V}, \cite{Duflo-Vergne2}).
For a $G$-equivariant form $\alpha(X),X\in \ggot$ on $Z$, the
form $\phi_\omega^Z(\alpha)\in\Acal(\Zcal)$ is equal to the
projection of $\alpha(\Omega)\in\Acal(P\times Z)^G$ on the basic
subspace $\Acal(P\times Z)_{\bas}\simeq \Acal(\Zcal)$.

In the case where  $Z$ is the $\{{\rm pt}\}$, $\phi^Z_\omega$ is
the usual Chern-Weil homomorphism which associates to a
$G$-invariant polynomial $Q$ the characteristic form $Q(\Omega)$.

The main property of the equivariant cohomology differential $D$
proved by Cartan (see \cite{B-G-V}, \cite{gui-ste99}) is the
following proposition.

\begin{prop}
$$
\phi_\omega^Z \circ D=d \circ \phi_\omega^Z.
$$
Thus a closed  equivariant form on $Z$  gives rise to a closed  de Rham form on $\Zcal$.
\end{prop}

We can repeat the construction above in the equivariant case.

Let $K$ and $G$  be two compact Lie groups.  Assume that $P$ is
provided with an action of $K\times G$: $(k,g)(y)=ky g^{-1}$, for
$k\in K,g\in G,y\in P$. We assume that $G$ acts freely. Thus the
manifold $P/G=B$ is provided with a left action of $K$. Let
$\omega$ be a $K$-invariant connection one form on $P$, with
curvature form $\Omega$. For $Y\in\kgot$, we denote by
$\mu(Y)=-\iota(VY)\omega\, \in\f(P)\otimes\ggot$ the moment of
$Y$. The equivariant curvature form is
$$\Omega(Y)=\Omega+\mu(Y),\, X\in\kgot.$$

Let $Z$ be  a $G$-manifold. We consider the Chern-Weil homomorphism
\begin{equation}\label{eq:chern-weil-K}
\phi_\omega^Z:\Acal^{\pol}(\ggot,Z) \longrightarrow
\Acal^{\pol}(\kgot, \Zcal).
\end{equation}
It is defined as follows (see \cite{B-G-V},\cite{Duflo-Vergne2}). For a
$G$-equivariant form $\alpha$ on $Z$, the
value of the equivariant form $\phi_\omega^Z(\alpha)$ at $Y\in \kgot$ is
equal to the projection of $\alpha(\Omega(Y))\in\Acal(P\times Z)^G$
onto the the basic subspace $\Acal(P\times Z)_{\bas}\simeq \Acal(\Zcal)$.
For $Z=\{{\rm pt}\}$, and $Q$ a $G$-invariant polynomial on $\ggot$, the form
$\phi_\omega^Z(Q)(Y)=Q(\Omega(Y))$ is the Chern-Weil characteristic
class constructed in \cite {ber-ver82}, see also \cite{bott-tu00}.

\begin{prop}

The map $\phi_\omega^Z:\Acal^{\pol}(\ggot, Z)\to \Acal^{\pol}(\kgot, \Zcal)$
satisfies
$$
\phi_\omega^Z\circ D=D\circ\phi_\omega^Z.
$$

Here, on the left side the equivariant differential is with
respect to the action of $G$ while, on the right side, the
equivariant differential  is with respect to the group $K$.

\end{prop}

\section{Relative equivariant
cohomology}\label{section:coho-relative}

Let $N$ be a manifold provided with an action of a compact Lie group $K$.

\subsection{Definition and basic properties}

Let $F$ be a closed $K$-invariant subset of $N$.  To an
equivariant cohomology class on $N$ vanishing on $N\setminus F$,
we associate a relative equivariant cohomology class. Let us
explain the construction (see \cite{bott-tu00},\cite{pep-vergne1} for
the non-equivariant case).

 Consider the complex
$\Acal^{\pol}(\kgot,N,N\setminus F)$ with
$$
\Acal^{\pol}(\kgot,N,N\setminus F):=\Acal^{\pol}(\kgot,N)\oplus
\Acal^{\pol}(\kgot,N\setminus F)
$$
and differential $D_{\rm rel}\left(\alpha,\beta\right)=
\left(D\alpha,\alpha|_{N\setminus F}- D\beta \right)$.

\begin{defi}\label{relcoh}
The cohomology of the complex $(\Acal^{\pol}(\kgot,N,N\setminus
F),D_{\rm rel})$ is the relative equivariant cohomology space
$\Hcal^{\pol}(\kgot,N,N\setminus F)$.
\end{defi}

In the case where $K$ is the identity, we skip the letter $\kgot$
in the notation. Then $D_{\rel}$ is the usual relative de Rham
differential and  $\Hcal(N,N\setminus F)$ is the usual de Rham relative
cohomology group \cite{Bott-Tu}.

The complex $\Acal^{\pol}(\kgot,N,N\setminus F)$ is $\Zbb$-graded :
for $k\in\Zbb$, we take
$$
\left[\Acal^{\pol}(\kgot,N,N\setminus
F)\right]^k=\left[\Acal^{\pol}(\kgot,N)\right]^{k}\oplus
\left[\Acal^{\pol}(\kgot,N\setminus F)\right]^{k-1}.
$$

Since $D_{\rm rel}$ sends $\left[\Acal^{\pol}(\kgot,N,N\setminus
F)\right]^{k}$ into $\left[\Acal^{\pol}(\kgot,N,N\setminus
F)\right]^{k+1}$, the $\Zbb$-grading descends to the relative
cohomology spaces $\Hcal^{\pol}(\kgot,N,N\setminus F)$. The class defined by a
$D_{\rm rel}$-closed element $(\alpha,\beta)\in\Acal^{\pol}(\kgot,N,N\setminus F)$ will be
denoted $[\alpha,\beta]$.

Remark that $\Hcal^{\pol}(\kgot,N,N\setminus F)$ is a module over $\Hcal^{\pol}(\kgot,N)$.
Indeed the multiplication by a closed equivariant form $\eta\in\Acal^{\pol}(\kgot,N)$ ,
$$
\eta\cdot(\alpha,\beta)=(\eta\wedge \alpha,\eta\vert_{N\setminus F}\wedge\beta),
$$
on $\Acal^{\pol}(\kgot,N,N\setminus F)$ commutes with $D_{\rm
rel}$.

If $S$ is a $K$-invariant  vector field on $N$, we define on $\Acal^{\pol}(\kgot,N,N\setminus F)$
the operations $\Lcal(S)(\alpha,\beta):= (\Lcal(S)\alpha,\Lcal(S)\beta)$ and
$\iota(S)(\alpha,\beta):= (\iota(S)\alpha,-\iota(S)\beta)$. It is immediate to
check that Cartan' relation (\ref{Cartan}) holds
\begin{equation}\label{cartanrelative}
\Lcal(S)=\iota(S)\circ D_{\rel}+D_{\rel}\circ \iota(S).
\end{equation}

\bigskip

We consider now the following maps.
\begin{itemize}
\item The projection $j:\Acal^{\pol}(\kgot,N,N\setminus F)\to\Acal^{\pol}(\kgot,N)$ is the degree $0$
map defined by $j(\alpha,\beta)=\alpha$.

\item The inclusion  $i:\Acal^{ \pol}(\kgot,N\setminus F)\to\Acal^{\pol}(\kgot,N,N\setminus~F~)$
is the degree $+1$ map defined by $i(\beta)=(0,\beta)$.

\item The restriction $r:\Acal^{\pol}(\kgot,N)\to\Acal^{\pol}(\kgot,N\setminus~ F)$
is the degree $0$ map defined by $r(\alpha)=\alpha|_{N\setminus~F}$.
\end{itemize}

It is easy to see that $i,j,r$ induce maps in cohomology that we still denote
by $i,j,r$.

\begin{prop}
$\bullet$  We have an exact triangle
$$
\xymatrix{
    & \Hcal^{\pol}(\kgot,N,N\setminus~F~)\ar@{<-}_{i}[ld]\ar@{->}^{j}[rd] &     \\
\Hcal^{\pol}(\kgot,N\setminus F)   &                &  \Hcal^{\pol}(\kgot,N).\ar@{->}^{r}[ll]
}
$$

$\bullet$ If $F\subset F'$ are closed $K$-invariant subsets of
$N$, the restriction map $(\alpha,\beta)\mapsto
(\alpha,\beta|_{N\setminus F'})$ induces a map
\begin{equation}\label{eq:res-F}
\res_{F',F}:\Hcal^{\pol}(\kgot,N,N\setminus F)\to
\Hcal^{\pol}(\kgot,N,N\setminus F').
\end{equation}

$\bullet$ If $g$ is a diffeomorphism of $N$ which preserves $F$
and commutes with the action of $K$, then
$g^*(\alpha,\beta)=(g^*\alpha,g^*\beta)$ induces a transformation
$g^*$ of $\Hcal^{ \pol}(\kgot,N,N\setminus F)$.
\end{prop}

\begin{proof}
This proof is the same than in the non equivariant-case
\cite{Bott-Tu} and we skip it.
\end{proof}\bigskip

\bigskip

The same statements hold in the $\Ccal^{\infty}$-case. Here we
consider the complex $\Acal^{\infty}(\kgot,N,N\setminus F)$ with
$\Acal^{\infty}(\kgot,N,N\setminus F)
:=\Acal^{\infty}(\kgot,N)\oplus \Acal^{\infty}(\kgot,N\setminus F)$
and differential $D_{\rm rel}\left(\alpha,\beta\right)=
\left(D\alpha,\alpha|_{N\setminus F}- D\beta \right)$.

\begin{defi}\label{relcohinfty}
The cohomology of the complex $(\Acal^{\infty}(\kgot,N,N\setminus
F),D_{\rm rel})$ is the relative equivariant cohomology spaces
$\Hcal^{\infty}(\kgot,N,N\setminus F)$.
\end{defi}

The complex $\Acal^{\infty}(\kgot,N,N\setminus F)$ is
$\Zbb_2$-graded by taking
$\left[\Acal^{\infty}(\kgot,N,N\setminus
F)\right]^\epsilon=$ $\left[\Acal^{\infty}(\kgot,N)\right]^{\epsilon}\oplus
\left[\Acal^{\infty}(\kgot,N\setminus F)\right]^{\epsilon+1}$.
Since $D_{\rm rel}$ send $\left[\Acal^{\infty}(\kgot,N,N\setminus
F)\right]^{\epsilon}$ into \break $\left[\Acal^{\infty}(\kgot,N,N\setminus F)\right]^{\epsilon+1}$,
the $\Zbb_2$-grading descends to the relative cohomology spaces
$\Hcal^{\infty}(\kgot,N,N\setminus F)$.

Here the space $\Hcal^{\infty}(\kgot,N,N\setminus F)$ is a module over
$\Hcal^{\infty}(\kgot,N)$.

\medskip

\begin{lem}\label{lem:basic-relative-coho}
$\bullet$ We have an exact triangle
$$
\xymatrix{
    & \Hcal^{\infty}(\kgot,N,N\setminus~F~)\ar@{<-}_{i}[ld]\ar@{->}^{j}[rd] &     \\
\Hcal^{\infty}(\kgot,N\setminus F)   &                &  \Hcal^{\infty}(\kgot,N).\ar@{->}^{r}[ll]
}
$$

$\bullet$ If $F\subset F'$ are closed $K$-invariant subsets of
$N$, the restriction $(\alpha,\beta)\mapsto
(\alpha,\beta|_{N\setminus F'})$ induces a map
$\res_{F',F}:\Hcal^{\infty}(\kgot,N,N\setminus F)\to
\Hcal^{\infty}(\kgot,N,N\setminus F')$.

$\bullet$ If $g$ is a diffeomorphism of $N$ which preserves $F$
and commutes with the action of $K$, then
$g^*(\alpha,\beta)=(g^*\alpha,g^*\beta)$ induces a transformation
$g^*$ of $\Hcal^{ \infty}(\kgot,N,N\setminus F)$.
\end{lem}

\subsection{Excision}\label{subsec:excision}

Let $\chi\in \f(N)^K$ be a $K$-invariant function such that $\chi$
is identically equal to $1$ on a neighborhood of $F$. If $\beta\in
\Acal^{\pol}(\kgot, N\setminus F)$, note that $d\chi\wedge \beta$
defines an equivariant form on $N$, since $d\chi$ is equal to $0$
in a neighborhood of $F$. We define
\begin{equation}\label{eq:I-chi}
I^{\chi}:\Acal^{\pol}(\kgot, N,N\setminus F) \to
\Acal^{\pol}(\kgot, N,N\setminus F)
\end{equation}
by $I^{\chi}(\alpha,\beta)=(\chi \alpha+d\chi \wedge \beta,\chi\beta)$. Then
\begin{equation}\label{Ichi}
I^{\chi}\circ D_{\rel}=D_{\rel}\circ I^{\chi},
\end{equation}
 so that $I^{\chi}$ defines a map $I^{\chi}:
\Hcal^{\pol}(\kgot, N,N\setminus F)\to \Hcal^{\pol}(\kgot,
N,N\setminus F)$.

\begin{lemm}\label{lem:excision}
The map $I^{\chi}$ is independent of $\chi$. In particular,
$I^{\chi}$ is  the identity in relative cohomology.
\end{lemm}

\begin{proof}
If $(\alpha,\beta)$ is $D_{\rel}$-closed, then
 $I^{\chi_1}(\alpha,\beta)-I^{\chi_2}(\alpha,\beta)=$ \break
 $D_{\rel}((\chi_1-\chi_2)\beta,0)$.
This shows that $I^{\chi}$ is independent of $\chi$. Choosing
$\chi=1$, we see that $I^{\chi}=\Id$ in cohomology.

\end{proof}

\bigskip

It follows from the above proposition that we can always choose a
representative  $(\alpha,\beta)$ of a relative cohomology class,
with $\alpha$ and $\beta$ supported in a neighborhood of $F$
as small as we want. This will be important to define the integral
over the fiber of a relative cohomology class with support
intersecting the fibers in compact subsets. The integration will
be defined in Section \ref{sec:integral-fiber}

In particular, if $F$ is compact, we define a map
\begin{equation}\label{eq:p-compact}
\p_c:\Hcal^{\pol}(\kgot, N,N\setminus F)\to \Hcal^{\pol}_c(\kgot, N).
\end{equation}
by setting $\p_c(\alpha,\beta)=\chi\alpha+d\chi\wedge \beta$,
where  $\chi\in \f(N)^K$ is a $K$-invariant function  with compact
support such that $\chi$ is identically equal to $1$ on a
neighborhood of $F$.

\medskip

An important property  of the relative cohomology group is the
{\em excision property}. Let $U$ be a $K$-invariant neighborhood of $F$. The
restriction $(\alpha,\beta)\mapsto (\alpha|_U,\beta|_{U\setminus F})$ induces a map
$$
\res^U:\Hcal^{\pol}(\kgot, N,N\setminus F)\to \Hcal^{\pol}(\kgot, U,U\setminus F).
$$

\begin{prop}\label{prop:excision}
The map $\res^U$ is an isomorphism.
\end{prop}
\begin{proof}
Let us choose $\chi\in \f(N)^K $ supported in $U$ and equal to $1$ in a neighborhood of
$F$. The map (\ref{eq:I-chi}) defines in this context three maps :
$I^{\chi}_N:\Acal^{\pol}(\kgot, N,N\setminus F) \to
\Acal^{\pol}(\kgot, N,N\setminus F)$, $I^{\chi}_U:\Acal^{\pol}(\kgot, U,U\setminus F) \to
\Acal^{\pol}(\kgot, U,U\setminus F)$ and
$I^{\chi}_{N,U}:\Acal^{\pol}(\kgot, U,U\setminus F) \to
\Acal^{\pol}(\kgot, N,N\setminus F)$.

We check easily that $I^{\chi}_{N,U}\circ \,\res^U=I^{\chi}_N$ on $\Acal^{\pol}(\kgot, N,N\setminus F)$, and
that $\res^U\circ \, I^{\chi}_{N,U}=I^{\chi}_U$ on  $\Acal^{\pol}(\kgot, U,U\setminus F)$.
From Lemma \ref{lem:excision}, we know  that $I^{\chi}_N$ and $I^{\chi}_U$ are the identity maps
in cohomology. This proves that $\res^U$ is an isomorphism in cohomology.



\end{proof}\bigskip

The same statements holds in the $\Ccal^{\infty}$-case.

\begin{prop}
$\bullet$ Let $U$ be a $K$-invariant neighborhood of $F$. The map
$\res^U: \Hcal^{\infty}(\kgot, N,N\setminus F)\to \Hcal^{\infty}(\kgot,U,U\setminus F)$
is an isomorphism.

$\bullet$ If $F$ is compact, there is a natural map $\p_c:
\Hcal^{\infty}(\kgot, N,N\setminus F)\to \Hcal^{\infty}_c(\kgot,N)$.
\end{prop}

\subsection{Product in relative equivariant cohomology}\label{section:coho-relproduct}

Let $F_1$ and $F_2$ be two closed $K$-invariant subsets of $N$. We
will now define  a graded product
\begin{eqnarray}\label{eq:produit-relatif}
\Hcal^{\pol}(\kgot,N,N\setminus F_1)
\times\Hcal^{\pol}(\kgot,N,N\setminus F_2) &\longrightarrow&
\Hcal^{\pol}(\kgot,N,N\setminus (F_1\cap F_2)) \nonumber\\
\quad\quad\quad\quad(\quad a \quad,\quad b\quad
)\quad\quad\quad\quad &\longmapsto& a\diamond b \ .
\end{eqnarray}

\medskip

Let $U_1:=N\setminus F_1$, $U_2:=N\setminus F_2$ so that
$U:=N\setminus (F_1\cap F_2)=U_1\cup U_2$. Let
$\Phi:=(\Phi_1,\Phi_2)$ be a partition of unity subordinate to the
covering $U_1\cup U_2$ of $U$. By averaging by $K$, we may suppose
that the functions $\Phi_k$ are invariant.

Since $\Phi_i\in\f(U)^K$ is supported in $U_i$, the product
$\gamma\mapsto \Phi_i\gamma$ defines a map
$\Acal^{\pol}(\kgot,N\setminus F_i) \to
\Acal^{\pol}(\kgot,N\setminus (F_1\cap F_2))$. Since
$d\Phi_1=-d\Phi_2\in \Acal(U)^K$ is supported in $U_1\cap U_2=
N\setminus (F_1\cup F_2)$ , the product $\gamma\mapsto
d\Phi_1\wedge\gamma$ defines a map \break
$\Acal^{\pol}(\kgot,N\setminus (F_1\cup F_2)) \to
\Acal^{\pol}(\kgot,N\setminus (F_1\cap F_2))$.

With the help of $\Phi$, we define a bilinear map
$\diamond_\Phi:\Acal^{\pol}(\kgot,N,N\setminus F_1)\times
\Acal^{\pol}(\kgot,N,N\setminus F_2) \to
\Acal^{\pol}(\kgot,N,N\setminus (F_1\cap F_2))$ as follows. For
$a_i:=(\alpha_i,\beta_i)\in \Acal^{{\pol}}(\kgot,N,N\setminus
F_i)$, $i=1,2$, we define
$$
a_1 \diamond_\Phi a_2 :=\Big(\alpha_1\wedge
\alpha_2,\beta_\Phi(a_1,a_2)\Big)
$$
with
$$
\beta_\Phi(a_1,a_2)= \Phi_1\beta_1 \wedge
\alpha_2+(-1)^{|a_1|}\alpha_1\wedge
\Phi_2\beta_2-(-1)^{|a_1|}d\Phi_1\wedge \beta_1\wedge \beta_2.
$$

Remark that all equivariant forms $\Phi_1\beta_1 \wedge \alpha_2$,
$\alpha_1\wedge \Phi_2\beta_2$ and $d\Phi_1\wedge \beta_1\wedge
\beta_2$ are well defined on $U_1\cup U_2$. So $a_1 \diamond_\Phi
a_2\in \Acal^{{\pol}}(\kgot,N,N\setminus (F_1\cap F_2))$.
It is immediate to verify that
$$
D_{\rm rel}(a_1 \diamond_\Phi a_2)=
(D_{\rm rel}a_1) \diamond_\Phi a_2+(-1)^{|a_1|}a_1 \diamond_\Phi (D_{\rm rel}a_2).
$$
Thus $\diamond_\Phi$ defines a bilinear map
$\Hcal^{\pol}(\kgot,N,N\setminus F_1)\times
\Hcal^{\pol}(\kgot,N,N\setminus F_2) \to\Hcal^{\pol}(\kgot,N,N\setminus (F_1\cap F_2))$.
Let us see that this product do not depend of the choice of the
partition of unity. If we have another partition
$\Phi'=(\Phi_1',\Phi_2')$, then $\Phi_1-\Phi'_1=-(\Phi_2-\Phi'_2)$.
It is immediate to verify that, if $D_{\rm rel}(a_1)=0$ and $D_{\rm
rel}(a_2)=0$, one has
$$
a_1 \diamond_\Phi a_2 - a_1 \diamond_{\Phi'}a_2=D_{\rm
rel}\Big(0,(-1)^{|a_1|} (\Phi_1-\Phi'_1)\beta_1\wedge \beta_2\Big).
$$
in $\Acal^{\pol}(\kgot,N,N\setminus (F_1\cap F_2))$. So the product on the relative
cohomology spaces will be denoted by $\diamond$.

\begin{rem}
The same formulae defines a $\Zbb_2$-graded product
\begin{eqnarray}\label{eq:produit-relatif-infty}
\Hcal^{\infty}(\kgot,N,N\setminus F_1) \times
\Hcal^{\infty}(\kgot,N,N\setminus F_2) &\longrightarrow&
\Hcal^{\infty}(\kgot,N,N\setminus (F_1\cap F_2)) \nonumber\\
\quad\quad\quad\quad(\quad a \quad,\quad b\quad
)\quad\quad\quad\quad &\longmapsto& a\diamond b \ .
\end{eqnarray}
\end{rem}

We note the following properties, which are well known in the
non-equivariant case.

\begin{prop}

\begin{itemize}

\item The relative product is compatible with restrictions:
if $F_1\subset F'_1$ and $F_2\subset F_2'$ are closed invariant
subsets of $N$, then the diagram
\begin{equation}\label{eq:fonctoriel-prod-relatif}
\xymatrix@C=6mm{ \Hcal^{\pol}(\kgot,N,N\setminus F_1)\ar[d]^{\res_1}
\!\!\!\!\! & \times \quad  \Hcal^{\pol}(\kgot,N,N\setminus F_2)\ar[d]^{\res_2}\ar[r]^-\lozenge &
\   \Hcal^{\pol}(\kgot,N,N\setminus (F_1\cap F_2)) \ar[d]^{\res_{12}}\\
\Hcal^{\pol}(\kgot,N,N\setminus F_1') \!\!\!\!\! & \times \quad
\Hcal^{\pol}(\kgot,N,N\setminus F_2') \ar[r]^-\lozenge & \
\Hcal^{\pol}(\kgot,N,N\setminus (F_1'\cap F_2'))}
\end{equation}

is commutative. Here the $\res_i$ are the restrictions maps defined in (\ref{eq:res-F}).\\

\item The relative product is graded commutative :
$a_1\diamond a_2 = (-1)^{|a_1|.|a_2|}a_2\diamond a_1$.\\

\item

The relative product is associative.\\

\end{itemize}

The same statements holds in the $\Ccal^{\infty}$-case.

\end{prop}

\begin{proof} The first point follows from the definition. Let $U_i=N\setminus F_i$,
$U'_i=N\setminus F'_i$. Let $\Phi_1+\Phi_2=1$ be a partition of
unity on $U_1\cup U_2$. Then $\Phi'=(\Phi'_1,\Phi'_2)$ with $\Phi'_i:=\Phi_i|_{U_i'}$
is a partition of unity on  $U'_1\cup U'_2$. Then, at the level of equivariant forms,
we have $\beta_{\Phi}(a_1,a_2)|_{N\setminus(F'_1\cap F'_2)}=
\beta_{\Phi'}(\res_1(a_1),\res_2(a_2))$. The commutative
diagram (\ref{eq:fonctoriel-prod-relatif}) follows.


The second point is immediate from the definition.

We now prove the third point. Let $F_1,F_2$ and $F_3$ be three
closed invariant subsets of $N$. Let
$a_i\in\Hcal^{\pol}(\kgot,N,N\setminus F_i)$ for $i=1,2,3$. In
order to prove that $(a_1\diamond a_2)\diamond a_3=a_1\diamond
(a_2\diamond a_3)$ in $\Hcal^{\pol}(\kgot,N,N\setminus (F_1\cap
F_2\cap F_3))$ we introduce a multi-linear map $\top: E_1\times
E_2\times E_3\longrightarrow \Hcal^{\pol}(\kgot,N,N\setminus
(F_1\cap F_2\cap F_3))$ where $E_i=\Hcal^{\pol}(\kgot,N,N\setminus
F_i)$.

Let $U_i=N\setminus F_i$ and $U=N\setminus
(F_1\cap F_2\cap F_3)$. Let $\Dcal$ be the data formed by :
\begin{itemize}
\item A partition of unity $\Phi_1+\Phi_2 +\Phi_3=1$ on $U_1\cup
U_2\cup U_3=U$, where the functions $\Phi_i$ are $K$-invariant.

\item Invariant one forms $\Lambda_1,\Lambda_2$ and $\Lambda_3$ on
$U$ supported respectively in $U_2\cap U_3$, $U_1\cap U_3$ and
$U_1\cap U_2$.
\end{itemize}
We suppose that the data $\Dcal$ satisfies the following
conditions
\begin{equation}\label{eq:condition1-prod}
d\Phi_1=\Lambda_2-\Lambda_3,\ d\Phi_2=\Lambda_3-\Lambda_1,\
d\Phi_3=\Lambda_1-\Lambda_2.
\end{equation}
Then we have
\begin{equation}\label{eq:condition2-prod}
D\Lambda_1(X)=D\Lambda_2(X)=D\Lambda_3(X).
\end{equation}
We denote $\Theta(X)$ the equivariant $2$-form equal to
$D\Lambda_1(X)$ : (\ref{eq:condition2-prod}) shows that
$\Theta(X)$ is supported in $U_1\cap U_2\cap U_3$.

With the help of $\Dcal$, we define a three-linear map
$\top_\Dcal$ from $\Acal^{\pol}(\kgot,N,N\setminus F_1)\times
\Acal^{\pol}(\kgot,N,N\setminus F_2) \times
\Acal^{\pol}(\kgot,N,N\setminus F_3)$ into
$\Acal^{\pol}(\kgot,N,N\setminus (F_1\cap F_2\cap F_3))$ as
follows. For $a_i:=(\alpha_i,\beta_i)\in
\Acal^{{\pol}}(\kgot,N,N\setminus F_i)$, $i=1,2,3$, we define
$\top_\Dcal(a_1,a_2,a_3) :=(\alpha_1\wedge
\alpha_2\wedge\alpha_3,\beta_\Dcal(a_1,a_2,a_3))$ with
\begin{eqnarray}\label{eq:beta123}
&&\beta_\Dcal(a_1,a_2,a_3)= \Phi_1\beta_1\alpha_2\alpha_3
+(-1)^{|a_1|}\Phi_2\alpha_1\beta_2 \alpha_3
+(-1)^{|a_1|+|a_2|}\Phi_3\alpha_1\alpha_2\beta_3\\
&&\qquad +(-1)^{|a_2|}\Lambda_1\alpha_1\beta_2 \beta_3
-(-1)^{|a_1|+|a_2|}\Lambda_2\beta_1\alpha_2\beta_3
+(-1)^{|a_1|}\Lambda_3\beta_1\beta_2 \alpha_3\nonumber\\
&&\qquad-(-1)^{|a_2|}\Theta\beta_1 \beta_2\beta_3.\nonumber
\end{eqnarray}
Remark that all equivariant forms which appears in the right hand
side of (\ref{eq:beta123}) are well defined on $U_1\cup U_2\cup
U_3$. So $\top_\Dcal(a_1,a_2,a_3)\in
\Acal^{{\pol}}(\kgot,N,N\setminus (F_1\cap F_2\cap F_3)).$

The following relation is ``immediate'' to verify:

$$ D_{\rm rel}(\top_\Dcal(a_1,a_2,a_3))=$$
$$ \top_\Dcal(D_{\rm rel}a_1,a_2,a_3)+(-1)^{|a_1|}
\top_\Dcal(a_1,D_{\rm rel}a_2,a_3)+
(-1)^{|a_1|+|a_2|}\top_\Dcal(a_1,a_2,D_{\rm rel}a_3).
$$
Thus $\top_\Dcal$ defines a three-linear map from
from $E_1\times E_2\times E_3$ into \break
$\Hcal^{\pol}(\kgot,N,N\setminus (F_1\cap F_2\cap F_3))$. Let us
see that this map do not depend of the choice of the data $\Dcal$.
Let $\Dcal'=\{\Phi'_i,\Lambda'_i\ {\rm for}\ i=1,2,3\}$ be another
data which satisfies conditions (\ref{eq:condition1-prod}).

We consider the functions $f_i=\Phi_i-\Phi'_i$ on $U$. If
$\{i,j,k\}=\{1,2,3\}$ the function $f_i$ is supported in $U_i\cap
(U_j\cup U_k)=(U_i\cap U_j)\cup (U_i\cap U_k)$. The relations
$f_1+f_2+f_3=0$ on $U$ shows that there exists $K$-invariant
functions $\theta_i$ on $U$ such that $\theta_i$ is supported in
$U_j\cap U_k$ and
$$
f_1=\theta_2-\theta_3,\ f_2=\theta_3-\theta_1,\
f_3=\theta_1-\theta_2.
$$
We see then that
$$
\Lambda_1-\Lambda'_1-d\theta_1=\Lambda_2-\Lambda'_2-d\theta_2=
\Lambda_3-\Lambda'_3-d\theta_3
$$
is an invariant one form on $U$ supported on $U_1\cap U_2\cap
U_3$: let us denote it by $\Delta$. We have
$\Theta(X)-\Theta'(X)=D\Delta(X)$.

Then for $D_{\rm rel}$-closed elements $a_i$, one checks that
$\top_\Dcal(a_1,a_2,a_3)-\top_{\Dcal'}(a_1,a_2,a_3)$ is equal to
$D_{\rm rel}(0,-\delta)=(0,D\delta)$ with
\begin{eqnarray*}
&&\delta=(-1)^{|a_2|}\theta_1\alpha_1\beta_2\beta_3
-(-1)^{|a_1|+|a_2|}\theta_2\beta_1\alpha_2 \beta_3
+(-1)^{|a_1|}\theta_3\beta_1\beta_2\alpha_3\\
&&\qquad -(-1)^{|a_2|}\Delta\beta_1\beta_2 \beta_3.
\end{eqnarray*}
Let us denote $\top$ the three-linear map induced by $\top_\Dcal$
in relative equivariant cohomology.

\medskip

Now we will see that the map $(a_1,a_2,a_3)\mapsto(a_1\diamond
a_2)\diamond a_3$ coincides with $\top$.
 Let $\phi_1+\phi_2=1$ be a partition of unity on $U_{12}:=U_1\cup U_2$, and
let $\varphi_{12} +\varphi_3=1$ be a partition of unity on
$U_{12}\cup U_3$ : all the functions are supposed $K$-invariant.
Then we take the data $\Dcal=\{\Phi_i,\Lambda_i\ {\rm for}\
i=1,2,3\}$ defined by the relations:
\begin{itemize}
\item $\Phi_1=\varphi_{12}\phi_1$,   $\Phi_2=\varphi_{12}\phi_2$,
$\Phi_3=\varphi_3$,

\item $\Lambda_1=-d(\varphi_{12})\phi_2$,
$\Lambda_2=d(\varphi_{12})\phi_1$,
$\Lambda_3=-\varphi_{12}d\phi_1$.
\end{itemize}
One checks that $\Dcal$ satisfies conditions
(\ref{eq:condition1-prod}), and that for
$a_i\in\Acal^{{\pol}}(\kgot,N,N\setminus~F_i~)$ the following
equality
$$
\top_\Dcal(a_1,a_2, a_3)=(a_1\diamond_\phi a_2)\diamond_\varphi
a_3
$$
holds in $\Acal^{{\pol}}(\kgot,N,N\setminus (F_1\cap F_2\cap F_3))$.

One proves in the same way that the map $(a_1,a_2,a_3)\mapsto
a_1\diamond (a_2\diamond a_3)$ coincides with $\top$. We have then
proved the associativity of the relative product $\diamond$.

\end{proof}

\subsection{Inverse limit of equivariant cohomology with support}
\label{sec:cohomologie-support}

Let $F$ be a closed $K$-invariant subset of $N$.
 We consider the
set $\Fcal_F$ of all open invariant neighborhoods $U$ of $F$ which
is ordered by the relation $U\leq V$ if and only if $V\subset U$.
For any $U\in \Fcal_F$, we consider the algebra
$\Acal_U^{\pol}(\kgot,N)$ of equivariant differential forms on $N$
with support contained in $U$: $\alpha\in \Acal_U^{\pol}(\kgot,N)$
if there exists a closed set $\Ccal_\alpha\subset U$ such that
$\alpha(X)|_n=0$ for all $X\in\kgot$ and all $n\in U\setminus
\Ccal_\alpha$. Note that the vector space $\Acal_U^{\pol}(\kgot,N)$ is naturally
a module over $\Acal^{\pol}(\kgot,N)$.

The algebra $\Acal_U^{\pol}(\kgot,N)$ is stable under the
differential $D$, and we denote by $\Hcal_U^{\pol}(\kgot,N)$ the
corresponding cohomology algebra. Note that $\Hcal_U^{\pol}(\kgot,N)$ is naturally
a module over $\Hcal^{\pol}(\kgot,N)$. If $U\leq V$, we have then an
inclusion map $\Acal_V^{\pol}(\kgot,N)\hookrightarrow
\Acal_U^{\pol}(\kgot,N)$ which gives rise to a map
$r_{U,V}:\Hcal_V^{\pol}(\kgot,N)\to \Hcal_U^{\pol}(\kgot,N)$ of
$\Hcal^{\pol}(\kgot,N)$-modules.

\begin{defi}\label{defiinductive}
 We denote by   $\Hcal_F^{\pol}(\kgot,N)$ the inverse limit of the inverse system
 $(\Hcal_{U}^{\pol}(\kgot,N),r_{U,V};U,V\in\Fcal_F)$. We will call
$\Hcal_F^{\pol}(\kgot,N)$ the equivariant cohomology of $N$
supported on $F$.
\end{defi}

Note that the vector space $\Hcal_F^{\pol}(\kgot,N)$ is naturally
a module over $\Hcal^{\pol}(\kgot,N)$. Let us give the following basic properties of the equivariant
cohomology spaces with support.

\begin{lem}\label{lem:basic-equi-coho-support-pol}

$\bullet$ $\Hcal_F^{ {\pol}}(\kgot,N)=\{0\}$ if $F=\emptyset$.

$\bullet$ There is a natural map $\Hcal^{{\pol}}_F(\kgot,N)\to
\Hcal^{{\pol}}(\kgot, N)$. If $F$ is compact, this map factors
through $\Hcal^{{\pol}}_F(\kgot,N)\to \Hcal_c^{{\pol}}(\kgot,N)$.

$\bullet$ If $F\subset F'$ are closed $K$-invariant subsets, there
is a restriction morphism
\begin{equation}\label{eq:res-F-local}
\res^{F',F}:\Hcal^{{\pol}}_F(\kgot,N)\to
\Hcal^{{\pol}}_{F'}(\kgot,N).
\end{equation}

$\bullet$ If $F_1$ and $F_2$ are two closed $K$-invariant subsets of
$N$, the wedge product of forms defines a natural product
\begin{equation}\label{eq:produit-equi-gene}
    \Hcal^{{\pol}}_{F_1}(\kgot,N)\times \Hcal_{F_2}^{ {\pol}}(\kgot,N)
    \stackrel{\wedge}{\longrightarrow}
    \Hcal_{F_1\cap F_2}^{ {\pol}}(\kgot,N).
\end{equation}

$\bullet$ If $F_1\subset F_1'$ and $F_2\subset F_2'$  are closed
$K$-invariant subsets, then the diagram
\begin{equation}\label{eq.2.reduction-gene}
\xymatrix@C=1cm{ \Hcal^{{\pol}}_{F_1}(\kgot,N)\ar[d]^{\res^1} \!&\!\times
\!&\! \Hcal^{ {\pol}}_{F_2}(\kgot,N)\ar[d]^{\res^2}\ar[r]^-{\wedge} &
\quad \Hcal_{F_1\cap F_2}^{\pol}(\kgot,N) \ar[d]^{\res^{12}}\\
\Hcal_{F_1'}^{{\pol}}(\kgot,N)   \!&\!\times \! &\!
\Hcal_{F_2'}^{\pol}(\kgot,N) \ar[r]^-{\wedge}   &
\quad \Hcal_{F_1'\cap F_2'}^{\pol}(\kgot,N)
  }
\end{equation}
is commutative. Here the $\res^i$ are the restriction morphisms defined in
(\ref{eq:res-F-local}).
\end{lem}

All the maps in the previous lemma preserve the structures of
$\Hcal^{\pol}(\kgot,N)$-module. The same definition and properties hold in the
$\Ccal^{\infty}$-case.
\begin{defi}\label{defiinductive-infty}
 We denote by   $\Hcal_F^{\infty}(\kgot,N)$ the inverse limit of the inverse system
 $(\Hcal_{U}^{\infty}(\kgot,N),r_{U,V};U,V\in\Fcal_F)$. We will call
$\Hcal_F^{\infty}(\kgot,N)$ the equivariant cohomology of $N$
supported on $F$.
\end{defi}

\begin{lem}\label{lem:basic-equi-coho-support}

$\bullet$ $\Hcal_F^{ {\infty}}(\kgot,N)=\{0\}$ if $F=\emptyset$.

$\bullet$ There is a natural map $\Hcal^{{\infty}}_F(\kgot,N)\to
\Hcal^{{\infty}}(\kgot, N)$. If $F$ is compact, this map factors
through $\Hcal^{{\infty}}_F(\kgot,N)\to
\Hcal_c^{{\infty}}(\kgot,N)$.

$\bullet$ If $F\subset F'$ are closed $K$-invariant subsets, there
is a restriction morphism  $\res^{F',F}:\Hcal^{{\infty}}_F(\kgot,N)\to
\Hcal^{{\infty}}_{F'}(\kgot,N)$.

$\bullet$ If $F$ and $R$ are two closed $K$-invariant subsets of
$N$, the wedge product of forms defines a natural product
\begin{equation}\label{eq:produit-equi-gene-infty}
    \Hcal^{{\infty}}_{F_1}(\kgot,N)\times \Hcal_{F_2}^{ {\infty}}(\kgot,N)
    \stackrel{\wedge}{\longrightarrow}
    \Hcal_{F_1\cap F_2}^{ {\infty}}(\kgot,N).
\end{equation}

$\bullet$ If $F_1\subset F_1'$ and $F_2\subset F_2'$  are closed
$K$-invariant subsets, then the diagram
\begin{equation}
\xymatrix@C=1cm{ \Hcal^{\infty}_{F_1}(\kgot,N)\ar[d]^{\res^1} \!&\!\times
\!&\! \Hcal^{ \infty}_{F_2}(\kgot,N)\ar[d]^{\res^2}\ar[r]^-{\wedge} &
\quad \Hcal_{F_1\cap F_2}^\infty(\kgot,N) \ar[d]^{\res^{12}}\\
\Hcal_{F_1'}^{\infty}(\kgot,N)   \!&\!\times \! &\!
\Hcal_{F_2'}^\infty(\kgot,N) \ar[r]^-{\wedge}   &
\quad \Hcal_{F_1'\cap F_2'}^\infty(\kgot,N)
  }
\end{equation}
is commutative.
\end{lem}

\subsection{Morphism $\p_F$}\label{section:morphism-pF}

If $F$ is any closed invariant subset of $N$, we define a morphism
\begin{equation}\label{eq:p-F}
\p_F:\Hcal^{\pol}(\kgot, N,N\setminus F)\longrightarrow
\Hcal^{\pol}_F(\kgot, N)
\end{equation}
of $\Hcal^{\pol}(\kgot, N)$-modules.

\begin{prop}\label{alphau}
For any open invariant neighborhood $U$ of $F$, we choose
$\chi\in\f(N)^K$ with support in $U$ and equal to $1$ in a
neighborhood of $F$.

$\bullet$ The map
\begin{equation}\label{eq:p-U-chi}
\p^{\chi}\left(\alpha,\beta\right)=\chi\alpha + d\chi\beta
\end{equation}
defines a morphism $\p^{\chi}:\Acal^{{\pol}}(\kgot,N,N\setminus F)
\to \Acal_U^{ {\pol}}(\kgot,N)$ of $\Acal^{\pol}(\kgot, N)$-modules.

In consequence, let $\alpha\in \Acal^{ {\pol}}(\kgot,N)$ be a
\emph{closed} equivariant form and $\beta\in \Acal^{
{\pol}}(\kgot,N\setminus F)$ such that $\alpha|_{N\setminus
F}=D\beta$, then $\p^{\chi}(\alpha,\beta)$ is a \emph{closed}
equivariant form supported in $U$.

$\bullet$ The  cohomology class of  $\p^{\chi}(\alpha,\beta)$ in
$\Hcal_U^{ {\pol}}(\kgot,N)$ does not depend of $\chi$. We denote
this class by $\p_U(\alpha,\beta)\in \Hcal_U^{ {\pol}}(\kgot,N)$.

$\bullet$ For any neighborhoods $V\subset U$ of $F$, we have
$r_{U,V}\circ \p_V=\p_U$.
\end{prop}
\begin{proof}
The proof is similar to the proof of Proposition 2.3 in
\cite{pep-vergne1}. We repeat the main arguments.

The equation $\p^{\chi}\circ \,D_{\rm rel}= D \circ \p^{\chi}$
is immediate to check. In particular $\p^{\chi}(\alpha,\beta)$ is
closed, if $D_{\rm rel}\left(\alpha,\beta\right)=0$. Directly:
$D(\chi\alpha+ d\chi \beta)=d\chi \alpha- d\chi D\beta=0$. This
proves the first point. For two different choices $\chi$ and
$\chi'$, we have
\begin{eqnarray*}
\p^{\chi}(\alpha,\beta)-\p^{\chi'}(\alpha,\beta)
&=& (\chi-\chi')\alpha + d(\chi-\chi')\beta \\
&=& D\left((\chi-\chi')\beta\right).
\end{eqnarray*}
Since $\chi-\chi'=0$ in a neighborhood of $F$, the right hand side
of the last equation is well defined, and is an element of
$\Acal_U^{ {\pol}}(\kgot,N)$. This proves the second  point.
Finally, the last point is immediate, since
$\p_U(\alpha,\beta)=\p_V(\alpha,\beta)=\p^\chi(\alpha,\beta)$ for
$\chi\in\f(N)^K$ with support in $V\subset U$.
\end{proof}

\bigskip

\begin{defi}\label{def-alpha-beta}
Let $\alpha\in \Acal^{ {\pol}}(\kgot,N)$ be a \emph{closed}
equivariant form and \break $\beta\in \Acal^{
{\pol}}(\kgot,N\setminus F)$ such that $\alpha|_{N\setminus
F}=D\beta$. We denote by $\p_F(\alpha,\beta)\in \Hcal_F^{
{\pol}}(\kgot,N)$ the element defined by the sequence
$\p_U(\alpha,\beta)\in \Hcal_U^{ {\pol}}(\kgot,N),\, U\in\Fcal_F$.
We have then a morphism
$\p_F:\Hcal^{ {\pol}}(\kgot,N,N\setminus F)\to \Hcal_F^{
{\pol}}(\kgot,N)$.
\end{defi}

The following proposition summarizes the functorial properties of
$\p$.

\begin{prop}

$\bullet$ If $F\subset F'$ are closed invariant subsets of $N$, then
the diagram
\begin{equation}\label{eq:fonctoriel-p-restriction}
\xymatrix@C=1cm{ \Hcal^{ {\pol}}(\kgot,N,N\setminus F)
\ar[d]^{\res_1}\ar[r]^-{\p_F} &
\Hcal_F^{ {\pol}}(\kgot,N) \ar[d]^{\res^1}\\
\Hcal^{ {\pol}}(\kgot,N,N\setminus F') \ar[r]^-{\p_{F'}} & \quad
\Hcal_{F'}^{ {\pol}}(\kgot,N) }
\end{equation}
is commutative. Here $\res_1$ and $\res^1$ are the restriction morphisms (see
(\ref{eq:res-F}) and (\ref{eq:res-F-local})).

$\bullet$ If $F_1,F_2$ are closed invariant subsets of $N$, then the
diagram 
\begin{equation}\label{eq:fonctoriel-p-produit}
\xymatrix@C=8mm{ \Hcal^{ {\pol}}(\kgot,N,N\setminus F_1)
\ar[d]^{\p_{F_1}} \!\!\!\!\!\!\!\!\!\!\!\! & \times \quad
\Hcal^{{\pol}}(\kgot,N,N\setminus F_2)\ \
\ar[d]^{\p_{F_2}}\ar[r]^-{\lozenge} &
\ \  \Hcal^{ {\pol}}(\kgot,N,N\setminus (F_1\cap F_2)) \ar[d]^{\p_{F_1\cap F_2}}\\
\Hcal_{F_1}^{ {\pol}}(\kgot,N)  \!\!\!\!\!   & \times \quad \quad
\Hcal_{F_2}^{{\pol}}(\kgot,N)\quad  \ar[r]^-{\wedge} & \quad
\Hcal_{F_1\cap F_2}^{ {\pol}}(\kgot,N)}
\end{equation}
is commutative.

\end{prop}

\begin{proof}
The proof of the first point is left to the reader. Let us prove
the second point. Let $W$ be an invariant open neighborhood of
$F_1\cap F_2$. Let $V_1,V_2$ be  invariant open neighborhoods
respectively  of $F_1$ and $F_2$ such that $V_1\cap V_2\subset W$.
Choose $\chi_i\in\f(N)^K$ supported in $V_i$ and equal to $1$ in a
neighborhood of $F_i$. Then $\chi_1\chi_2$ is supported in $W$ and
equal to $1$ in a neighborhood of $F_1 \cap F_2$. Let $\Phi_1+
\Phi_2={\rm 1}_{N\setminus (F_1\cap F_2)}$ be a partition of unity
relative to the decomposition $N\setminus (F_1\cap
F_2)=(N\setminus F_1)\cup (N\setminus F_2)$, and where the
function $\Phi_i$ are $K$-invariant.

Then one checks easily that
\begin{eqnarray}\label{eq:p-chi-morphism}
\lefteqn{\p^{\chi_1}(a_1)\wedge\ \p^{\chi_2}(a_2)-
\p^{\chi_1\chi_2}(a_1\diamond_\Phi a_2) =}\\
& &D\Big((-1)^{|a_1|}\chi_2
d\chi_1(\beta_1\Phi_2\beta_2)-(-1)^{|a_1|}\chi_1 d\chi_2
(\Phi_1\beta_1\beta_2)\Big)\nonumber
\end{eqnarray}
for any $D_{\rm rel}$-closed forms
$a_i=(\alpha_i,\beta_i)\in\Acal^{ {\pol}}(\kgot,N,N\setminus F_i)$.
Remark that $\Phi_1\beta_1\beta_2$ is defined on $N\setminus F_2$,
so that $d\chi_2 (\Phi_1\beta_1\beta_2)$ is well defined on $N$
and supported in $V_2$. Thus the equivariant form $(-1)^{|a_1|}\chi_2
d\chi_1(\beta_1\Phi_2\beta_2)-(-1)^{|a_1|}\chi_1 d\chi_2
(\Phi_1\beta_1\beta_2) $ is well defined on $N$ and supported in
$V_1\cap V_2\subset W$. Thus $\p^{\chi_1}(a_1)\wedge\
\p^{\chi_2}(a_2)$ and $\p^{\chi_1\chi_2}(a_1\diamond_\Phi a_2)$
are equal  in $\Hcal^{\pol}_W(\kgot,N)$.  As this holds for any
neighborhood $W$ of $F_1\cap F_2$, this proves that
$\p_{F_1}(a_1)\wedge\p_{F_2}(a_2)=\p_{F_1\cap F_2}(a_1\diamond
a_2)$.

\end{proof}

\bigskip

If we take $F'=N$ in (\ref{eq:fonctoriel-p-restriction}), we see
that the  map $\p_F:\Hcal^{ {\pol}}(\kgot,N,N\setminus F)\to
\Hcal_F^{ {\pol}}(\kgot,N)$ factors the natural map $\Hcal^{
{\pol}}(\kgot,N,N\setminus F)\to \Hcal^{ {\pol}}(\kgot,N)$.

\bigskip

By the same formulae, we define the morphism of $\Zbb_2$-graded
algebras:
\begin{equation}\label{eq:p-F-infty}
\p_F:\Hcal^{ \infty}(\kgot,N,N\setminus F)\to \Hcal_F^{
\infty}(\kgot,N),
\end{equation}

which enjoys the same properties:
\begin{prop}

$\bullet$ If $F\subset F'$ are closed invariant subsets of $N$, then
the diagram
\begin{equation}\label{eq:fonctoriel-p-restriction-infty}
\xymatrix@C=1cm{ \Hcal^{ \infty}(\kgot,N,N\setminus F)
\ar[d]^{\res_1}\ar[r]^-{\p_F} &
\Hcal_F^{ \infty}(\kgot,N) \ar[d]^{\res^1}\\
\Hcal^{ \infty}(\kgot,N,N\setminus F') \ar[r]^-{\p_{F'}} & \quad
\Hcal_{F'}^{ \infty}(\kgot,N) }
\end{equation}
is commutative. Here $\res_1$ and $\res^1$ are the restriction morphisms (see
(\ref{eq:res-F}) and (\ref{eq:res-F-local})).

$\bullet$ If $F_1,F_2$ are closed invariant subsets of $N$, then the
diagram 
\begin{equation}\label{eq:fonctoriel-p-produit-infty}
\xymatrix@C=8mm{ \Hcal^{ \infty}(\kgot,N,N\setminus F_1)
\ar[d]^{\p_{F_1}} \!\!\!\!\!\!\!\!\!\!\!\! & \times \quad
\Hcal^{\infty}(\kgot,N,N\setminus F_2)\ \
\ar[d]^{\p_{F_2}}\ar[r]^-{\lozenge} &
\ \  \Hcal^{ \infty}(\kgot,N,N\setminus (F_1\cap F_2)) \ar[d]^{\p_{F_1\cap F_2}}\\
\Hcal_{F_1}^{ \infty}(\kgot,N)  \!\!\!\!\!   & \times \quad \quad
\Hcal_{F_2}^{\infty}(\kgot,N)\quad  \ar[r]^-{\wedge} & \quad
\Hcal_{F_1\cap F_2}^{ \infty}(\kgot,N)}
\end{equation}
is commutative.

\end{prop}

\medskip

If $F$ is  a compact $K$-invariant subset of $N$, we have a natural morphism
\begin{equation}\label{eq:map-p-compact}
\Hcal^{\pol}_F(\kgot,N)\to \Hcal_c^{\pol}(\kgot,N).
\end{equation}

The composition of $\p_F$ with this morphism is the morphism
$\p_c$ defined in (\ref{eq:p-compact})~:
$\p_c(\alpha,\beta)$ is the class of
$\p^{\chi}(\alpha,\beta)=\chi\alpha+d\chi \beta$, where
 $\chi\in\f(N)^K$ has a compact support and is equal to $1$ in a
neighborhood of $F$.

\subsection{Integral over the fiber in relative cohomology}\label{sec:integral-fiber}

Let $\pi: N\to B$ be a $K$-invariant fibration, with oriented
fibers : the orientation is assumed to be invariant relatively to the action of
$K$. Recall the definition of  $\Acal^{\pol}_{\cf}(\kgot,N)$,
the sub-algebra of $\Acal^{\pol}(\kgot,N)$ formed by the
equivariant forms which have a support that intersects the fibers
of $\pi$ in compact subsets, and of $\Hcal^{\pol}_{\cf}(\kgot,N)$
the corresponding cohomology space.

We have an integration morphism
$\pi_*:\Acal^{\pol}_{\cf}(\kgot,N)\to\Acal^{\pol}(\kgot,B)$
satisfying the following rules:
\begin{equation}\label{eq:D-pi-star}
\pi_*(D\alpha)=D(\pi_*(\alpha)),
\end{equation}
\begin{equation}\label{eq:pi-pi-star}
\pi_*(\pi^*(\gamma)\wedge\alpha)=
\gamma\wedge\pi_*(\alpha),
\end{equation}
for $\alpha\in \Acal^{\pol}_{\cf}(\kgot,N)$ and
$\gamma\in \Acal^{\pol}(\kgot,B)$. Thanks to
(\ref{eq:D-pi-star}) the integration morphism descends to the
cohomology :
$$
\pi_*:\Hcal^{\pol}_{\cf}(\kgot,N)\to \Hcal^{\pol}_c(\kgot,B).
$$
Note that $\pi_*$ sends $\Hcal^{\pol}_c(\kgot,N)$ to $\Hcal^{\pol}_c(\kgot,B)$
and that (\ref{eq:pi-pi-star}) is still valid in cohomology.

\begin{rem}\label{rem:pi-star-U}
If $U$ is an invariant open subset of $B$, we have an integration
morphism $\pi_*:\Acal^{\pol}_{\cf}(\kgot,\pi^{-1}(U))\to
\Acal^{\pol}(\kgot,U)$.
\end{rem}

Let $F$ be a $K$-invariant closed subset of $N$ which is {\bf compact}.
We will define an  integration morphism
$\pi_*:\Hcal^{\pol}(\kgot,N,N\setminus F)\to \Hcal^{ \pol}(\kgot,B,B\setminus\pi(F))$
which makes the following diagram
\begin{equation}\label{eq:fonctoriel-integration-fibre}
\xymatrix@C=8mm{
\Hcal^{ \pol}(\kgot,N,N\setminus F)\ar[d]^{\p_{c}}\ar[r]^-{\pi_*} &
\Hcal^{ \pol}(\kgot,B,B\setminus \pi(F)) \ar[d]^{\p_{c}}\\
\Hcal_{c}^{ \pol}(\kgot,N)  \quad  \ar[r]^-{\pi_*} &
\quad \Hcal_{c}^{ \pol}(\kgot,B)
}
\end{equation}
commutative.

To perform the integration, it is natural to choose a
representative $[\alpha,\beta]$ of a relative cohomology class $a$
where $\alpha$ and $\beta$ are  supported near $F$. This can be
done via the Lemma \ref{lem:excision}: let us choose  a
$K$-invariant function $\chi$ which is compactly supported and is equal to $1$ on a neighborhood
of $F$. Let $(\alpha,\beta)\in
\Acal^{\pol}(\kgot, N, N\setminus F)$. Then the equivariant form $\chi \alpha
+d\chi\beta$ is compactly supported and can be integrated over the
fiber. Similarly, the form $\chi \beta$ belongs to
$\Acal^{\pol}_{\cf}(\kgot,\pi^{-1}(B\setminus \pi(F)))$ and can be
integrated over the fiber. The expression
$$
\pi_*^\chi(\alpha,\beta):=(\pi_*(\chi
\alpha +d\chi\beta), \pi_*(\chi \beta))
$$
defines an element in $\Acal^{\pol}(\kgot, B, B\setminus \pi(F))$. Since we have the relation
$\pi_*^\chi\circ D_{\rel}= D_{\rel}\circ\pi_*^\chi$ the map $\pi_*^\chi$
descends to cohomology.

Furthermore, if $(\alpha,\beta)$ is $D_{\rm rel}$-closed, and
$\chi_1$, $\chi_2$ are two  different choices  of functions $\chi$, we
verify that
$$
\pi_*^{\chi_1}(\alpha,\beta)-\pi_{*}^{\chi_2}(\alpha,\beta)=
D_{\rel}\Big(\pi_*((\chi_1-\chi_2)\beta),0\Big)
$$
so that the map $\pi^\chi_*$ is independent of the choice of $\chi$.

This allows us to make the following definition.

\begin{defi}
Let  us choose $\chi\in \f(N)$ a $K$-invariant function
identically equal to $1$ on a neighborhood of
$F$, and with compact support. Then we define
$$
\pi_*: [\Hcal^{\pol}(\kgot, N,N\setminus F)]^{k}\longrightarrow
[\Hcal^{\pol}(\kgot, B,B\setminus \pi(F))]^{k-\dim N+\dim B}
$$
by
the formula : $\pi_*([\alpha,\beta])=[\pi_*(\chi \alpha +d\chi\beta), \pi_*(\chi
\beta)]$.
\end{defi}

If $B$ is a point, then the integral of $a\in\Hcal^{\pol}(\kgot, N,N\setminus F)$
over the fiber is just
the integral over $N$ of the class $\p_c(a)\in \Hcal_c^{\pol}(\kgot,N)$.

The same definition makes sense for equivariant forms with
$\Ccal^{\infty}$-coefficients and defines a map
$\pi_*:\Hcal^{\infty}(\kgot, N,N\setminus F)\to \Hcal^{\infty}(\kgot,
B,B\setminus \pi(F))$.

We now prove

\begin{prop}\label{prop:intfib}
\begin{itemize}
\item The diagram (\ref{eq:fonctoriel-integration-fibre}) is commutative

\item Let $F_2$ be a {\em compact} $K$-invariant set of $N$. Let $F_1$ be a closed
$K$-invariant set in $B$. Then, for any $a \in \Hcal^{\pol}(\kgot, B,B\setminus
F_1)$ and $b\in \Hcal^{\pol}(\kgot,N,N\setminus F_2)$, we have
$$\pi_*(\pi^*a\diamond b)=a\diamond \pi_*(b)$$
in $\Hcal^{\pol}(\kgot, B,B\setminus (F_1\cap \pi(F_2))).$
\end{itemize}
\end{prop}

\begin{proof} Let $\chi\in\f(N)^K$ and $\chi'\in\f(B)^K$ be two compactly supported functions :
$\chi$ is identically equal to $1$ on a neighborhood of $F$ and $\chi'$ is
identically equal to $1$ on a neighborhood of $\pi(F)$. For $[\alpha,\beta]\in\Hcal^{\pol}(\kgot,
N,N\setminus F)$, the equivariant class $\p_c\circ\,\pi_*[\alpha,\beta]$ is represented by
$$
\p^{\chi'}(\pi_*^\chi(\alpha,\beta))=\chi'\pi_*(\chi\alpha+d\chi\beta)+d\chi'\pi_*(\chi\beta).
$$
On the other hand, the equivariant class $\pi_*\circ\p_c[\alpha,\beta]$ is represented by
$\pi_*(\chi\alpha+d\chi\beta)$. We check that
$$
\p^{\chi'}(\pi_*^\chi(\alpha,\beta))-\pi_*(\chi\alpha+d\chi\beta)=D\Big((\chi'-1)\pi_*(\chi\beta)\Big),
$$
where $(\chi'-1)\pi_*(\chi\beta)$ is an equivariant form on $B$ with compact support.
Then, the first point is proved.

Let us prove the second point. We work with the invariant open subsets
$U=N\setminus ( \pi^{-1}(F_1)\cap F_2), U_1= N\setminus \pi^{-1}(F_1), U_2= N\setminus F_2$ of $N$,
and the invariant open subsets $U'=B\setminus (F_1\cap \pi(F_2)), U'_1=B\setminus F_1,
U'_2=B\setminus \pi(F_2)$ of $B$.  Note that $\pi^{-1}U'_1=U_1$ and
$\pi^{-1}U'_2\subset U_2$ : hence $\pi^{-1}U'\subset U$.

Let $\Phi_1+\Phi_2=1$ be a partition of unity on $U=U_1\cup U_2$,
and let $\Phi'_1+\Phi'_2=1$ a partition of unity on $U'=U'_1\cup
U'_2$ : all the functions are supposed invariant. Let
$(\alpha_1,\beta_1)$ and $(\alpha_2,\beta_2)$ be respectively the
representatives of $a\in\Hcal^{\pol}(\kgot,B,B\setminus F_1)$ and
$b\in \Hcal^{\pol}(\kgot, N,N\setminus F_2)$. The equivariant
forms $\alpha_2,\beta_2$ are chosen so that their supports belong
to a compact neighborhood of $F_2$.

Then $\pi_*(\pi^*(a)\diamond b)$ is represented by $(\alpha_1\pi_*(\alpha_2),\beta)$ with
\begin{equation*}\label{eq:beta-pi-star}
\beta=\underbrace{\beta_1\pi_*(\Phi_1\alpha_2)}_{\beta(1)}+
\underbrace{(-1)^{|a|}\alpha_1\pi_*(\Phi_2\beta_2)}_{\beta(2)}+
\underbrace{\beta_1\pi_*(d\Phi_1\beta_2)}_{\beta(3)}.
\end{equation*}
On the other hand $a\diamond \pi_*(b)$ is represented by $(\alpha_1\pi_*(\alpha_2),\beta')$ with
\begin{equation*}\label{eq:betaprime-pi-star}
\beta'=
\underbrace{\Phi'_1\beta_1\pi_*(\alpha_2)}_{\beta'(1)}+
\underbrace{(-1)^{|a|}\Phi'_2\alpha_1\pi_*(\beta_2)}_{\beta'(2)}+
\underbrace{\beta_1d\Phi'_1\pi_*(\beta_2)}_{\beta'(3)}.
\end{equation*}
Note that the equivariant forms $\beta(i),\beta'(i)$ are well defined on
$B\setminus (F_1\cap \pi(F_2))$.

\begin{lem}
The equivariant form $\delta=\Phi'_2\pi_*(\beta_2)-\pi_*(\Phi_2\beta_2)$
is defined on $B\setminus (F_1\cap \pi(F_2))$ and supported on $B\setminus (F_1\cup \pi(F_2))$.
We have
$$
\beta-\beta'= D\left(\beta_1\wedge\delta\right),
$$
where the equivariant form $\beta_1\wedge\delta$  is defined on $B\setminus (F_1\cap \pi(F_2))$.
It gives the
following relation in $\Acal^{\pol}(\kgot, B,B\setminus (F_1\cap \pi(F_2)))$ :
$(\alpha_1\pi_*(\alpha_2),\beta')-(\alpha_1\pi_*(\alpha_2),\beta)=
D_{\rm rel}\left(0,\beta_1\delta\right)$.
\end{lem}

\begin{proof} The invariant function  $\Phi_2$ is defined on $U$, supported on
$U_2$, and equal to $1$ on $U_2\setminus(U_1\cap U_2)$. Then its restriction
$\Phi_2\vert_{\pi^{-1}(U')}$ is supported on $\pi^{-1}(U'_2)$ and equal to $1$ on
$\pi^{-1}(U'_2)\setminus(U_1\cap \pi^{-1}(U'_2))$. Similarly the function $\pi^*\Phi'_2$ is defined on
$\pi^{-1}(U')$,  supported on $\pi^{-1}(U'_2)$ and equal to $1$ on
$\pi^{-1}(U'_2)\setminus(U_1\cap \pi^{-1}(U'_2))$. Hence the difference
$\pi^*\Phi'_2-\Phi_2\vert_{\pi^{-1}(U')}$ is defined on $\pi^{-1}(U')$ and supported on
$U_1\cap \pi^{-1}(U'_2)$. This shows that
$\delta=\pi_*((\pi^*\Phi'_2-\Phi_2\vert_{\pi^{-1}(U')})\beta_2)$ is defined on $U'$
and supported on $U'_1\cap U'_2$.

We have
\begin{eqnarray*}
\beta_1\wedge D\left(\delta\right)
&=& \beta_1d\Phi'_2\pi_*(\beta_2)-\beta_1\pi_*(d\Phi_2\beta_2) +
\beta_1\wedge\Big(\Phi'_2\pi_*(\alpha_2)-\pi_*(\Phi_2\alpha_2)\Big) \\
&=& -\beta'(3)+\beta(3)+
\beta_1\wedge\Big((1-\Phi'_1)\pi_*(\alpha_2)-\pi_*((1-\Phi_1)\alpha_2)\Big) \\
&=&-\beta'(3)+\beta(3) -\beta'(1)+ \beta(1).
\end{eqnarray*}
We see also that
\begin{eqnarray*}
(-1)^{|a|-1}D(\beta_1)\wedge \delta
&=& -(-1)^{|a|}\Phi'_2\alpha_1\pi_*(\beta_2)
+(-1)^{|a|}\alpha_1\pi_*(\Phi_2\beta_2)\\
&=& -\beta'(2)+\beta(2).
\end{eqnarray*}

Finally we have proved that $D(\beta_1\wedge \delta)=
\beta(1)+\beta(2)+\beta(3) -\beta'(1) -\beta'(2)-\beta'(3)$.

\end{proof}
\end{proof}\bigskip

\begin{rem}\label{rem:integration-F-non-compatc}
If $F$ is a {\em closed but non-compact} subset of $N$, the morphism
$\pi_*:\Hcal^{\pol}(\kgot,N,N\setminus F) \to \Hcal^{\pol}(\kgot,B,B\setminus\pi(F))$
is still defined in the case where the image $\pi(F)$ is closed and the
intersection of $F$ with each fiber is compact. The second point of Proposition
(\ref{prop:intfib}) still holds in this case.
\end{rem}

\subsection{The Chern-Weil construction}

Let $\pi:P\to B$ be a principal bundle with structure group $G$.
Let $\omega\in (\Acal^{1}(P)\otimes \ggot)^G$ be a  connection one
form on $P$, with curvature form
$\Omega=d\omega+\frac{1}{2}[\omega,\omega]$.

For any $G$-manifold $Z$, we define $\Zcal=P\times_G Z$. If $F$ is
a $G$ invariant closed subset of $Z$, then $\Fcal:=P\times_G F$ is
a closed subset of $\Zcal$. We consider the Chern-Weil
homomorphism
$$
\phi_\omega:\Acal^{\pol}(\ggot,Z,Z\setminus F)) \longrightarrow
\Acal(\Zcal,\Zcal\setminus \Fcal)
$$
defined by $\phi_{\omega}(\alpha,\beta)=(\phi_\omega^{Z}(\alpha),
\phi_\omega^{Z\setminus F}\beta)$. We have the relation:
\begin{prop}
We have $d_{\rel}\circ\phi_{\omega}=\phi_{\omega}\circ D_{\rel}$.
\end{prop}

We can repeat the construction above in the equivariant case. If
$P$ is a $G$-principal bundle with left action of $K$, and $F$
is a  $G\times K$-invariant closed subset of $P$, we define
$$
\phi_\omega:\Acal^{\pol}(\ggot,Z, Z\setminus F) \longrightarrow
\Acal^{\pol}(\kgot, \Zcal,\Zcal\setminus \Fcal)
$$
by $\phi_\omega(\alpha,\beta)= (\phi_\omega^Z\alpha,
\phi_\omega^{Z\setminus F}\beta)$.

\begin{prop}
The map $\phi_\omega:\Acal^{\pol}(\ggot, Z,Z\setminus F)\to
\Acal^{\pol}(\kgot, \Zcal,\Zcal\setminus \Fcal)$ satisfies
$$
\phi_\omega\circ D_{\rel}=D_{\rel}\circ\phi_\omega.$$

On the left side the equivariant differential is with respect to
the action of $G$ while on the right side, this is with respect to
the group $K$.
\end{prop}

\section{ Explicit formulae for  Thom  Classes in relative cohomology}
\label{sec:chg-sigma-equi}

\subsection{The  equivariant Thom forms of a vector space}
\label{subsec:Thom-vector}

Let $V$ be an  Euclidean oriented vector space of dimension $d$.
Consider the group ${\rm SO}(V)$ of orthogonal transformations of
$V$ preserving the orientation. Let $\sogot(V)$ be its Lie
algebra.  Consider the projection
$\pi:V\to \{{\rm pt}\}$ and the closed subset $F=\{0\}\subset V$. We denote
$\int_V$ the integration morphism $\pi_*:\Hcal^{\pol}(\sogot(V),V,V\setminus\{0\})\to
\Ccal^{\pol}(\sogot(V))^{\So(V)}$ defined in Section \ref{sec:integral-fiber}.

In this section, we will describe a generator over
$\Ccal^{\pol}(\sogot(V))^{\So(V)}$ of  the equivariant relative
cohomology space of the pair $(V,V\setminus \{0\})$. This explicit
representative is basically due to Mathai-Quillen
\cite{Mathai-Quillen}. As a consequence, we obtain the following
theorem.
\begin{theo}
There exists a unique class  $\tur(V)$  in
$\Hcal^{\pol}(\sogot(V),V,V\setminus\{0\})$ such that $\int_V {\tur}(V)=1$. This class is called the
relative Thom class.
\end{theo}

Before establishing the unicity, a {\em closed} form  $\alpha\in
\Acal^{\pol}(\sogot(V),V,V\setminus \{0\})$, or in
$\Acal_c^{\pol}(\sogot(V),V)$, or in $\Acal_{\dr}^{\pol}(\sogot(V),V)$, of
integral $1$ will be called a Thom form. A Thom class will be the class defined
by a Thom form.

\bigskip

We start by constructing Thom forms in
the spaces $\Acal^{\pol}(\sogot(V),V,V\setminus
\{0\})$, $\Acal_c^{\pol}(\sogot(V),V)$ and
$\Acal_{\dr}^{\pol}(\sogot(V),V)$.
The Lie algebra $\sogot(V)$ is identified with the Lie algebra of
antisymmetric endomorphisms  of $V$.
Let $(e_1,\ldots,e_d)$ be an oriented orthonormal basis of $V$.
Denote by $\bere:\wedge V\to \Rbb$ the Berezin integral normalized
by $\bere(e_1\wedge \cdots\wedge e_d)=1$.

\begin{defi}
The pfaffian of $X\in \sogot(V)$ is defined by
$$\Pf(X)=\bere(\e^{\sum_{k<l}\langle X e_k,e_{\ell}\rangle e_k\wedge
e_{\ell}}).$$
Here the exponential is computed in the algebra $\wedge V$.
\end{defi}

Recall that  $\Pf(X)$ is an ${\rm SO}(V)$-invariant polynomial on
$\sogot(V)$ such that $(\pf(X))^2=\det_V(X).$ In particular this
polynomial is identically equal to $0$, if $d$ is odd.
We also denote by $\Pf(X)\in\Acal^{\pol}(\sogot(V),V)$
the function on $V$ identically equal
to $\Pf(X)$.

Let $x_k=(e_k,x)$ be the  coordinates on $V$. We consider the
equivariant map $f_t:\sogot(V)\to\Acal(V)\otimes\wedge V$ given by the
formula
\begin{equation}\label{eq-map-f-t-defi}
f_t(X)= -t^2\|x\|^2 + t\sum_k dx_k e_k + \frac{1}{2}\sum_{k<l}
\langle X e_k,e_{\ell}\rangle e_k\wedge e_\ell.
\end{equation}

For any (real or complex) vector space $\Acal$, the Berezin integral is extended to a map
$\bere : \Acal\otimes\wedge V\to\Acal$ by $\bere(\alpha\otimes \xi)=\alpha\bere(\xi)$ for
$\alpha\in\Acal$ and $\xi\in\wedge V$. If ${\rm d}: \Acal\to \Acal$ is a linear map,
we extend it on $\Acal\otimes\wedge V$ by ${\rm d}(\alpha\otimes \xi)={\rm d}(\alpha)\otimes\xi$.
Note that this extension satisfies $\bere\circ \,{\rm d}={\rm d}\circ \bere$.

The Berezin integral is the ``super-commutative" analog of the
super-trace for endomorphisms of a super-space.  We now  give a
construction for the relative Thom form, analogous to Quillen's
construction of the Chern character. We will discuss Quillen's
construction of the Chern character in the next chapter (the
formulae we give here for the ``super-commutative case" are
strongly inspired by the formulae for the curvature of the
super-connection attached to the Bott symbol: see Section
\ref{sec:bott-class-V}).

 We consider the $\So(V)$-equivariant forms on $V$
defined by
\begin{eqnarray}\label{eq:C-eta-wedge}
\mathrm{C}_{\wedge}^t(X)&:=&\bere\left(\e^{f_t(X)}\right),\\
\eta_{\wedge}^t(X)&:=& - \bere\left((\sum_k x_k e_k)\e^{f_t(X)}\right),
\end{eqnarray}
for $X\in\sogot(V)$. Here the exponential is computed in the super-algebra
$\Acal(V)\otimes\wedge V$.

\begin{lem}\label{thomgaussianclosed}
The equivariant form $ \mathrm{C}_{\wedge}^t(X)$ is closed.
Furthermore,
\begin{equation}\label{eq:Thom}
\frac{d}{dt}\mathrm{C}_{\wedge}^t=-D(\eta_{\wedge}^t).
\end{equation}

\end{lem}

\begin{proof}
The proof of the first point is given in \cite{B-G-V} (Chapter 7,
Theorem 7.41). We recall the proof. If $e\in V$, we denote by
$\iota_{\wedge}(e)$ the derivation of $\wedge V$ such that
$\iota_{\wedge}(e)v=\langle e,v\rangle$ when $v\in \wedge^1 V=V$.
We extend it to a derivation of $\Acal(V)\otimes\wedge V$. We denote by
$\iota_{\wedge}(\x)$ the operator
$\sum_{k}x_k\iota_{\wedge}(e_k)$. The equivariant derivative $D$ which is
defined on $\Ccal^{\rm pol}(\sogot(V))\otimes \Acal(V)$ is extended to a derivation of
$\Ccal^{\rm pol}(\sogot(V))\otimes \Acal(V)\otimes \wedge V$. We have
$\bere\circ \, D= D\circ \bere$.

It is easy to verify that
\begin{equation}\label{eft}
(D -2t \iota_{\wedge}(\x))f_t(X)=0.
\end{equation}
The exponential $\e^{f_t(X)}$ satisfies also
$(D-2t\iota_{\wedge}(\x))(\e^{f_t(X)})=0$, since $D$ and
$\iota_{\wedge}(\x)$ are derivations. The Berezin integral is such
that $\bere(\iota_{\wedge}(\x)\alpha)=0$ for any $\alpha\in
\Acal(V)\otimes\wedge V$. This shows that $D\left(\bere(\e^{f_t(X)})\right)=0$.

Let us prove the second point. As $(D -2t\iota_{\wedge}(\x))f_t(X)=0$, we have
\begin{eqnarray*}
D\circ \bere \left((\sum x_k e_k)\e^{f_t(X)}\right)
&=&\bere\circ\, (D-2t\iota_{\wedge}(\x))\left((\sum x_k e_k) \e^{f_t(X)}\right)\\
&=& \bere\left(\big((D-2t \iota_{\wedge}(\x))\cdot (\sum x_k e_k)\big) \e^{f_t(X)}\right)\\
&=&\bere\left( (\sum dx_k e_k-2t\|x\|^2)\e^{f_t(X)}\right)\\
&=&\bere\left(\frac{d}{dt}\e^{f_t(X)}\right).
\end{eqnarray*}
\end{proof}

\bigskip

When $t=0$, then $\mathrm{C}_{\wedge}^0(X)$ is just equal to
$\Pf(\frac{X}{2})=\frac{1}{2^{d/2}}\Pf(X).$ When $t=1$, then
$\mathrm{C}_{\wedge}^1(X)=\bere(\e^{f_1(X)})=e^{-\|x\|^2}Q(dx,X)$ is a
closed equivariant class with a Gaussian look on $V$ (with
$Q(dx,X)$  a polynomial in $dx,X$ that we will write more
explicitly in a short while). This form was considered by
Mathai-Quillen in \cite{Mathai-Quillen}.

\begin{defi}
The Mathai-Quillen form is the closed equivariant form on $V$
defined by
$$
\mathrm{C}_{\wedge}^1(X):=\bere\left(\e^{f_1(X)}\right),\quad X\in\sogot(V).
$$
\end{defi}

\bigskip

 We have $\eta_{\wedge}^t=\e^{-t^2\|x\|^2}Q(t,X,x,dx)$ where $Q(t,X,x,dx)$
depends polynomially of $(t,X,x,dx)$ where $t\in \Rbb, X\in
\sogot(V),x\in V, dx\in \Acal^1(V)$. Thus, if $x\neq 0$, when $t$
goes to infinity, $\eta_{\wedge}^t$ is an exponentially decreasing
function of $t$. We can thus define the following equivariant form
on $V\setminus\{0\}$ :
\begin{equation}\label{def:beta-form}
\beta_{\wedge}(X)=\int_{0}^{\infty}\eta_{\wedge}^t(X)dt,\quad
X\in\sogot(V).
\end{equation}
If we integrate (\ref{eq:Thom}) between  $0$ and $\infty$, we get
$$
\mathrm{C}_{\wedge}^0=D(\beta_{\wedge})\quad \mathrm{on} \quad
V\setminus\{0\}.
$$

Thus the couple $(\mathrm{C}^0_\wedge,\beta_\wedge)$ defines a canonical
relative class
\begin{equation}\label{def:pf.beta}
\left[\Pf(\hbox{$\frac{X}{2}$}),\beta_\wedge(X)\right]\in\Hcal^{\pol}(\sogot(V),V, V\setminus\{0\})
\end{equation}
of degree equal to $\dim(V)$.


Consider now the equivariant cohomology with compact support of
$V$. Following (\ref{eq:p-compact}) we have a map $\p_c:\Hcal^{\pol}(\sogot(V),V,
V\setminus\{0\})\to \Hcal^{\pol}_c(\sogot(V),V)$. Thus
$\left[\Pf(\hbox{$\frac{X}{2}$}),\beta_\wedge(X)\right]$
provides us a class
\begin{equation}\label{def:pf.beta.compact}
\mathrm{C}_V:=
\p_c\left(\left[\Pf(\hbox{$\frac{X}{2}$}),\beta_\wedge(X)\right]\right)\in \Hcal^{\pol}_c(\sogot(V),V).
\end{equation}
As the map $\p_c$ commutes with the
integration over the fiber (see (\ref{eq:fonctoriel-integration-fibre})), it is the same thing
to compute the integral of (\ref{def:pf.beta}) or of (\ref{def:pf.beta.compact}).

The orientation on $V$ is given by $dx_1\wedge\cdots\wedge
dx_d$.

\begin{prop}\label{prop:Thom}
Let $\chi\in\f(V)$ be an ${\rm SO}(V)$-invariant function with
compact support and equal to $1$ in a neighborhood of $0$. The
form
$$
\mathrm{C}_V^\chi(X)= \chi \Pf(\hbox{$\frac{X}{2}$})
 +d\chi \int_{0}^\infty \eta_{\wedge}^t(X) dt
$$
is a closed equivariant form with compact support on $V$.  Its cohomology class
in $\Hcal^\infty_c(\sogot(V),V)$ coincides with
$\mathrm{C}_V$~: in particular, it does
not depend of the choice of $\chi$. For every $X\in \sogot(V)$,
$$
\frac{1}{\epsilon_d}\int_V\mathrm{C}^\chi _V(X)=1,
$$
with $\epsilon_d=(-1)^{\frac{d(d-1)}{2}}\pi^{d/2}$. Thus
$\frac{1}{\epsilon_d}\mathrm{C}^\chi_V(X)$ is a Thom form in $\Acal_c^{\pol}(\sogot(V),V)$.

\end{prop}
\begin{proof}
The first assertions are  immediate to prove, since
$\mathrm{C}_V=$\break $\p_c\left(\left[\Pf(\hbox{$\frac{X}{2}$}),\beta_\wedge(X)\right]\right)$.
The class of $\mathrm{C}_V^\chi$ in compactly supported
cohomology do not depend of $\chi$. To compute its integral, we
may choose  $\chi= f(\|x\|^2)$ where $f\in \f(\Rbb)$ has a compact
support and is equal to $1$ in a neighborhood of $0$. The
component of maximal degree of the differential form $d\chi\wedge
\eta^t_\wedge(X)$ is $$-2(-1)^{\frac{d(d-1)}{2}}
t^{d-1}f'(\|x\|^2)\|x\|^2 e^{-t^2\|x\|^2} dx_{1}\wedge\cdots\wedge
dx_{d}.$$ Hence, using the change of variables $x\to
\frac{1}{t}x$,
\begin{eqnarray*}
  \int_V \mathrm{C}_V^\chi(X) &=&
  -2(-1)^{\frac{d(d-1)}{2}}\int_{0}^\infty t^{d-1}\left(\int_V
  f'(\|x\|^2)\|x\|^2 e^{-t^2\|x\|^2} dx\right)dt\\
   &=&(-1)^{\frac{d(d-1)}{2}}\int_{0}^\infty \underbrace{\left(\int_V
  f'(\frac{\|x\|^2}{t^2})(\frac{-2\|x\|^2}{t^3})e^{-\|x\|^2}
  dx\right)}_{I(t)}dt.
\end{eqnarray*}
Since for $t>0$, $I(t)=\frac{d}{dt}\left(\int_V
  f(\frac{\|x\|^2}{t^2})e^{-\|x\|^2} dx\right)$, we have
$\int_V \mathrm{C}_V^\chi(X)=$ \break $(-1)^{\frac{d(d-1)}{2}}\int_V e^{-\|x\|^2}
dx= (-1)^{\frac{d(d-1)}{2}}\pi^{d/2}.$

\end{proof}

\bigskip

Using Mathai-Quillen form, it is possible to construct
representatives of a Thom form with Gaussian look.

\begin{prop}\label{thomgaussian}
The Mathai-Quillen form $\mathrm{C}_{\wedge}^1(X)$ is a closed equivariant form which belongs to
$\Acal_{\dr}^{\pol}(\sogot(V),V)$. For every $X\in \sogot(V)$,
$$
\frac{1}{\epsilon_d}\int_V \mathrm{C}_{\wedge}^1(X)=1,
$$
with $\epsilon_d=(-1)^{\frac{d(d-1)}{2}}\pi^{d/2}$. Thus
$\frac{1}{\epsilon_d}\mathrm{C}_{\wedge}^1(X)$ is a Thom form in
$\Acal_{\dr}^{\pol}(\sogot(V),V)$.
\end{prop}

\begin{proof}
Indeed, only the term in $dx_1\wedge dx_2\wedge \cdots\wedge dx_d$
will contribute to the integral. This highest term is
$(-1)^{\frac{d(d-1)}{2}} e^{-\|x\|^2}dx$.
\end{proof}
\bigskip

We summarize Propositions \ref{prop:Thom} and \ref{thomgaussian}
in the following theorem.

\begin{theo}\label{theo:thom-form}
Let $V$ be an oriented Euclidean vector space with oriented
orthonormal basis $(e_1,\ldots,e_d)$. Let
$\epsilon_d:=(-1)^{\frac{d(d-1)}{2}}\pi^{d/2}$.
Let $\bere:\wedge V\to \Rbb$
be the Berezin integral. Let, for $X\in \sogot(V)$,
\begin{eqnarray*}
f_t(X)&=& -t^2\|x\|^2 + t\sum_k dx_k e_k + \frac{1}{2}\sum_{k<l}
\langle X e_k,e_{\ell}\rangle e_k\wedge e_\ell,\\
\eta_{\wedge}^t(X)&=&- \bere\left((\sum_k x_k
e_k)\e^{f_t(X)}\right),\\
\beta_{\wedge}(X)&=&\int_{0}^{\infty}\eta_{\wedge}^t(X)dt.
\end{eqnarray*}

$\bullet$ Let $$\tur(V) =
\frac{1}{\epsilon_d}
\left(\Pf(\hbox{$\frac{X}{2}$}),\beta_{\wedge}(X)\right),\quad
X\in\sogot(V).
$$
Then  $\tur(V)\in\Acal^{\pol}(\sogot(V),V,V\setminus \{0\})$ is  a Thom form. It defines a
Thom class, still denoted   $\tur(V)$, in $\Hcal^{\pol}(\sogot(V),V,V\setminus \{0\})$.

$\bullet$
Let $$\tuc(V)=
\frac{1}{\epsilon_d}\mathrm{C}_{V}^\chi(X)=\frac{1}{\epsilon_d}\left( \chi \Pf(\hbox{$\frac{X}{2}$})
 +d\chi \beta_\wedge(X)\right),\quad
X\in\sogot(V).
$$
Then  $\tuc(V)\in\Acal_{c}^{\pol}(\sogot(V),V)$  is a  Thom form.
It defines a
Thom class, still denoted   $\tuc(V)$, in $\Hcal^{\pol}_c(\sogot(V),V)$.

Here $\chi\in\f(V)^{\So(V)}$ is an invariant
function with compact support and equal to $1$ in a neighborhood
of $0$.

$\bullet$
Let
$$\tumq(V)
=\frac{1}{\epsilon_d}\mathrm{C}_{\wedge}^1(X)=\frac{1}{\epsilon_d}\bere\left(\e^{f_1(X)}\right),\quad
X\in\sogot(V)
$$
be the Mathai-Quillen form.
Then  $\tumq(V)\in\Acal_{\dr}^{\pol}(\sogot(V),V)$ is a Thom  form.
It defines a
Thom class, still denoted   $\tumq(V)$, in $\Hcal^{\pol}_{\dr}(\sogot(V),V)$.

\end{theo}

Thus the use of the Berezin integral allowed us to give slim
formulae for Thom forms in relative cohomology, as well as in
compactly supported cohomology or in rapidly decreasing
cohomology.

\bigskip

Let us explain the relation between the Thom forms
$\frac{1}{\epsilon_d}\mathrm{C}_V^\chi$ and
$\frac{1}{\epsilon_d}\mathrm{C}_{\wedge}^1$. For
$t>0$, it is easy to see that the forms $\mathrm{C}_V^\chi(X)$ and
$\mathrm{C}_{\wedge}^t(X)$ differ by the differential of an
equivariant form with Gaussian decay. We could deduce this fact as
a corollary of the unicity theorem, that we will prove soon, but
we prove it here directly by giving explicit transgression
formulae.

\begin{prop}
For any $t>0$, $\mathrm{C}_V^\chi=\mathrm{C}_{\wedge}^t$ in
$\Hcal_{\dr}^{\pol}(\sogot(V),V)$.
\end{prop}

\begin{proof}
Fix $t>0$. The form $\mathrm{C}_\wedge^t=\e^{-t^2\|x\|^2}P(t,dx,X)$ is
rapidly decreasing on $V$, thus  $\mathrm{C}_\wedge^t\in
\Acal_{\dr}^{\pol}(\sogot(V),V)$. Define on $V\setminus\{0\}$ the equivariant form
$\beta_\wedge^t=\int_t^{\infty}\eta_\wedge^s ds$.  We have
$\eta_\wedge^s=\e^{-s\|x\|^2}P(s,X,x,dx)$ where $P(s,X,x,dx)$ depends
polynomially of $(s,X,x,dx)$. For $s=t+u$, $e^{-s\|x\|^2}=
\e^{-t\|x\|^2}\e^{-u\|x\|^2}$, thus $\beta_\wedge^t$ is rapidly decreasing
when $\|x\|$ tends to $\infty$. The
transgression formula (\ref{eq:Thom}) integrated between $t$ and
$\infty$ shows that $\mathrm{C}_\wedge^t=D(\beta_\wedge^t)$ on $V\setminus\{0\}$.
Thus $\mathrm{C}_V^{\chi,t}:=\chi \mathrm{C}_\wedge^t+ d\chi \beta_\wedge^t$ is a closed
equivariant form belonging to $\Acal_{\dr}^{\pol}(\sogot(V),V)$.
 The two following
transgression formulae are evident to prove:
$$
\mathrm{C}_V^\chi-\mathrm{C}_{V}^{\chi,t}=D\left(\chi\int_0^t \eta_\wedge^s ds\right)\quad
{\rm and}\quad
\mathrm{C}_{V}^{\chi,t}-\mathrm{C}_\wedge^t=D\left((\chi-1)\beta_\wedge^t\right).
$$

Thus we obtain $\mathrm{C}_V^\chi-\mathrm{C}_{\wedge}^t=D(\delta)$, where
$$
\delta:=\chi \int_0^t \eta_\wedge^s ds +(\chi-1)\beta_\wedge^t
$$
is an equivariant form defined on $V$ which is rapidly decreasing. So we have
$\mathrm{C}_V^\chi=\mathrm{C}_{\wedge}^t$ in $\Hcal_{\dr}^{\pol}(\sogot(V), V)$.
\end{proof}

\bigskip

\begin{rem}
In the next sections, we will keep the same notations for the Thom classes and for their
representatives defined in Theorem \ref{theo:thom-form}.
\end{rem}

Before going on in proving the unicity, let us give more explicit
formulae for the Thom forms we have constructed.

\subsection{Explicit formulae for the Thom forms of a vector space}\label{subsec:explicit-Thom-vector}

If $I=[i_1,i_2,\ldots,i_p]$ (with $i_1<i_2<\cdots<i_p$) is a
subset of $[1,2,\ldots,d]$, we denote by $e_I=e_{i_1}\wedge e_{i_2}\wedge \cdots \wedge e_{i_p}$.
If $X$ is an antisymmetric matrix, we denote by $X_I$ the sub-matrix
$\left(\langle X e_i,e_j \rangle\right)_{i,j\in I}$, which is viewed as an antisymmetric endomorphism
of the vector space $V_I$ generated by $e_i,i\in I$: let $\Pf(X_I)$ be its pfaffian, where $V_I$ is
oriented by $e_I$. One sees easily that
\begin{equation}\label{eq:pfaffian-gene}
\e^{\sum_{k<l}\langle X e_k,e_{\ell}\rangle e_k\wedge
e_{\ell}}=\sum_{I}\Pf(X_I)\, e_I \quad {\rm in}\quad \wedge V.
\end{equation}
Only those $I$ with $|I|$ even will
contribute to the sum (\ref{eq:pfaffian-gene}),  as otherwise the pfaffian of $X_I$
vanishes.

If $I$ and $J$ are two disjoint subsets of $\{1,2,\ldots,d\}$, we denote by
$\epsilon(I,J)$ the sign such that $e_I\wedge e_J=\epsilon(I,J)\, e_{I\cup J}$.

\begin{prop}\label{prop:explicit}

$\bullet$ We have
$\tur(V)(X)=\frac{1}{\epsilon_d}\left(\Pf(\hbox{$\frac{X}{2}$}),\beta_{\wedge}(X)\right)$
with
$$
\beta_\wedge(X)=\sum_{k,I,J} \gamma_{(k,I,J)}
 \Pf\left(\hbox{$\frac{X_I}{2}$}\right)\frac{x_k dx_J}{\|x\|^{|J|+1}},
$$
with
$$
\gamma_{(k,I,J)}=-\frac{1}{2}(-1)^{\frac{|J|(|J|+1)}{2}}\Gamma\left(\hbox{$\frac{|J|+1}{2}$}\right)
\epsilon(I,J)\epsilon(\{k\},I\cup J).
$$
Here for $1\leq k\leq d$, the sets $I,J$ vary over the subsets of
$\{1,2,\ldots,d\}$ such that  $\{k\}\cup I\cup
J$ is a partition of $\{1,2,\ldots,d\}$. Only those $I$ with $|I|$ even will
contribute to the sum.

$\bullet$ We have
$$
\tuc(V)(X)=\frac{1}{\epsilon_d}\Big(f(\|x\|^2) \Pf(\hbox{$\frac{X}{2}$})+ 2f'(\|x\|^2) (\sum
x_i dx_i)\beta_{\wedge}(X)\Big)
$$
where $f$ is a compactly supported function on $\Rbb$, identically equal to $1$ in a neighborhood of
$0$.

$\bullet$ We have
$$
\tumq(V)(X)=\frac{1}{(\pi)^{d/2}}e^{-\|x\|^2} \sum_{I} (-1)^{\frac{| I |}{2}}\epsilon(I,I')
\Pf(\hbox{$\frac{X_I}{2}$}) dx_{I'}.
$$
Here $I$ runs over the subset of $\{1,2,\ldots,d\}$ with an even number of elements,
and $I'$ denotes the complement of $I$.
\end{prop}

\begin{proof}Since $dx_ke_k$ and $dx_le_l$ commute we have
\begin{eqnarray*}
\e^{t\sum_{k}dx_ke_k}&=&\prod_k(1+tdx_ke_k)\\
&=& \sum_J (-1)^{\frac{|J|(|J|-1)}{2}}dx_J\,e_J\, t^{|J|}.
\end{eqnarray*}
If we use (\ref{eq:pfaffian-gene}), we get
\begin{eqnarray}\label{eq:exp-f-t}
\e^{f_t(X)}&=&\e^{-t^2\|x\|^2}\left(\sum_{I}\Pf(\hbox{$\frac{X_I}{2}$})\, e_I \right)
\left(\sum_J (-1)^{\frac{|J|(|J|-1)}{2}}dx_J\,e_J\, t^{|J|}\right) \nonumber\\
&=& \e^{-t^2\|x\|^2}
\sum_{I,J} \epsilon(I,J)(-1)^{\frac{|J|(|J|-1)}{2}}\Pf(\hbox{$\frac{X_I}{2}$})\,dx_J\,e_{I\cup J}\,
t^{|J|},
\end{eqnarray}
where the $I,J$ run over the subsets of $\{1,2,\ldots,d\}$ which are disjoint. If we take $t=1$ in
(\ref{eq:exp-f-t}), we see that
$$
\bere\left(\e^{f_1(X)}\right)=\e^{-\|x\|^2}\sum_{I} \epsilon(I,I')(-1)^{\frac{|I'|(|I'|-1)}{2}}
\Pf(\hbox{$\frac{X_I}{2}$})\,dx_{I'}.
$$
The third point then follows since
$(-1)^{\frac{|I'|(|I'|-1)}{2}}=(-1)^{\frac{d(d-1)}{2}}(-1)^{\frac{| I |}{2}}$.

Equation (\ref{eq:exp-f-t}) gives also
\begin{eqnarray*}
\eta_\wedge^t(X)&=&-\e^{-t^2\|x\|^2}\bere\left((\sum_k x_ke_k)
\Big(\sum_{I,J} \epsilon(I,J)(-1)^{\frac{|J|(|J|-1)}{2}}
\Pf(\hbox{$\frac{X_I}{2}$})\,dx_J\,e_{I\cup J}\, t^{|J|}\Big)\right)\\
&=&-\e^{-t^2\|x\|^2}\sum_{k,I,J} \epsilon(I,J)\epsilon(\{k\},I\cup J)
(-1)^{\frac{|J|(|J|+1)}{2}}\Pf(\hbox{$\frac{X_I}{2}$})\,x_k\,dx_J\, t^{|J|},
\end{eqnarray*}
where the sum runs over the partitions $\{k\}\cup I\cup J=\{1,\ldots,d\}$.
If we integrate the last equality between $0$ and infinity, and
use the formula $\int_0^{\infty} e^{-t^2} t^adt=\frac{1}{2}
\Gamma(\frac{a+1}{2})$, we get the first point.

The second point follows from the definition.
\end{proof}

\bigskip

Let us give the formulae for $\tur,\tuc,\tumq$  for small
dimensions. Here $f$ is a compactly supported function on $\Rbb$,
identically equal to $1$ in a neighborhood of $0$.

\begin{exam}
 If $\dim V=1$, then ${\rm SO}(V)$ is reduced to the
 identity. We have
$$
\beta_{\wedge}=-\frac{1}{2}\pi^{1/2}\frac{x}{|x|}
 = -\frac{1}{2}\pi^{1/2}\sign(x),
$$
so that
\begin{eqnarray*}
\tur(V)&=&\Big(0,-\frac{1}{2}\sign(x)\Big),\\
\tuc(V)&=&- f'(x^2) |x| dx,\\
\tumq(V)&=&\frac{1}{\pi^{1/2}}\e^{-\|x\|^2} dx.
\end{eqnarray*}

\end{exam}

\begin{exam}\label{exa:Vdim2}
Let $V=\Rbb e_1\oplus \Rbb e_2$ of dimension $2$. We have
$$
\beta_{\wedge}= \frac{x_1dx_2-x_2 dx_1}{2\|x\|^2}
$$
so that
\begin{eqnarray*}
\tur(V)(X)&=&\frac{-1}{\pi}\Big(\Pf(\hbox{$\frac{X}{2}$}),\beta_\wedge\Big),\\
\tuc(V)(X)&=& -\frac{1}{\pi}\Big(f(\|x\|^2) \Pf(\hbox{$\frac{X}{2}$}) +
  f'(\|x\|^2) dx_1\wedge dx_2\Big),\\
\tumq(V)(X)&=&\frac{1}{\pi}\e^{-\|x\|^2}(-\Pf(\hbox{$\frac{X}{2}$})+ dx_1\wedge dx_2).
\end{eqnarray*}

\end{exam}

\begin{exam} Let $V=\Rbb e_1\oplus \Rbb e_2\oplus \Rbb e_3$ of dimension $3$.
We have
%
\begin{eqnarray*}
&&\beta_{\wedge}(X)=-\frac{\pi^{1/2}}{4\|x\|}\Big(x_1 \langle X e_2,e_3
 \rangle + x_2 \langle X e_3,e_1\rangle +x_3 \langle X e_1,e_2\rangle\Big) \\
&&\qquad\qquad+ \frac{\pi^{1/2}}{4\|x\|^3}\Big(x_1 dx_2\wedge dx_3+ x_2 dx_3\wedge dx_1
 +x_3 dx_1\wedge dx_2\Big)
\end{eqnarray*}
 so that
$$
\tur(V)(X)=\frac{-1}{\pi^{3/2}}(0,\beta_\wedge(X)).
$$

The equivariant form  $\tuc(V)(X)$ is equal to
$$
\frac{1}{4\pi} \frac{f'(\|x\|^2)}{\|x\|}
\Big( (x_1\langle X e_2,e_3\rangle + x_2 \langle X e_3,e_1\rangle +
x_3 \langle X e_1,e_2 \rangle)(d\|x\|^2) -
2dx_1\wedge dx_2\wedge dx_3\Big).
$$

The equivariant form  $\tumq(V)(X)$ is equal to
$$
-\frac{1}{2\pi^{3/2}}\e^{-\|x\|^2}\Big(
\langle X e_2,e_3\rangle dx_1 + \langle X
e_3,e_1\rangle dx_2+\langle X e_1,e_2\rangle dx_3
- 2 dx_1\wedge dx_2\wedge dx_3\Big).
$$

\end{exam}

\medskip

\subsection{Unicity of the Thom forms of a vector space}\label{subsec:unicity-Thom-vector}

The following theorem is well known in the non equivariant case.

\begin{theo}\label{theo:thomvectorspace}

$\bullet$ The relative class $\tur(V)$ is a free generator of\break
$\Hcal^{\pol}(\sogot(V),V,V\setminus \{0\})$  over
$\Ccal^{\pol}(\sogot(V))^{\So(V)}$.
%
%
Thus $\tur(V)$ is the unique class in $\Hcal^{\pol}(\sogot(V),V,
V\setminus \{0\})$ of integral $1$.

\medskip

$\bullet$ The  equivariant class $\tuc(V)$ is a free generator of
$\Hcal^{\pol}_c(\sogot(V),V)$ over
$\Ccal^{\pol}(\sogot(V))^{\So(V)}$.
Thus $\tuc(V)$ is the unique class in
$\Hcal^{\pol}_c(\sogot(V),V)$ of integral  $1$.


\medskip

$\bullet$ The Mathai-Quillen class $\tumq(V)$  is a free generator
of $\Hcal^{\pol}_{\dr}(\sogot(V),V)$ over
$\Ccal^{\pol}(\sogot(V))^{\So(V)}$. Thus $\tumq(V)$ is the unique class in
$\Hcal^{\pol}_{\dr}(\sogot(V),V)$ of integral $1$.


\end{theo}

The same theorem holds in the $\Ccal^{\infty}$-category.

\begin{proof}
We give the proof of the second  point first. Let $t_V(X)$ be a closed equivariant form
with compact support, so that $\int_V t_V=1$. We want to prove that any
closed equivariant form $\alpha(X)$ in $\Acal^{\pol}_c(\sogot(V),V)$ is proportional to
$t_V(X)$ in cohomology : $\alpha=Q\, t_V$ in $\Hcal^{\pol}_c(\sogot(V),V)$ with $Q\in
\Ccal^{\pol}(\sogot(V))^{\So(V)}$.

 Let $\Ex(v)=-v$ be the
transformation $-\Id_V$. We consider the space $V\times V$. The
transformation $g(t)(v_1, v_2)=(\cos(t) v_1+\sin (t)
v_2,-\sin(t) v_1+\cos(t) v_2)$  is a one-parameter
transformation of $V\times V$ which commutes with the diagonal
action of $\So(V)$ and preserves the pair $(V\times V, V\times
V\setminus\{(0,0)\})$. The transformation $g(0)$ is the identity,
while $g(\frac{\pi}{2})(v_1,v_2)=(v_2,-v_1)$. Let $S: V\times V\to
V\times V$ be the vector field on $V\times V$ associated to the
action of the one-parameter subgroup $g(t)$. Thus
$S_{v_1,v_2}=(-v_2,v_1)$.

We denote by $\pi_1,\pi_2$ the first and second projections of
$V\times V$ on $V$. Let $\alpha_1,\alpha_2$ be closed equivariant
forms in $\Acal^{\pol}_c(\sogot(V),V)$. The exterior product
$A_{12}=\pi_1^*\alpha_1\wedge  \pi_2^*\alpha_2$ belongs to
$\Acal^{\pol}_c(\sogot(V),V\times V)$.

Let us apply the transformation $g(t)$ to  $A_{12}$. Define
$$A(t)= g(t)^*(A_{12}).$$
Then $A(0)=\pi_1^*\alpha_1\wedge \pi_2^*\alpha_2$ while
$A(\frac{\pi}{2})=\pi_2^* \alpha_1\wedge \pi_1^*\Ex^*\alpha_2$.
If the equivariant form $\alpha_i$ are supported in the balls
$B(0,r_i)\subset V$, we see that the equivariant $A(t)$ is supported
in the ball $B(0,r_1+r_2)\subset V\times V$ for any $t\in \Rbb$.

We have $\frac{d}{dt}A(t)=\Lcal(S) A(t)=D(\iota(S) A(t))$ from
Cartan's relation (\ref{Cartan}). Thus the cohomology class of
$A(t)$ remains constant. We obtain that, for  homogeneous elements
$\alpha_1,\alpha_2$,
$$
\pi_1^*\alpha_1\wedge  \pi_2^*\alpha_2
=\pi_2^*\alpha_1\wedge \pi_1^*\Ex^*\alpha_2=
(-1)^{|\alpha_1||\alpha_2|} \pi_1^*\Ex^*\alpha_2\wedge
\pi_2^*\alpha_1
$$
in $\Hcal^{\pol}_c(\sogot(V),V\times V)$.

We now consider the first
projection $\pi_1: V\times V\to V$ and the corresponding integral over the fiber
$(\pi_1)_*$ (see Section \ref{sec:integral-fiber}). Note that the map
$(\pi_1)_*\circ \pi_2^*: \Hcal^{\pol}_c(\sogot(V),V)\to
\Ccal^{\pol}(\sogot(V))^{\So(V)}$ corresponds to the
integration map $\int_V$. Thus we obtain the relation
$$
 \left(\int_V \alpha_2\right) \alpha_1=(-1)^{|\alpha_1||\alpha_2|}
 \left(\int_V \alpha_1\right) \Ex^*\alpha_2
$$
in $\Hcal^{\pol}_c(\sogot(V),V)$.

Let us apply this relation to the couple $(\alpha,\Ex^*t_V)$. As
$\int_V \Ex^*t_V=(-1)^{\dim V}$, we obtain that $\alpha=
Q(X)t_V$, with $Q(X)=\epsilon \int_V\alpha(X)$.
 Thus, if $\alpha$ is not zero, we can assume that its degree is
of same parity that $\dim V$, and  we obtain the relation
\begin{equation}\label{eq:alpha-tV}
\alpha=\left(\int_V\alpha\right) t_V
\end{equation}
in $\Hcal^{\pol}_c(\sogot(V),V)$. This proves that $t_V$ is unique in cohomology, and is a generator
of $\Hcal^{\pol}_c(\sogot(V),V)$.

The third point is proved in exactly the same way.

Let us prove the first point. To prove the fact that $\tur(V)$ is
a generator of $\Hcal^{\pol}(\sogot(V),V,V\setminus \{0\})$, we
apply the same one parameter  transformation $g(t)$ which acts on
$V\times V$ and preserves the subset $\{(0,0)\}\subset V\times V$.
Here we use the product $\diamond$ from
$$
\Hcal^{\pol}(\sogot(V),V\times V,(V\setminus \{0\})\times V)
\times \Hcal^{\pol}(\sogot(V),V\times V,V\times (V\setminus \{0\}))
$$
into $\Hcal^{\pol}(\sogot(V),V\times V,V\times V\setminus \{(0,0)\})$.
For a couple $(a_1,a_2)$ of $D_{\rm rel}$-closed
elements in $\Acal^{\pol}(\sogot(V),V,V\setminus \{0\})$, we consider the
product
$$
\pi_1^*(a_1)\diamond  \pi_2^*(a_2)\in
\Acal^{\pol}(\sogot(V),V\times V,V\times V\setminus \{(0,0)\})
$$
and their transformations $g(t)^*(\pi_1^*a_1\diamond  \pi_2^*a_2)$ which are in the
same cohomology class. We need the
\begin{lemm} We have the following equality
$$
g(\hbox{$\frac{\pi}{2}$})^*\left(\pi_1^*\alpha_1\diamond \pi_2^*(\alpha_2)\right)=
(-1)^{|\alpha_1||\alpha_2|} \pi_1^*\Ex^*\alpha_2\diamond
\pi_2^* \alpha_1
$$
in $\Hcal^{\pol}(\sogot(V),V\times V,V\times V\setminus \{(0,0)\})$.
\end{lemm}
\begin{proof} We consider the covering $V\times V\setminus \{(0,0)\}=U_1\cup U_2$
where $U_1:= (V\setminus \{0\})\times V$ and
$U_2=V\times (V\setminus \{0\})$. Let $\Phi=(\Phi_1,\Phi_2)$ a partition of unity
subordinate to this covering: the functions $\Phi_k$ are supposed ${\rm SO}(V)$-invariant.
We have also the group of symmetry generated by $\theta:=g(\hbox{$\frac{\pi}{2}$})$.
We have $\theta(U_1)=U_2$, $\theta(U_2)=U_1$ and $\theta^2(x,y)=(-x,-y)$.
We can suppose that the functions $\Phi_k$ are invariant under $\theta^2$, and that
$\theta^*(\Phi_1)=\Phi_2$.

Let $a_k:=(\alpha_k,\beta_k)\in \Acal^{\pol}(\sogot(V),V,V\setminus \{0\})$, $k=1,2$. Recall that
\break $\pi_1^*a_1\diamond_\Phi  \pi_2^*a_2$ is equal to
$(\pi_1^*\alpha_1\wedge \pi_2^*\alpha_2,\beta_\Phi(\pi_1^*a_1,\pi_2^*a_2)) $ with
\begin{eqnarray*}
\lefteqn{\beta_\Phi(\pi_1^*a_1,\pi_2^*a_2)= } \\
 & &\Phi_1\pi_1^*\beta_1\wedge
\pi_2^*\alpha_2+(-1)^{|a_1|}\pi_1^*\alpha_1\wedge \Phi_2 \pi_2^*\beta_2
 -(-1)^{|a_1|}d\Phi_1\wedge \pi_1^*\beta_1\wedge \pi_2^*\beta_2.
\end{eqnarray*}
Then $\theta^*(\pi_1^*a_1\diamond_\Phi  \pi_2^*a_2)$ is equal
to $(\theta^*(\pi_1^*\alpha_1\wedge \pi_2^*\alpha_2),\theta^*(\beta_\Phi(a_1,a_2)))$.
We know already that $\theta^*(\pi_1^*\alpha_1\wedge \pi_2^*\alpha_2)$ is equal to
$(-1)^{|a_1||a_2|} \pi_1^*\Ex^*\alpha_2\wedge \pi_2^*\alpha_1$. Let us compute
\begin{eqnarray*}
\lefteqn{\theta^*(\beta_\Phi(\pi_1^*a_1,\pi_2^*a_2))}\\
&=&\theta^*\left(\Phi_1\pi_1^*\beta_1\wedge \pi_2^*\alpha_2+(-1)^{|a_1|}\pi_1^*\alpha_1\wedge
\Phi_2 \pi_2^*\beta_2-(-1)^{|a_1|}d\Phi_1\wedge \pi_1^*\beta_1\wedge \pi_1^*\beta_2\right)\\
&=&\Phi_2 \pi_2^*\beta_1\wedge \pi_1^*\Ex^*\alpha_2+(-1)^{|a_1|}\pi_2^*\alpha_1\wedge
\Phi_1 \pi_1^*\Ex^*\beta_2\\
&& \qquad -(-1)^{|a_1|}d\Phi_2\wedge \pi_2^*\beta_1\wedge \pi_1^*\Ex^*\beta_2\\
&=&(-1)^{|a_1||a_2|} \beta_\Phi(\pi_1^*\Ex^*a_2,\pi_2^*a_1).
\end{eqnarray*}
So the elements $\theta^*(\pi_1^*a_1\diamond_\Phi  \pi_2^*a_2)$ and $(-1)^{|a_1||a_2|}
\pi_1^*\Ex^*a_2\,\diamond_\Phi  \pi_2^*a_1$ coincide in
$\Acal^{\pol}(\sogot(V),V\times V,V\times V\setminus \{(0,0)\})$.
\end{proof}

\bigskip

We have proved that
$$
\pi_1^*a_1\diamond \pi_2^*a_2=(-1)^{|a_1||a_2|}\pi_1^*\Ex^*a_2\diamond\pi_2^*a_1
$$
holds in $\Hcal^{\pol}(\sogot(V),V\times V,V\times V\setminus \{(0,0)\})$ for
any couple $(a_1,a_2)$ in \break $\Hcal^{\pol}(\sogot(V),V,V\setminus \{0\})$.
Then we consider the integral over the fiber of $\pi_1$ and we apply
Proposition \ref{prop:intfib}. The rest of the proof is the same.

\end{proof}

\bigskip

\begin{rem}
At the level of equivariant forms, Equality (\ref{eq:alpha-tV}) can be made more specific as follows. We have
the following equality in $\Acal^{\pol}_c(\sogot(V),V)$:
\begin{equation}\label{eq:alpha-tV-precis}
\alpha(X)=\left(\int_V\alpha(X)\right) t_V(X)+ D(\delta)(X)
\end{equation}
where $\delta=(-1)^{{\rm dim}(V)+1}(\pi_1)_*\left(\int_0^1\iota(S)A(t)\right)$, and
$A(t)=g(t)^*(\pi_1^*\alpha\wedge \pi_2^*\Ex^*t_V)$. We have also a control
on the support of the equivariant form $\delta$. If $\alpha$ and $t_V$ are supported respectively
in the balls $B(0,r)$ and $B(0,\epsilon)$ of $V$, the form $\delta$ is supported in the ball
$B(0,r+\epsilon)$.
\end{rem}

\bigskip

Since the map $\p_c: \Hcal^{\pol}(\sogot(V), V, V\setminus \{0\})\to\Hcal^{\pol}_c(\sogot(V),V)$
sends the class $\tur(V)$ to the class  $\tuc(V)$,
Theorem \ref{theo:thomvectorspace} shows that $\p_c$ is an {\em isomorphism}.

We can specify this property as follows. We have $\p_c=j\circ \p_{\{0\}}$ where
$\p_{\{0\}}: \Hcal^{\pol}(\sogot(V), V, V\setminus \{0\})\to
\Hcal^{\pol}_{\{0\}}(\sogot(V), V)$ and
$$
j:\Hcal^{\pol}_{\{0\}}(\sogot(V), V)\to
\Hcal^{\pol}_c(\sogot(V),V).
$$
is a natural map (see Lemma \ref{lem:basic-equi-coho-support-pol}): let us
recall its definition. A class in $\Hcal^{\pol}_{\{0\}}(\sogot(V), V)$
is defined by  a collection $[\gamma_U]\in \Hcal^{\pol}_{U}(\sogot(V), V)$, where
$U$ runs over the open neighborhood of $\{0\}$, and such that
$r_{U,U'}[\gamma_{U'}]=[\gamma_U]$ for $U'\subset U$. The image of $([\gamma_U])_U$ by $j$
is the class defined by the closed the equivariant form $\gamma_U$ in
$\Hcal^{\pol}_{c}(\sogot(V), V)$, for
any {\em relatively compact} open neighborhood $U$.

\begin{theo}
$\bullet$  The  maps $j$ and $\p_{\{0\}}$ are isomorphisms.

$\bullet$ Similarly, the  maps $\p_{\{0\}}:\Hcal^{\infty}(\sogot(V), V, V\setminus \{0\})\to
\Hcal^{\infty}_{\{0\}}(\sogot(V), V)$ and $j:\Hcal^{\infty}_{\{0\}}(\sogot(V), V)\to
\Hcal^{\infty}_c(\sogot(V), V)$ are isomorphisms.
\end{theo}

\begin{proof} Since $\p_c=j\circ \p_{\{0\}}$ is an isomorphism, it is enough to prove that
$j$ is one to one. Let $([\gamma_U])_U$ an element in the kernel of $j$ : for any
{\em relatively compact} open neighborhood $U$ of $\{0\}$, we have $\gamma_U=0$ in
$\Hcal^{\pol}_c(\sogot(V), V)$. If we show that $\gamma_U=0$ in
$\Hcal^{\pol}_U(\sogot(V), V)$, it gives that  $([\gamma_U])_U=0$.
Let $U'=B(0,r)$ be a ball and $0<\epsilon<<1$ such that $\overline{B(0,r+\epsilon)}$
is included in $U$.

We use (\ref{eq:alpha-tV-precis}) with the closed equivariant form $\gamma_{U'}$
supported in the ball $B(0,r)$, and a Thom form $t_V$ supported in the ball
$B(0,\epsilon)$.  We have the following relation in $\Acal^{\pol}_{U}(\sogot(V), V)$
\begin{eqnarray*}
\gamma_{U'}(X)&=&\left(\int_V\gamma_{U'}(X)\right) t_V(X)+ D(\delta)(X)\\
&=& D(\delta)(X)
\end{eqnarray*}
since $\gamma_{U'}=\gamma_U=0$ in $\Hcal^{\pol}_c(\sogot(V), V)$.
This proves that $\gamma_U=r_{U,U'}(\gamma_{U'})=0$ in $\Hcal^{\pol}_U(\sogot(V), V)$.

\end{proof}

\subsection{Explicit relative Thom form of a vector bundle and Thom isomorphism}
\label{subsection:thomvectorbundle}

Let $M$ be a manifold. Let $p:=\Vcal\to M$ be a real oriented
Euclidean vector bundle over $M$ of rank $d$. In this section, we
will describe a generator over $\Hcal(M)$ of the  relative
cohomology space of the pair $(\Vcal,\Vcal\setminus M)$. We will
use Chern-Weil construction.

Recall the sub-space $\Acal_{\cf}(\Vcal)\subset
\Acal(\Vcal)$ of differential forms on $\Vcal$
which have a compact support in the fibers of $p:\Vcal\to M$. We
have also defined the sub-space $\Acal_{\dr}(\Vcal)$. The
integration over the fiber is well defined on the three spaces
$\Acal^{\pol}(\kgot,\Vcal,\Vcal\setminus M)$,
$\Acal^{\pol}_{\cf}(\kgot,\Vcal)$ and
$\Acal^{\pol}_{\dr}(\kgot,\Vcal)$ and take values in $\Acal(M)$.
A Thom form on $\Vcal$ will be
a closed element which integrates to the constant function $1$ on
$M$.


\bigskip

The bundle $\Vcal$ is associated to a principal bundle $P\to M$
with structure group $G=\So(V)$, where $V=\Rbb^d$. We denote
$\sogot(V)$ by $\ggot$.
 An element $y\in
P$ above $x\in M$ is by definition a map $y:V\to \Vcal_x$
conserving the inner product and the orientation. Thus $P$ is
equipped with an action of $\So(V)$: $g\cdot y=yg^{-1}: V\to
\Vcal_{x}$.

\begin{defi}
Let $\omega$ be a  connection one form on $P$, with curvature form
$\Omega$. The Euler form $\Eul(\Vcal,\omega)\in \Acal(M)$ of $\Vcal\to M$
is the closed differential form on $M$ defined by
$$
\Eul(\Vcal,\omega):=\Pf\left(-\frac{\Omega}{2\pi}\right).
$$
The class of $\Eul(\Vcal,\omega)$, which does not depend of $\omega$, is denoted
$\Eul(\Vcal)\in\Hcal^{d}(M)$.
\end{defi}

\begin{rem}
Since the pfaffian vanishes on $\so(V)$ when $\dim V$ is odd,
the Euler class $\Eul(\Vcal)\in\Hcal^{d}(M)$
is identically equal to $0$ when the rank of $\Vcal$ is odd.
\end{rem}

Consider the Chern-Weil map
$$
\phi_\omega^Z:\Acal^{\pol}(\ggot,Z) \longrightarrow \Acal( \Zcal)
$$
where the manifold $Z$ is the $\{ {\rm pt} \}$, $V$, or
$Z=V\setminus\{0\}$: hence the quotient manifold
$\Zcal=P\times_{G}Z$ are respectively $M$, $\Vcal$ and
$\Vcal\setminus M$. In order to simplify the notations we denote
all these maps by $\phi_\omega$. The map $\phi_\omega$ maps
$\Acal_c^{\pol}(\ggot,V)$ to $\Acal_{\cf}(\Vcal)$ and
$\Acal_{\dr}^{\pol}(\ggot,V)$ to $\Acal_{\dr}(\Vcal).$

Recall the equivariant differential forms defined in Theorem \ref{theo:thom-form} :
\begin{eqnarray*}
\tur(V)&\in& \Acal^{\pol}(\ggot,V,V\setminus \{0\}) ,\\
\tuc(V)&\in& \Acal_c^{\pol}(\ggot,V),\\
\tumq(V)&\in& \Acal_{\dr}^{\pol}(\ggot,V).
\end{eqnarray*}

We can take the image of these forms via the Chern-Weil map. We consider the transgression form
$\beta_{\Vcal,\omega}\in\Acal^{d-1}(\Vcal\setminus M)$ defined by
$\beta_{\Vcal,\omega}=\phi_\omega(\beta_\wedge)$.

\begin{defi}

$\bullet$ The relative form
$\tur(\Vcal,\omega)\in\Acal^d(\Vcal,\Vcal\setminus M)$ is the image of
$\tur(V)\in\Acal^{\pol}(\ggot,V,V\setminus \{0\})$ by the
Chern-Weil map $\phi_\omega$~:
$$
\tur(\Vcal,\omega)=\frac{1}{\epsilon_d}
\left(p^*\Pf(\hbox{$\frac{\Omega}{2}$}),\beta_{\Vcal,\omega}\right).
$$

$\bullet$  The form
$\tucf(\Vcal,\omega)\in\Acal_{\cf}^d(\Vcal)$ is defined as the image of
$\tuc(V)\in\Acal^{\pol}_c(\ggot,V)$ by the Chern-Weil map
$\phi_\omega$.

$\bullet$ The Mathai-Quillen form
$\tumq(\Vcal,\omega)\in\Acal_{\dr}^d(\Vcal)$ is defined as the image of
$\tumq(V)\in\Acal^{\pol}_{\dr}(\ggot,V)$ by the Chern-Weil map
$\phi_\omega$:


\end{defi}

\medskip

The form $\tur(\Vcal,\omega)$  is relatively closed, since
$\beta_{\Vcal}$ is defined on $\Vcal\setminus M$ and satisfies
$d\beta_{\Vcal,\omega}=p^*\Pf(\frac{\Omega}{2})$. The forms
$\tucf(\Vcal,\omega)$ and $\tumq(\Vcal,\omega)$ are closed de Rham differential
forms. We denote $\tur(\Vcal)$, $\tucf(\Vcal)$ and $\tumq(\Vcal)$ the corresponding
cohomology classes. Since
the map $\phi_\omega$ commutes with the integration over the
fiber, these cohomology classes  images are Thom classes.

We obtain the analog of Theorem \ref{theo:thomvectorspace}.

\begin{theo}\label{theo:thomvectorbundle}.

$\bullet$ The relative class $\tur(\Vcal)$ is a free generator of
$\Hcal(\Vcal,\Vcal\setminus M)$  over $\Hcal(M)$. Thus $\tur(\Vcal)$ is the unique class in
$\Hcal(\Vcal,\Vcal\setminus M)$ with  integral $1$ along
  the fiber. We say that $\tur(\Vcal)$ is the Thom class in
$\Hcal(\Vcal,\Vcal\setminus M)$.

\medskip

$\bullet$ The  class $\tucf(\Vcal)$ is a free generator of
$\Hcal_{\cf}(\Vcal)$ over $\Hcal(M)$. Thus $\tucf(\Vcal)$ is the
unique class  in $\Hcal_{\cf}(\Vcal)$ with  integral $1$ along the
fiber. We say that $\tucf(\Vcal)$ is the Thom class in
$\Hcal_{\cf}(\Vcal)$.

\medskip

$\bullet$ The Mathai-Quillen class $\tumq(\Vcal)$ is a free generator
of $\Hcal_{\dr}(\Vcal)$ over $\Hcal(M)$. Thus $\tumq(\Vcal)$ is
the unique class in $\Hcal_{\dr}(\Vcal)$ with  integral $1$ along
the fiber. We say that $\tumq(\Vcal)$ is the
Thom class in $\Hcal_{\dr}(\Vcal)$.

\end{theo}

\begin{proof}
The proof is the same than the proof of Theorem
\ref{theo:thomvectorspace}. We work on the sum of vector bundles
$\Vcal\oplus \Vcal$ over $M$ and we relate the identity on
$\Vcal\oplus \Vcal$ to the exchange of the two copies of $\Vcal$
via a one-parameter transformation group.
\end{proof}

\bigskip

From the explicit formulae for $\tucf,\tumq$, we obtain the
following proposition.

\begin{prop}
The restriction of the Thom classes $\tucf(\Vcal)$ or $\tumq(\Vcal)$  to $M$ is the
Euler class $\Eul(\Vcal)$.
\end{prop}

Let us give explicit expressions for the forms $\tur(\Vcal,\omega),
\tucf(\Vcal,\omega)$ as well as $\tumq(\Vcal,\omega)$. Let us fix an oriented orthonormal basis
$(e_1,\ldots,e_d)$ of $V$. We write
$\omega=\sum_{k,\ell}\omega_{ k\ell}E_{\ell}^k$ where $E_{\ell}^k$
is the  transformation of $V$ such that
$E_{\ell}^k(e_i)=\delta_{i,k}e_{\ell}$. Here $\omega_{k\ell}$ are $1$-forms on
$P$ and $\omega_{k\ell}=-\omega_{\ell k}$. The curvature $\Omega$
is  $\Omega=\sum_{k,\ell}\Omega_{k\ell}E_{\ell}^k$ where
$\Omega_{k\ell}=d\omega_{k\ell}+\sum_{j}\omega_{j\ell}\wedge\omega_{k
j}$.


The connection $\nabla$ on $\Vcal\to M$ induced by
the connection form $\omega$ is the operator $\nabla=d+\omega$
acting on $(\f(P)\otimes V)^G$ with values in $(\Acal^1(P)_{\hor} \otimes V)^G$.
The Chern-Weil map $\phi_\omega:\Acal^{\infty}(\ggot,V)\to (\Acal(P\times V)_{\hor})^G$ admit
a natural extension
$$
\phi_\omega^\wedge:(\f(\ggot)\otimes\Acal(V)\otimes \wedge V)^G\longrightarrow
(\Acal(P\times V)_{\hor} \otimes \wedge V)^G
$$
such that $\bere\circ \,\phi_\omega^\wedge=\phi_\omega\circ\bere$.

Let $f_t\in(\Ccal^{\pol}(\ggot)\otimes\Acal(V)\otimes \wedge V)^G$ be the map defined in
(\ref{eq-map-f-t-defi}). The element $f_t^\omega:=\phi_\omega^\wedge(f_t)$  is
defined by the equation
\begin{equation}
f_t^\omega= -t^2\|x\|^2 + t\sum_i \eta_i e_i + \frac{1}{2}\sum_{k<l}
\Omega_{k\ell} e_k\wedge e_\ell.
\end{equation}
where $\eta_i= dx_i +\sum_{\ell}\omega_{\ell i}x_{\ell}$
is a horizontal $1$-form on $P\times V$. If $I=\{i_1<\cdots<i_k\}$ is a subset
of $[1,2,\ldots,d]$, we denote by $\eta_I$ the product $\eta_{i_1}\wedge\cdots\wedge\eta_{i_k}$.

Let $\eta_{\wedge}^t\in(\Ccal^{\pol}(\ggot)\otimes\Acal(V))^G$ be the transgression forms defined in
(\ref{eq:C-eta-wedge}). The element $\phi_\omega(\eta_{\wedge}^t)\in  \Acal(P\times V)_{\bas}$
satisfies
\begin{eqnarray*}
\phi_\omega(\eta_{\wedge}^t)&=&
- \phi_\omega\circ\bere\Big((\sum_k x_k e_k)\e^{f_t(X)}\Big)\\
&=& -\bere\circ\,\phi_\omega^\wedge\Big((\sum_k x_k e_k)\e^{f_t(X)}\Big)\\
&=& -\bere\Big((\sum_k x_k e_k)\e^{f_t^\omega}\Big).
\end{eqnarray*}

We use the same notations as in Proposition \ref{prop:explicit}.

\begin{prop}\label{prop:explicit-2}

$\bullet$  The form $\beta_{\Vcal,\omega}=\phi_\omega(\beta_\wedge)$ is defined
in $\Acal(P\times (V\setminus\{0\}))_{\bas}$ by the relation
$$
\beta_{\Vcal,\omega}=-\int_0^{\infty} \bere\Big((\sum x_k e_k)\e^{-t^2\|x\|^2+
t\sum_k \eta_k e_k+\frac{1}{2}\sum_{k<\ell}\Omega_{k\ell}e_k\wedge e_\ell}\Big ).
$$
More explicitly,
$$
\beta_{\Vcal,\omega}=
\sum_{k,I,J} \gamma_{(k,I,J)}
 \Pf\left(\hbox{$\frac{\Omega_I}{2}$}\right)\frac{x_k \eta_J}{\|x\|^{|J|+1}}\cdot
$$

$\bullet$ The form $\tucf(\Vcal,\omega)$ is defined in $\Acal(P\times V)_{\bas}$ by the
relation
$$
\tucf(\Vcal,\omega)=
\frac{1}{\epsilon_d}(\chi \Pf(\hbox{$\frac{\Omega}{2}$})+d\chi \beta_{\Vcal,\omega}).
$$
where $\chi$ is a $\So(V)$-invariant function on $V$,
identically equal to $1$ in a neighborhood of $0$.

$\bullet$ The form $\tumq(\Vcal,\omega)$ is defined in $\Acal(P\times V)_{\bas}$ by the
relation
$$
\tumq(\Vcal,\omega)=\frac{1}{\epsilon_d}\bere\Big (\e^{-\|x\|^2+
\sum_k\eta_k  e_k+\frac{1}{2}\sum_{k<\ell}
\Omega_{k\ell} e_k\wedge e_\ell}\Big).
$$
More explicitly
$$
\tumq(\Vcal,\omega)=\frac{1}{(\pi)^{d/2}}\e^{-\|x\|^2} \sum_{I\ {\rm even}} (-1)^{\frac{| I |}{2}}\epsilon(I,I')
\Pf(\hbox{$\frac{\Omega_I}{2}$}) \eta_{I'}.
$$

\end{prop}

This gives an explicit expression of $\beta_{\Vcal,\omega}$ and thus of $\tur(\Vcal,\omega)$ in
functions of the variables $x_j$, $dx_j$ on $V$, the connection
one forms $\omega_{k\ell}$ and  the curvature forms
$\Omega_{k\ell}$ of an Euclidean  connection on $M$.

The form $\tumq(\Vcal,\omega)$ is the representative with Gaussian look
of the Thom form constructed by Mathai-Quillen
\cite{Mathai-Quillen}.

We also obtain:

\begin{theo}

$\bullet$ The map  $\Hcal(\Vcal, \Vcal\setminus M)\to \Hcal_M(\Vcal) $ is an
isomorphism.

$\bullet$ If $M$ is compact, the map $\Hcal(\Vcal, \Vcal\setminus M)\to
\Hcal_c(\Vcal) $ is an isomorphism.
\end{theo}

\subsection{Explicit equivariant relative Thom form of a vector bundle and Thom isomorphism}

We can repeat the construction above in the equivariant case.
Assume now that  $\Vcal$ is a $K$-equivariant vector bundle. Then
$\Vcal$  is associated to a $K$-equivariant principal bundle $P\to
M$ with structure group $G=\So(V)$. The principal bundle $P$
 is provided with an action of $K\times G$. If $y:V\to \Vcal_x$ is
 an orthonormal frame, then $(g,k)\cdot y=kyg^{-1}$ is a frame
 over $kx$.

 Let $\omega$ be a $K$-invariant
connection one form on $P$, with curvature form $\Omega$.
 For $Y\in\kgot$, we denote by $\mu(Y)=-\iota(VY)\omega\,
\in\f(P)\otimes\ggot$ the moment of $Y\in \kgot$. The equivariant
curvature form is $\Omega(Y)=\Omega+\mu(Y),\, Y\in\kgot$.

Consider the Chern-Weil maps
$$
\phi_\omega:\Acal^{\pol}(\ggot,Z) \longrightarrow
\Acal^{\pol}(\kgot, \Zcal)
$$
where the manifold $Z$ is the $\{{\rm pt}\}$, $V$, or
$Z=V\setminus\{0\}$. See (\ref{eq:chern-weil-K}) for the definition.

\begin{defi}
Let $\omega$ be a  connection one form on $P$, with curvature form
$\Omega$. The equivariant Euler form $\Eul(\Vcal,\omega)\in \Acal^{\pol}(\kgot,M)$ of $\Vcal\to M$
is the closed equivariant form on $M$ defined by
$$
\Eul(\Vcal,\omega)(Y):=\Pf\left(-\frac{\Omega(Y)}{2\pi}\right),\ Y\in\kgot.
$$
In other words, $\Eul(\Vcal,\omega)$ is the image by the Chern-Weil map $\phi_\omega$ of the invariant
polynomial $X\mapsto\Pf(-\frac{X}{2\pi})$. The class of $\Eul(\Vcal,\omega)$, which
does not depend of $\omega$, is denoted
$\Eul(\Vcal)\in\Hcal^{\pol}(\kgot,M)$.
\end{defi}

As integration over the fiber commutes with the Chern-Weil
construction,  the image by $\phi_\omega$ of Thom classes are Thom
classes. We define the transgression form $\beta_{\Vcal,\omega}(Y)\in\Acal^{\pol}(\kgot,\Vcal\setminus M)$ by
$$
\beta_{\Vcal,\omega}(Y)=(\phi_\omega(\beta_\wedge))(Y).
$$

\begin{defi}
$\bullet$ The relative equivariant form
$\tur(\Vcal,\omega)\in\Acal^{\pol}(\kgot,\Vcal,\Vcal\setminus M)$ is the
image of $\tur(V)\in\Acal^{\pol}(\ggot,V,V\setminus \{0\})$ by the
Chern-Weil map $\phi_\omega$. More explicitly, for $Y\in \kgot$
$$
\tur(\Vcal,\omega)(Y)=
\frac{1}{\epsilon_d}\left(p^*\Pf(\hbox{$\frac{\Omega(Y)}{2}$}),\beta_{\Vcal,\omega}(Y)\right).
$$

$\bullet$  We define the equivariant form
$\tucf(\Vcal,\omega)(Y)\in\Acal_{\cf}^{\pol}(\kgot,\Vcal)$ as the image of
$\tuc(V)\in\Acal^{\pol}_c(\ggot,V)$ by the Chern-Weil map
$\phi_\omega$:
$$
\tucf(\Vcal,\omega)(Y):=(\phi_{\omega}\tuc(V))(Y).
$$

$\bullet$ We define the Mathai-Quillen form
$\tumq(\Vcal,\omega)(Y)\in\Acal_{\dr}^{\pol}(\kgot,\Vcal)$ as the image
of $\tumq(V)\in\Acal^{\pol}_{\dr}(\ggot,V)$ by the Chern-Weil map
$\phi_\omega$:
$$
\tumq(\Vcal,\omega)(Y):=(\phi_{\omega}\tumq(V))(Y).
$$
\end{defi}

The equivariant form $\tur(\Vcal,\omega)$  is relatively closed, since
$\beta_{\Vcal,\omega}$ is defined on $\Vcal\setminus M$ and satisfies
$D(\beta_{\Vcal,\omega})=p^*\Pf(\frac{\Omega}{2})$. The equivariant forms
$\tucf(\Vcal,\omega)$ and $\tumq(\Vcal,\omega)$ are closed equivariant differential
forms. We denote $\tur(\Vcal)$, $\tucf(\Vcal)$ and $\tumq(\Vcal)$ the corresponding
equivariant cohomology classes. Since
the map $\phi_\omega$ commutes with the integration over the
fiber, these cohomology classes  images are equivariant Thom classes.

With the same proof as the non equivariant case, we obtain the
following Theorem.

\begin{theo}\label{theo:thomvectorbundleequi}

$\bullet$ The relative class $\tur(\Vcal)(Y)$ is a free generator of \break
$\Hcal_{\dr}^{\pol}(\kgot,\Vcal,\Vcal\setminus M)$  over
$\Hcal^{\pol}(\kgot,M)$. Thus $\tur(\Vcal)$ is the unique class in
$\Hcal_{\rel}^{\pol}(\kgot,\Vcal,\Vcal\setminus M)$ with integral $1$ along the fiber. We say that
$\tur(\Vcal)$ is the equivariant Thom form in
$\Hcal_{\dr}^{\pol}(\kgot,\Vcal,\Vcal\setminus M)$.

\bigskip

 $\bullet$ The  equivariant class $\tucf(\Vcal)(Y)$ is a free generator of
$\Hcal^{\pol}_{\cf}(\kgot,\Vcal)$ over $\Hcal^{\pol}(\kgot, M)$.
Thus $\tucf(\Vcal)$ is the unique class in
$\Hcal_{\cf}^{\pol}(\kgot,\Vcal)$ with integral $1$ along the
fiber. We say that $\tucf(\Vcal)$ is the
equivariant Thom form in $\Hcal_{\cf}^{\pol}(\kgot,\Vcal)$.

\bigskip

$\bullet$ The Mathai-Quillen class $\tumq(\Vcal)(Y)$ is a free
generator of $\Hcal_{\dr}^{\pol}(\kgot,\Vcal)$ over
$\Hcal^{\pol}(\kgot,M)$. Thus $\tumq(\Vcal)$ is the unique
class in $\Hcal_{\dr}(\kgot,\Vcal)$  with integral $1$
along the fiber. We say that
$\tumq(\Vcal)$ is the equivariant Thom form in $\Hcal_{\dr}(\kgot,\Vcal)$.

\end{theo}

\begin{proof}
The proof is the same than the proof of Theorem
\ref{theo:thomvectorbundle}.

\end{proof}

\bigskip

From the explicit formulae for $\tucf(\Vcal,\omega)(Y),\tumq(\Vcal,\omega)(Y)$,
we obtain the following proposition.

\begin{prop}
The restriction of the equivariant Thom class $\tucf(\Vcal)$ (or
$\tumq(\Vcal)$) to $M$ is the equivariant Euler class $\Eul(\Vcal)$.
\end{prop}

Let us give explicit expressions for the equivariant forms
$\tucf(\Vcal,\omega)(Y)$ as well as $\tumq(\Vcal,\omega)(Y)$ : we will express them as
$K$-equivariant map from $\kgot$ into $\Acal(P\times V)_{\bas}$. We use the same
notations as in Proposition \ref{prop:explicit-2}.

We denote by $\mu(Y)$ the moment of $Y\in \kgot$ with respect to
the connection $\nabla=d+\omega$. By definition
$\mu(Y)=\Lcal(Y)-\nabla_{VY}=-\iota(VY)\omega$ where $VY$ is the vector field on $P$ associated to
$Y\in\kgot$. Thus $\mu(Y)$, viewed as a $\So(V)$-invariant map from $P$ into $\so(V)$, satisfies
$$
\mu(Y)= - \sum_{k,\ell} \langle\omega_{k\ell},VY\rangle E_{\ell}^{k}.
$$
For the equivariant curvature we have $\Omega(Y)=\sum_{k,\ell}\Omega_{k\ell}(Y)E_{\ell}^k$ with
$$
\Omega_{k\ell}(Y)=d\omega_{k\ell}+\sum_{j}\omega_{j\ell}\wedge\omega_{kj}
 -\langle\omega_{k\ell},VY\rangle.
$$

As usual in equivariant cohomology, formulae for equivariant
classes are obtained from the non equivariant case by replacing
the curvature $\Omega$ by the equivariant curvature.

\begin{prop}

 $\bullet$  The transgression form $\beta_{\Vcal,\omega}(Y)$ is defined by the following formula :
 for $Y\in\kgot$ we have
$$
\beta_{\Vcal,\omega}(Y)=-\int_0^{\infty} \bere\Big((\sum x_k e_k)\e^{-t^2\|x\|^2+
t\sum_i \eta_i e_i+\frac{1}{2}\sum_{k<\ell} \Omega_{k,\ell}(Y)e_k\wedge e_\ell}\Big ).
$$
More explicitly,
$$
\beta_{\Vcal,\omega}(Y)=
\sum_{k,I,J} \gamma_{(k,I,J)}
 \Pf\left(\hbox{$\frac{\Omega_I(Y)}{2}$}\right)\frac{x_k \eta_J}{\|x\|^{|J|+1}}\cdot
$$


$\bullet$ Let $\chi$ be a $\So(V)$ invariant function on $V$,
identically equal to $1$ in a neighborhood of $0$.  We have thus
$$
\tucf(\Vcal,\omega)(Y)=\frac{1}{\epsilon_d}
\left(\chi \Pf(\hbox{$\frac{\Omega(Y)}{2}$})+d\chi \beta_{\Vcal,\omega}(Y)\right),\quad
Y\in\kgot.
$$

$\bullet$ We have
$$
\tumq(\Vcal,\omega)(Y)=\frac{1}{\epsilon_d}\bere\Big (\e^{-\|x\|^2+
\sum_k\eta_k  e_k+\frac{1}{2}\sum_{k<\ell}\Omega_{k\ell}(Y)e_k\wedge e_\ell}\Big),\quad
Y\in\kgot.
$$
More explicitly
$$
\tumq(\Vcal,\omega)(Y)=\frac{1}{(\pi)^{d/2}}\e^{-\|x\|^2} \sum_{I\ {\rm even}}
(-1)^{\frac{| I |}{2}}\epsilon(I,I')
\Pf(\hbox{$\frac{\Omega_I(Y)}{2}$})\, \eta_{I'},\quad
Y\in\kgot.
$$

\end{prop}

\medskip

We obtain similarly:

\begin{theo}
$\bullet$ The map  $\p_M: \Hcal^{\pol}(\kgot,\Vcal, \Vcal\setminus M)\to
\Hcal^{\pol}_M(\kgot,\Vcal) $ is an isomorphism.

$\bullet$ If $M$ is compact, the map $\p_c:\Hcal^{\pol}(\kgot,\Vcal, \Vcal\setminus
M)\to \Hcal^{\pol}_c(\kgot,\Vcal) $ is an isomorphism
\end{theo}

\section{The relative Chern character}\label{subsec:chern-character}

Let $N$ be a manifold equipped with an action of a compact Lie
group $K$.

\subsection{Quillen's Chern form of a super-connection}

For an introduction to the Quillen's notion of super-connection, see
\cite{B-G-V}.

Let $\Fcal$ be a complex vector bundle on $N$. Let $\nabla$ be a
connection on $\Fcal$. The curvature $\nabla^2$ of $\nabla$ is a
$\End(\Fcal)$-valued two-form on $N$. Recall that the
(non-normalized)  Chern character of $\Fcal$ is the de Rham
cohomology class $\ch(\Fcal)\in \Hcal(N)^+$ of the closed
differential form $\ch(\nabla):=\tr(\e^{\nabla^2})$.

Suppose now that $\Fcal$ is a $K$-equivariant vector bundle, and
suppose that $\nabla$ is $K$-invariant. For any $X\in\kgot$, we
consider $\mu^{\nabla}(X)=\Lcal(X)-\nabla_{VX}$ which is an
$\End(\Fcal)$-valued function on $N$: here $\Lcal(X)$ is the Lie
derivative of $X$, acting on $\Acal(N,\Fcal)$, and $\nabla_{VX}$
is  equal to $\iota(VX)\nabla$. Then $\ch(\nabla)(X)=
\tr(\e^{\nabla^2+\mu^{\nabla}(X)})$ is a closed equivariant form
on $N$: its class $\ch(\Fcal)\in\Hcal^\infty(\kgot,N)$ is the
equivariant Chern character of $\Fcal$.

More generally, let $\Ecal=\Ecal^+\oplus\Ecal^-$ be an equivariant
$\Zbb_2$-graded complex vector bundle on $N$. We denote by
$\Acal(N,\End(\Ecal))$ the algebra of $\End(\Ecal)$-valued
differential forms on $N$. Taking in account the $\Zbb_2$-grading of
$\End(\Ecal)$, the algebra $\Acal(N,\End(\Ecal))$ is a
$\Zbb_2$-graded algebra. The super-trace on $\End(\Ecal)$ extends to
a map $\str:\Acal(N,\End(\Ecal))\to \Acal(N)$.

Let $\A$ be a $K$-invariant super-connection  on $\Ecal$, and
$\F=\A^2$ its curvature, an element of $\Acal(N,\End(\Ecal))^+$.
Recall that, for $X\in \kgot$, the \emph{moment} of $\A$ is the
equivariant map
\begin{equation}
  \label{eq:moment-phi}
  \mu^{\A}:\kgot \longrightarrow \Acal(N,\End(\Ecal))^+
\end{equation}
defined by the relation $\mu^{\A}(X)=\Lcal(X)-[\iota(VX),\A]$. We
define the \emph{equivariant curvature} of $\A$ by
\begin{equation}\label{eq:equi-curvature}
    \F(X)=\A^2+\mu^{\A}(X),\quad X\in\kgot.
\end{equation}

We usually denote simply by $\F$ the equivariant curvature, keeping
in mind that in the equivariant case, $\F$ is a function from
$\kgot$ to $\Acal(N,\End(\Ecal))^+$.

\begin{defi}\label{defi:chernequibete}(Quillen)
  The equivariant Chern character of  $(\Ecal,\A)$
  is the equivariant differential form on
  $N$ defined by
   $\ch(\A)=\str\left(\e^{\F}\right)$
   (e.g. $\ch(\A)(X)=\str\left(\e^{\F(X)}\right)$).
\end{defi}
The form $\ch(\A)$ is equivariantly closed. We will use the
following transgression formulas:

\begin{prop}\label{transequi}

$\bullet$ Let $\A_t$, for $t\in \Rbb$, be  a one-parameter family
of $K$-invariant super-connections on $\Ecal$, and let
$\frac{d}{dt}\A_t\in \Acal(N,\End(\Ecal))^-$. Let $\F_t$ be the
equivariant curvature of $\A_t$. Then one has
\begin{equation}\label{transgression-equi}
\frac{d}{dt}\ch(\A_t)=D\left(\str\Big((\frac{d}{dt}\A_t)
\e^{\F_t}\Big)\right).
\end{equation}

$\bullet$ Let $\A(s,t)$ be a two-parameters family of
$K$-invariant super-connections. Here $s,t \in \Rbb$. We denote by
$\F(s,t)$ the equivariant curvature of $\A(s,t)$. Then:
\begin{eqnarray*}
\lefteqn{\frac{d}{ds}\str\Big((\frac{d}{dt}\A(s,t))\,
\e^{\F(s,t)}\Big)-
\frac{d}{dt}\str\Big((\frac{d}{ds}\A(s,t))\, \e^{\F(s,t)}\Big)} \nonumber\\
& &= D\left(\int_{0}^1\str\Big((\frac{d}{ds}\A(s,t)) e^{u
\F(s,t)}(\frac{d}{dt}\A(s,t))\,\e^{(1-u)\F(s,t)}\Big)du\right).
\end{eqnarray*}

\end{prop}

These formulae are the consequences of the two following identities.
See \cite{B-G-V}, chapter 7.

$\bullet$ If we denote by  $\A_{\kgot}(X)$ the operator
$\A-\iota(VX)$ on $\Acal(N,\End(\Ecal))$, then we have the relation:
$$\F(X)=\A_{\kgot}(X)^2+\Lcal(X).$$

$\bullet$ If $\alpha$ is an equivariant map from $\kgot$ to
$\Acal(N,\End(\Ecal))$, then one has

$$D\left(\str(\alpha(X))\right)=\str[\A_{\kgot}(X),\alpha(X)].$$

\medskip

 Then, using the invariance by $K$ of all terms involved,
  the proof of  Proposition \ref{transequi} is entirely similar to the
proof of Proposition 3.1 in \cite{pep-vergne1}.

In particular, the cohomology class defined by $\ch(\A)$  in
$\Hcal^\infty(\kgot,N)$ is independent of the choice of the
super-connection $\A$ on $\Ecal$. By definition, this is the
equivariant Chern character $\ch(\Ecal)$ of $\Ecal$. By choosing
$\A=\nabla^+\oplus\nabla^-$ where $\nabla^{\pm}$ are connections on
$\Ecal^{\pm}$, this class is just $\ch(\Ecal^+)-\ch(\Ecal^-)$.
However, different choices of $\A$ define very different looking
representatives of $\ch(\Ecal)$.

\subsection{The relative Chern character of a morphism}
\label{subsec:chern-character-rel}

Let $\Ecal=\Ecal^+\oplus \Ecal^-$ be an equivariant $\Zbb_2$-graded
complex vector bundle on $N$ and $\sigma: \Ecal^+ \to \Ecal^-$ be a
smooth morphism which commutes with the action of $K$. At each point
$n\in N$, $\sigma(n): \Ecal^+_n\to \Ecal^-_n$ is a linear map.
 The \emph{support} of $\sigma$ is the $K$-invariant closed subset
$$
\supp(\sigma)=\{n\in N\mid \sigma(n)\ {\rm is\ not \
invertible}\}.
$$

\begin{defi}\label{elliptic}
The morphism $\sigma$ is elliptic if $\,\supp(\sigma)$ is
compact.
\end{defi}

Recall that the data $(\Ecal^+,\Ecal^-,\sigma)$ defines an element
of the equivariant  $K$-theory $\KK_K(N)$ of $N$ when $\sigma$ is
elliptic. In the following, {\bf we do not} assume $\sigma$
elliptic. Inspired by  Quillen \cite{Quillen85}, we construct  a
cohomology class $\ch_{\rm rel}(\sigma)$ in
$\Hcal^\infty(\kgot,N,N\setminus \supp(\sigma))$.

The definition will involve several choices. We choose invariant Hermitian
structures on $\Ecal^{\pm}$ and  an invariant super-connection $\A$ on $\Ecal$
{\bf without $0$ exterior degree term}. This means that
$\A=\sum_{j\geq 1}\A_{[j]}$, where $\A_{[1]}$ is a connection on
the bundle $\Ecal$ which preserves the grading, and for $j\geq 2$,
the operator $\A_{[j]}$ is given by the action of a differential
form $\omega_{[j]}\in \Acal^j(N,\End(\Ecal))$ on $\Acal(N,\Ecal)$.
 Furthermore, $\omega_{[j]}$ lies in $\Acal^j(N,\End(\Ecal)^-)$ if $j$ is
even, and in $\Acal^j(N,\End(\Ecal)^+)$ if $j$ is odd.

We define, with the help of the invariant Hermitian metric on $\Ecal^\pm$, the dual of the morphism
$\sigma$ as an equivariant morphism $\sigma^*:\Ecal^-\to\Ecal^+$.
Introduce the odd Hermitian endomorphism of $\Ecal$ defined by
\begin{equation}\label{eq:v-sigma}
v_\sigma=
\left(\begin{array}{cc}
0 & \sigma^*\\
\sigma & 0\\
\end{array}\right).
\end{equation}
Then $ v_\sigma^2=\left(\begin{array}{cc}
\sigma^*\sigma& 0\\
0 &\sigma\sigma^*\\
\end{array}\right)
$ is a non negative even Hermitian  endomorphism of $\Ecal$. The
support of $\sigma$ coincides with the set of elements $n\in N$
where the spectrum of  $v_\sigma^2(n)$ contains $0$.

\begin{defi}\label{defi:h-sigma}
We denote by $h_\sigma(n)\geq 0$ the smallest eigenvalue
of $v_\sigma^2(n)$.
\end{defi}



\medskip

 Consider the family of super-connections
$\A^{\sigma}(t)=\A+i t\, v_\sigma,\ t\in\Rbb$ on $\Ecal$. The
equivariant curvature of $\A^{\sigma}(t)$ is thus the map
\begin{equation}\label{eq:F-A-sigma}
    \F(\sigma,\A,t)(X)=-t^2 v_\sigma^2+it
[\A,v_\sigma]+\A^2+\mu^\A(X),\quad X\in\kgot
\end{equation}
with $0$ exterior degree term equal to $-t^2 v_\sigma^2+
\mu^\A_{[0]}(X)$. As the super-connection $\A$ do not have $0$
exterior degree term, both elements $it [\A,v_\sigma]$ and $\A^2$
are sums of terms with strictly positive exterior degrees. For
example, if $\A=\nabla^+\oplus \nabla^-$ is a direct sum of
connections, then $it [\A,v_\sigma]\in \Acal^1(N,\End(\Ecal)^-)$ and
$\A^2\in \Acal^2(N,\End(\Ecal)^+)$.

Consider the equivariant closed form
$\ch(\sigma,\A,t)(X):=\str\left(\e^{\F(\sigma,\A,t)(X)}\right)$.

\begin{defi}\label{defi:quillen}
The Quillen Chern character form $\ch_Q(\sigma)$ attached to the
connection $\A$ is the closed equivariant form
$$\ch(\sigma,\A,1)(X):=\str\left(\e^{- v_\sigma^2+i
[\A,v_\sigma]+\A^2+\mu^\A(X)}\right).$$
\end{defi}

Consider  the  \emph{transgression form }
\begin{equation}\label{eq:trangression-form}
    \eta(\sigma,\A,t)(X):=-
\str\left(iv_\sigma \,\e^{\F(\sigma,\A,t)(X)}\right).
\end{equation}
As $iv_\sigma=\frac{d}{dt}\A^\sigma(t)$, we have
$\frac{d}{dt}\ch(\sigma,\A,t)=-D(\eta(\sigma,\A,t))$. After
integration, it gives the following equality of equivariant
differential forms on $N$
\begin{equation}\label{eq:transgression-integral}
    \ch(\A)-\ch(\sigma,\A,t)=D\left(\int_0^t\eta(\sigma,\A,s)ds\right),
\end{equation}
since $\ch(\A)=\ch(\sigma,\A,0)$.

We choose a metric on the tangent bundle to $N$. Thus we obtain a
norm $\|-\|$ on $\wedge \T_n^*N\otimes \End(\Ecal_n)$
which varies smoothly with $n\in N$.

\begin{prop}\label{estimates}
Let $\Kcal_1\times \Kcal_2$ be a compact subset of $N\times \kgot$.

$\bullet$ There exists $\cst>0$ such that,
if $(n,X)\in\Kcal_1\times\Kcal_2$,
\begin{equation}\label{est1}
\Big|\!\Big| \e^{\F(\sigma,\A,t)(X)}\Big|\!\Big|(n)\leq \cst \, (1+t)^{\dim N}
\e^{-h_\sigma(n)t^2 },\quad \mathrm{for\ all\ }t\geq 0.
\end{equation}
$\bullet$ The differential forms $\ch( \sigma,\A,t)(X)$ and
$\eta(\sigma,\A,t)(X)$ (and all their partial derivatives) tends to $0$ exponentially
fast when $t\to\infty$ uniformly on compact subsets of $(N\setminus
\supp(\sigma))\times\kgot$.
\end{prop}

\begin{proof}
 To estimate $\|\e^{\F(\sigma,\A,t)(X)}\|$, we employ Proposition
\ref{prop-estimation-generale} of the Appendix, with the variable
$x=(n,X)$ and the maps  $R(n,X)=v_\sigma^2(n)$,
$S(n,X)=\mu^\A(n)(X)$, and $T(t,n,X)=it [\A,v_\sigma](n)+\A^2(n)$.
The same estimate holds for $\| D(\partial)\cdot
\e^{\F(\sigma,\A,t)(X)}\|$, when $D(\partial)$ is differential
operator on $N\times \kgot$. Hence the second point follows from
the fact that $\inf_{n\in\Kcal_1}h_\sigma(n)  >0$ when the compact
subset $\Kcal_1$ lies inside $N\setminus \supp(\sigma)$.
\end{proof}\bigskip

\medskip

The former estimates allows us to take the limit $t\to\infty$ in
(\ref{eq:transgression-integral}) on the open subset $N\setminus
\supp(\sigma)$. We get the following important lemma (see
\cite{Quillen85,pep-vergne1} for the non-equivariant case).

\begin{lem}\label{lem:quillen-equi}
We can define on $N\setminus \supp(\sigma)$ the equivariant
differential form with smooth coefficients
\begin{equation}
  \label{eq:beta}
  \beta(\sigma,\A)(X)=\int_{0}^\infty\eta(\sigma,\A,t)(X)dt,\quad
  X\in\kgot.
\end{equation}
We have $\ch(\A)|_{N\setminus
\supp(\sigma)}=D\left(\beta(\sigma,\A)\right)$.
\end{lem}

\medskip

We are in the situation of Section \ref{section:coho-relative}.
The closed equivariant form $\ch(\A)$ on $N$ and the equivariant
form $\beta(\sigma,\A)$ on $N\setminus\supp(\sigma)$ define an
even relative cohomology class $[\ch(\A),\beta(\sigma,\A)]$ in
$\Hcal^\infty(\kgot,N,N\setminus \supp(\sigma))$.

\begin{prop}\label{inds}

$\bullet$
 The class
 $[\ch(\A),\beta(\sigma,\A)]\in\Hcal^\infty(\kgot,N,N\setminus \supp(\sigma))$
 does not depend of the choice of $\A$, nor on the Hermitian
 structure on $\Ecal$. We denote it by $\chr(\sigma)$ and we call it Quillen's
 {\bf relative} equivariant Chern character.

 $\bullet$
 Let $F$ be an invariant closed subset of $N$. For $s\in [0,1]$,
 let $\sigma_s:\Ecal^+\to \Ecal^-$ be a
family of equivariant smooth morphisms such that
$\supp(\sigma_s)\subset F$. Then all classes $\chr(\sigma_s)$
coincide in $\Hcal^\infty(\kgot,N,N\setminus F)$.
\end{prop}

\begin{proof}
The proof is identical to the proof of Proposition 3.8 in
\cite{pep-vergne1} for the non equivariant case.
\end{proof}\bigskip

\subsection{Tensor product}\label{tensor}

Let $\Ecal_1, \Ecal_2$ be two  equivariant $\Zbb_2$-graded vector bundles
on $N$. The space $\Ecal_1\otimes \Ecal_2$ is a
$\Zbb_2$-graded vector bundle with even part $\Ecal_1^+\otimes
\Ecal_2^+\oplus \Ecal_1^-\otimes \Ecal_2^-$ and odd part
$\Ecal_1^-\otimes \Ecal_2^+\oplus \Ecal_1^+\otimes \Ecal_2^-$.

\begin{rem}
If $E_1$ and $E_2$ are super vector spaces, the super-algebra
$\End(E_1)\otimes \End(E_2)$ is identified  with the super-algebra
$\End(E_1\otimes E_2)$ via the following rule. For $v_1\in
E_1,v_2\in E_2, A\in \End(E_1), B\in \End(E_2)$ homogeneous
$$(A\otimes B)(v_1\otimes v_2)=(-1)^{|B||v_1|}Av_1\otimes Bv_2.$$
\end{rem}

The super-algebra $\Acal(N,\End(\Ecal_1\otimes \Ecal_2))$ can be
identified with \break
$\Acal(N,\End(\Ecal_1))\otimes\Acal(N,\End(\Ecal_2))$ where
the tensor is taken in the sense of super-algebras. Then, if $A\in
\Acal^0(N,\End(\Ecal_1)^-)$ and $B\in \Acal^0(N,\End(\Ecal_2)^-)$
are odd endomorphisms, we have $(A\otimes{\rm Id}_{\Ecal_2}+ {\rm
Id}_{\Ecal_1}\otimes B)^2=A^2\otimes {\rm Id}_{\Ecal_2}+ {\rm
Id}_{\Ecal_1}\otimes B^2$.

\bigskip

Let $\sigma_1:\Ecal_1^+\to \Ecal_1^-$ and $\sigma_2:\Ecal_2^+\to
\Ecal_2^-$ be two smooth equivariant morphisms. With the help of
$K$-invariant Hermitian structures, we define the morphism
$$
\sigma_1\odot \sigma_2: \left(\Ecal_1\otimes \Ecal_2\right)^+
\longrightarrow \left(\Ecal_1\otimes \Ecal_2\right)^-
$$
by $\sigma_1\odot \sigma_2:= \sigma_1\otimes {\rm
Id}_{\Ecal_2^+}+{\rm Id}_{\Ecal_1^+}\otimes \sigma_2+{\rm
Id}_{\Ecal_1^-}\otimes\sigma_2^*+ \sigma_1^*\otimes {\rm
Id}_{\Ecal_2^-}$.

\bigskip

Let $v_{\sigma_1}$ and $v_{\sigma_2}$ be the odd Hermitian
endomorphisms of $\Ecal_1,\Ecal_2$
associated to $\sigma_1$ and $\sigma_2$ (see (\ref{eq:v-sigma})).
Then $v_{\sigma_1\odot
\sigma_2}= v_{\sigma_1}\otimes {\rm Id}_{\Ecal_2}+{\rm
Id}_{\Ecal_1}\otimes v_{\sigma_2}$ and $v_{\sigma_1\odot
\sigma_2}^2=v_{\sigma_1}^2\otimes {\rm Id}_{\Ecal_2}+{\rm
Id}_{\Ecal_1}\otimes v_{\sigma_2}^2$. Thus the  square
$v_{\sigma_1\odot \sigma_{2}}^2$ is the sum of two commuting non
negative Hermitian endomorphisms $v_{\sigma_1}^2\otimes {\rm
Id}_{\Ecal_2}$ and ${\rm Id}_{\Ecal_1}\otimes v_{\sigma_2}^2$. It follows
that
$$
\supp(\sigma_1\odot \sigma_2)= \supp(\sigma_1)\cap \supp(\sigma_2).
$$

We can now state the main result of this section.

\begin{theo}\label{theo:chrel-produit}{\rm \bf (The relative Chern character is
multiplicative)} Let $\sigma_1,\sigma_2$ be two equivariant morphisms over $N$.
The relative equivariant cohomology classes
\begin{itemize}
  \item $\chr(\sigma_k)\in\Hcal^\infty(\kgot,N,N\setminus\supp(\sigma_k))$,
  $k=1,2$,
  \item
  $\chr(\sigma_1\odot\sigma_2)\in
  \Hcal^\infty(\kgot,N,N\setminus(\supp(\sigma_1)\cap\supp(\sigma_2)))$
\end{itemize}
satisfy the following equality
$$
\chr(\sigma_1\odot\sigma_2)=\chr(\sigma_1)\diamond\chr(\sigma_2)
$$
in
$\Hcal^\infty(\kgot,N,N\setminus(\supp(\sigma_1)\cap\supp(\sigma_2)))$.
Here $\diamond$ is the product of relative classes (see
(\ref{eq:produit-relatif-infty})).
\end{theo}

\begin{proof} The proof is identical to the proof of Theorem 4.3 in \cite{pep-vergne1}
for the non equivariant case.
\end{proof}\bigskip

\subsection{The equivariant Chern character of a morphism}

Let $\sigma:\Ecal^+\to\Ecal^-$ be an equivariant morphism on $N$.
Following  Subsection \ref{section:morphism-pF}, we consider the
image of $\chr(\sigma)$ through the map $\Hcal^\infty(\kgot,N,N\setminus
\supp(\sigma))\to \Hcal_{\supp(\sigma)}^\infty(\kgot,N)$: the following
theorem summarizes the construction of the image.

\begin{theo}\label{theo:chgood}
$\bullet$ For any invariant neighborhood $U$ of $\supp(\sigma)$, take
$\chi\in\f(N)^K$ which is equal to 1 in a neighborhood of
$\supp(\sigma)$ and with support contained in $U$. The equivariant
differential form
\begin{equation}\label{sigmaA}
c(\sigma,\A,\chi)=\chi\, \ch(\A) + d\chi \,\beta(\sigma,\A)
\end{equation}
is {\em equivariantly closed} and supported in $U$.
 Its cohomology class $c_U(\sigma)\in \Hcal_U^\infty(\kgot,N)$ does not depend of
the choice of $(\A,\chi)$ and the Hermitian structures on
$\Ecal^{\pm}$. Furthermore, the inverse family $c_U(\sigma)$ when
$U$ runs over the neighborhoods of $\supp(\sigma)$ defines a class
$$
\chg(\sigma)\in \Hcal_{\supp(\sigma)}^\infty(\kgot,N).
$$
The image of this class in $\Hcal^\infty(\kgot,N)$ is the  Chern character
$\ch(\Ecal)$ of $\Ecal$.

$\bullet$ Let $F$ be an invariant closed subset of $N$. For $s\in [0,1]$, let
$\sigma_s:\Ecal^+\to \Ecal^-$ be a family of smooth morphisms such
that $\supp(\sigma_s)\subset F$. Then all classes $\chg(\sigma_s)$
coincide in $\Hcal_F^\infty(\kgot,N)$.

\end{theo}

\bigskip

\begin{defi}\label{def:ch-compact}
    When $\sigma$ is elliptic, we denote by
\begin{equation}\label{eq:ch-c}
    \ch_c(\sigma)\in \Hcal_c^\infty(\kgot,N)
\end{equation}
the cohomology class with compact support which is the image of
$\chg(\sigma)\in \Hcal_{\supp(\sigma)}^\infty(\kgot,N)$ through the canonical
map $\Hcal_{\supp(\sigma)}^\infty(\kgot,N)\to \Hcal_{c}^\infty(\kgot,N)$.
\end{defi}

A representative of $\ch_c(\sigma)$ is given by
$c(\sigma,\A,\chi)$, where $\chi\in\f(N)^K$ is chosen with
compact support, and equal to $1$ in a neighborhood of
$\supp(\sigma)$ and $c(\sigma,\A,\chi)$ is given by Formula
(\ref{sigmaA}).

\bigskip

We will now rewrite Theorem \ref{theo:chrel-produit} for the
equivariant Chern characters  $\chg$ and $\ch_c$. Let
$\sigma_1:\Ecal_1^+\to \Ecal_1^-$ and $\sigma_2:\Ecal_2^+\to
\Ecal_2^-$ be two smooth equivariant morphisms on $N$. Let
$\sigma_1\odot \sigma_2: \left(\Ecal_1\otimes \Ecal_2\right)^+ \to
\left(\Ecal_1\otimes \Ecal_2\right)^-$ be their product.

Following (\ref{eq:produit-equi-gene-infty}), the product of the
elements $\chg(\sigma_k)\in
\Hcal_{\supp(\sigma_k)}^\infty(\kgot,N)$ for $k=1,2$ belongs to
$\Hcal_{\supp(\sigma_1)\cap\supp(\sigma_2)}^\infty(\kgot,N)=
\Hcal_{\supp(\sigma_1\odot\sigma_2)}^\infty(\kgot,N)$.

\begin{theo}\label{theo:produit-chg}
$\bullet$ We have the equality
$$
\chg(\sigma_1)\wedge \chg(\sigma_2)=
\chg(\sigma_{1}\odot\sigma_2)\quad {\rm in}\quad  \Hcal_{\supp(\sigma_1\odot
\sigma_2)}^\infty(\kgot,N).
$$

$\bullet$ If the morphisms $\sigma_1,\sigma_2$ are elliptic, we have
$$
\ch_c(\sigma_1)\wedge \ch_c(\sigma_2)=
\ch_c(\sigma_{1}\odot\sigma_2)\quad {\rm in} \quad \Hcal_{c}^\infty(\kgot,N).
$$
\end{theo}

\begin{proof} This follows from Theorem \ref{theo:chrel-produit} and the diagram
(\ref{eq:fonctoriel-p-produit}).
\end{proof}

\bigskip

The second point of Theorem \ref{theo:produit-chg} has the
following interesting refinement. Let $\sigma_1,\sigma_2$ be two
equivariant morphisms on $N$ which are {\bf not elliptic}, and
assume that the product $\sigma_1\odot\sigma_2$ is {\bf elliptic}.
Since $\supp(\sigma_1)\cap\supp(\sigma_2)$ is compact, we consider
equivariant neighborhoods  $U_k$ of $\supp(\sigma_k)$ such that
$\overline{U_1}\cap\overline{U_2}$ is compact. Choose
$\chi_k\in\f(N)^K$ supported on $U_k$ and equal to $1$ in a
neighborhood of $\supp(\sigma_k)$.  Then, the  equivariant
differential form $c(\sigma_1,\A_1,\chi_1)\wedge
c(\sigma_2,\A_2,\chi_2)$ is compactly supported on $N$, and we
have
$$
\ch_c(\sigma_1\odot\sigma_2)=c(\sigma_1,\A_1,\chi_1)\wedge
c(\sigma_2,\A_2,\chi_2)\quad
{\rm in} \quad \Hcal_c^\infty(\kgot,N).
$$
Note that the equivariant differential forms $c(\sigma_k,\A_k,\chi_k)$ are not compactly supported.

\subsection{Retarded  construction}\label{sec:retarded}

We have defined a representative of the Chern characters
$\chr(\sigma)$ and $\chg(\sigma)$ using the one-parameter family
$\A^\sigma(t)$ of super-connections, for $t$ varying between $0$
and $\infty$. Quillen's Chern character
$\ch_Q(\sigma)=\ch(\sigma,\A,1)$ is another representative. We
will compare them in appropriate cohomology spaces in the next
section.

 Consider any
$T\in\Rbb$. We have $\ch(\sigma,\A,T)=D(\beta(\sigma,\A,T))$ with
$ \beta(\sigma,\A,T)=\int_T^\infty\eta(\sigma,\A,t)dt$. It is easy
to check that the following equality
\begin{eqnarray}\label{eq-chrel-t-delta}
\lefteqn{\!\!\!\!\!\!\!\!\!\!\!\!\!\!\!\!\!\!\!\!\!\!\!
\Big(\ch(\A),\beta(\sigma,\A)\Big)-\Big(\ch(\sigma,\A,T),\beta(\sigma,\A,T)\Big)
=}\\
& &D_{\rm rel}\left(\int_0^T\!\!\!\eta(\sigma,\A,t)dt\, ,\,
0\right)\nonumber
\end{eqnarray}
holds in $\Acal^\infty(\kgot,N,N\setminus \supp(\sigma))$. Hence we get
the following
\begin{lem}\label{lem:other-representant}
 For any $T\in\Rbb$, the relative Chern character $\chr(\sigma)$ satisfies
$$
\chr(\sigma)=\Big[\ch(\sigma,\A,T),\beta(\sigma,\A,T)\Big]\quad {\rm
in}\quad \Hcal^\infty(\kgot,N,N\setminus \supp(\sigma)).
$$
\end{lem}

Using Lemma \ref{lem:other-representant}, we get
\begin{lem}\label{lem:c-sigma-t}
For any $T\geq 0$, the class $\chg(\sigma)$ can be defined with the
forms $c(\sigma,\A,\chi,T)=\chi\, \ch(\sigma,\A,T) +
d\chi\,\beta(\sigma,\A,T)$.
\end{lem}
\begin{proof} It is due to the following transgression
\begin{equation}\label{eq:trangression-c-t}
    c(\sigma,\A,\chi)-c(\sigma,\A,\chi,T)=D(\chi\int_0^T\!\!\!\eta(\sigma,\A,t)dt),
\end{equation}
which follows from (\ref{eq-chrel-t-delta}).
\end{proof}\bigskip

In some situations the Quillen's Chern form
$\ch_Q(\sigma)=\ch(\sigma,\A,1)$ enjoys good properties relative to
the integration. So it is natural to compare the equivariant differential form
$c(\sigma,\A,\chi)$ and $\ch(\sigma,\A,1)$.

\begin{lem}\label{retardord}
We have
$$c(\sigma,\A,\chi)-\ch_Q(\sigma)=
D\left(\chi\int_{0}^1\eta(\sigma,\A,s)ds\right)+ D\Big((\chi-1)
\beta(\sigma,\A,1)\Big).$$
\end{lem}
\begin{proof} This follows immediately from the transgressions
(\ref{eq:transgression-integral}) and (\ref{eq:trangression-c-t}).

\end{proof}\bigskip

\subsection{Quillen Chern character with Gaussian look}\label{sec:gaussian}

As we have seen, Mathai-Quillen  gives an explicit representative
with ``Gaussian look''  of the Thom class of  a Euclidean  vector
bundle $\Vcal\to M$. Similarly, they give an explicit representative
with ``Gaussian look''  of the Bott class of  a complex
equivariant vector bundle $\Vcal\to M$. The purpose of this paragraph
is to compare the Mathai-Quillen construction  of Chern characters
with ``Gaussian look'' and the  relative construction.

Let $\Vcal$ be a real $K$-equivariant vector bundle over a manifold
$M$. We denote by $p: \Vcal\to M$ the projection. We denote by
$(x,\xi)$ a point of $\Vcal$ with $x\in M$ and $\xi\in \Vcal_x$.
 Let $\Ecal^\pm\to M$ be two $K$-equivariant Hermitian vector bundles. We
consider a $K$-equivariant morphism $\sigma: p^*\Ecal^+\to
p^*\Ecal^-$ on $\Vcal$.

We choose a metric on the fibers of the fibration $\Vcal\to M$. We
work under the following assumption on $\sigma$.

\begin{assu}\label{assum:dec-rap-ord}
The morphism $\sigma: p^*\Ecal^+\to p^*\Ecal^-$ and all its
partial derivatives have at most a polynomial growth along the
fibers of $\Vcal\to M$. Moreover we assume  that, for any compact
subset $\Kcal$ of $M$, there exist $R\geq 0$ and $c>0$ such
that\footnote{This inequality means that $\|\sigma(x,\xi)w\|^2\geq
c\|\xi\|^2\|w\|^2$ for any $w\in \Ecal_x$.} $v_\sigma^2(x,\xi)\geq
c\|\xi\|^2$ when $\|\xi\|\geq R$ and $x\in \Kcal$.
\end{assu}

Let $\nabla=\nabla^+\oplus \nabla^-$ be a $K$-invariant connection
on $\Ecal\to M$, and consider the super-connection
$\A=p^*\nabla$ so that $\A^\sigma(t)=p^*\nabla+i tv_\sigma$.
Then, the form $\ch(\sigma,\A,1)(X)$ has a ``Gaussian'' look.
\begin{lemm}\label{prop:ch=dec-rapid}
    The equivariant differential forms $\ch(\sigma,\A,1)(X)$ and $\beta(\sigma,\A,1)(X)$ are
    rapidly decreasing along the fibers.
\end{lemm}

\begin{proof} The equivariant curvature of $\A^\sigma(t)$ is
$$
\F(t)(X)= p^*\F(X) -t^2v_\sigma^2+it[p^*\nabla,v_\sigma].
$$
Here $\F\in\Acal^2(M,\End(\Ecal))$ is the equivariant curvature of
$\nabla$.

To estimate $\e^{\F(t)}(X)$, we apply Lemma \ref{suffit} of the
Appendix, with $H=t^2 v_{\sigma}^2,$ and $R=-p^*\F-it
[p^*\nabla, v_\sigma]$, $S=\mu^{\A}(X)$. The proof is very
similar to the proof of Lemma  5.17 in \cite{pep-vergne1} for the
non equivariant case and we skip it.

\end{proof}\bigskip

\begin{theo}

Quillen's Chern character  form $\ch_Q(\sigma)\in\Acal^{\infty}_{\dr}(\kgot,\Vcal)$ represents
the image of the class $\chg(\sigma)\in\Hcal^{\infty}_{\supp(\sigma)}(\kgot,\Vcal)$ in
$\Hcal^{\infty}_{\dr}(\kgot,\Vcal)$.
\end{theo}

\begin{proof}
The proof is entirely analogous to Proposition 5.18 in
\cite{pep-vergne1} and we skip it.
\end{proof}\bigskip

\subsection{A simple example}

Thus there are  three useful representations of the equivariant
Chern character of a morphism $\sigma$: the relative Chern
character $\chr(\sigma)\in
\Hcal^\infty(\kgot,\Vcal,\Vcal\setminus{\supp(\sigma)})$, the Chern
character with support $\chg(\sigma)\in
\Hcal_{\supp(\sigma)}^\infty(\kgot,\Vcal)$ and (in the case of vector
bundles) the Chern character with Gaussian look. We will describe
these explicit representatives  in the three cohomology spaces in
a very simple example.

Recall the following convention. Let $V=V^+\oplus V^-$ be a
$\Zbb_2$-graded finite dimensional complex vector space and
$\Acal$ a super-commutative ring (the ring of differential forms
on a manifold for example). Consider the ring $\End(V)\otimes
\Acal$. Let $\{e_i\}_{i=1}^{\dim V^+}$ be a basis of $V^+$ and
$\{f_j\}_{j=1}^{\dim V^-}$ a basis of $V^-$. Consider the odd
endomorphism $M_j^i:V\to V$ such that $M_j^i(e_i)=f_j$, and
sending all other basis elements of $V$ to $0$. Similarly,
consider the odd endomorphism $R_i^j:V\to V$ such that
$R_i^j(f_j)=e_i$, and sending all other basis elements of $V$ to
$0$.
\begin{con}\label{co}
 Let $a\in \Acal$. The matrix written with $a$ in column
$i$ and row $j$, and $0$ for all other entries, represents
$M_{j}^i\otimes a$ in the ring $End(V)\otimes \Acal$. The matrix
written with $a$ in column $\dim V^++j$ and row $i$, and $0$ for
all other entries, represents $R_{i}^j\otimes a$ in the ring
$End(V)\otimes \Acal$.
\end{con}

Let $U(1)=K$ be the circle group.
We identify the Lie algebra $\ugot(1)$ of $U(1)$ with $\Rbb$ so that the
exponential map is $\theta\mapsto \e^{i\theta}$.

 Consider the case where $\Vcal=\Rbb^2\simeq \Cbb$ and $K=U(1)$ acts
by rotation: $t\cdot z= tz$ for $z\in \Cbb$ and $t\in U(1)$. Take
$E^+=\Vcal\times \wedge^0 \Vcal$ and $E^-=\Vcal\times \wedge^1 \Vcal$. The
action of $U(1)$ on $E^+\equiv\Cbb\times \Cbb$ is $t\cdot
[z,u]=[tz,u].$ The action of $U(1)$ on $E^-\equiv\Cbb\times \Cbb$
is $t\cdot [z,u]=[tz,tu].$

We consider the Bott symbol $\sigma_b(z)=z$ which produces the map
$\sigma_b([z,u])=[z,zu]$ from $E^+$ to $E^-$. Then, the bundle map
$\sigma_b$ commutes with the action of $U(1)$ and defines an
element of $\KK_{U(1)}(\Rbb^2)$.
 Recall that the Bott
isomorphism tell us that $\KK_{U(1)}(\Rbb^2)$ is a free module over
$\KK_{U(1)}(\textrm{pt})=R(U(1))$ with base $\sigma_b$.

 We choose on $E^\pm$
the trivial connections $\nabla^+=\nabla^-=d$. Let
$\A=\nabla^+\oplus \nabla^-$. The moment of $\A$ is the map
$$
\mu^\A(\theta)= \left(\begin{array}{cc}0& 0\\ 0 &
{i\theta}\end{array}\right).
$$

The equivariant curvature  $\F(\theta)$ of $\A$ is equal to
$\mu^\A(\theta)$, thus we have the following formula for the
equivariant Chern character
\begin{equation}\label{chbete}
 \ch(\A)(\theta)=1-\e^{{i\theta}}.
\end{equation}

With the conventions of \ref{co},
 the equivariant curvature $\F(t)$ of
the super-connection $\A^{\sigma_b}(t):=\A +it v_{\sigma_b}$ is
written in matrix form as
$$\F(t)(\theta)=\left(\begin{array}{cc}
-t^2|z|^2 & 0\\
0 & -t^2|z|^2\\
\end{array}\right)+\left(\begin{array}{cc}
0 & -itd\overline{z}\\
-itdz & 0\\
\end{array}\right)+ \left(\begin{array}{cc}
0& 0\\ 0 & {i\theta}\end{array}\right),
$$
for $\theta\in \ugot(1)\simeq \Rbb$. We compute
$\e^{\F(t)(\theta)}$ using Volterra's formula. We obtain
$$
\e^{\F(t)(\theta)}=\e^{-t^2 |z|^2} \left(\begin{array}{cc}
1+(g'({i\theta})-g({i\theta})) t^2 dzd\overline{z}& it g({i\theta})
d\overline{z}
\\ it g({i\theta}) dz & \e^{i\theta} +g'({i\theta}) t^2 dzd\overline{z}\end{array}\right)
$$
where $g({i\theta})=\frac{\e^{i\theta} -1}{{i\theta}}$. Hence
$\str(\e^{\F(t)({\theta})}) =-g({i\theta})
({i\theta} + t^2 dzd\overline{z})\e^{-t^2 |z|^2}$. Here
\begin{eqnarray*}
 \eta(\sigma_b,\A,t)({\theta}) &=& -i \str\left(
\left(\begin{array}{cc}0 & \overline{z}\\z & 0\end{array}\right)
\e^{\F(t)({\theta})}\right) \\
   &=&g({i\theta})(zd\overline{z}-\overline{z}dz)\, t \, \e^{-t^2|z|^2}.
\end{eqnarray*}
For $z\neq 0$, we can integrate $t\mapsto \eta(\sigma_b,\A,t)(\theta)$ from
$0$ to $\infty$ and we obtain
$$
\beta(\sigma_b,\A)(\theta)=g({i\theta})\frac{
z\,d\overline{z}-{\overline{z}}\,dz}{2|z|^2}.
$$

Take $f\in\f(\Rbb)$ with compact support and equal to $1$ in a
neighborhood of $0$. Let $\chi(z):= f(|z|^2)$.

Similarly to the Thom form, we can give formulae for the three
different representatives of the Chern character.

\begin{prop}
$\bullet$ The class $\chr(\sigma_b)\in \Hcal^\infty(\ugot(1),\Rbb^2,\Rbb^2\setminus\{0\})$
is represented by the couple of equivariant differential forms:
$$
\left(1-\e^{i\theta},  \frac{\e^{i\theta} -1}{{i\theta}}\frac{
z\,d\overline{z}-{\overline{z}}\,dz}{2|z|^2}\right).
$$

$\bullet$ The Chern character with compact support is represented
by the  equivariant differential form
\begin{eqnarray*}
c(\sigma_b,\A,\chi)&=&\chi(1-\e^{i\theta})+d\chi\beta(\sigma_b,\A)\\
&=&\frac{\e^{i\theta} -1}{{i\theta}}\left(- f(|z|^2){i\theta} +
f'(|z|^2)dz\wedge d\overline{z}\right).
\end{eqnarray*}

$\bullet$ Quillen's Chern character with Gaussian look
$\ch_Q(\sigma_b)$ is represented by the equivariant differential form

$$\frac{\e^{i\theta} -1}{{i\theta}}\e^{-\|z\|^2}\left({i\theta} +
dz\wedge d\overline{z}\right).
$$
\end{prop}

Comparing with Example \ref{exa:Vdim2},
 we see that the Chern character form in all of these different
versions is proportional to the Thom form

$$\ch(\sigma_b)(\theta)=(2i\pi)\frac{\e^{i\theta} -1}{{i\theta}} {\rm Th}(V)(\theta).$$

We will see in the next section that this identity generalizes to
any Euclidean vector spaces.

\section{Comparison between relative Thom classes and Bott classes:Riemann Roch formula}
\label{comparedgaussian}

Let $p:\Vcal\to M$ be a $K$-equivariant Euclidean vector bundle
over $M$ of {\em even rank}. Here we compare the relative  Chern
character of the Bott symbol and the relative Thom class. Both
classes live in the relative equivariant cohomology space
$\Hcal^\infty(\kgot,\Vcal, \Vcal\setminus M)$. The formulae relating
them is an important step  in the proof
of the Grothendieck-Riemann-Roch relative theorem \cite{Borel-Serre}, as well
as the Atiyah-Singer theorem \cite{Atiyah,Atiyah-Singer-1}. As usual, the relation is
deduced from an explicit computation in equivariant cohomology of
a vector space.

\subsection{The relative  equivariant Bott class of a vector space}
\label{sec:bott-class-V}

We consider first the case of an oriented Euclidean vector space
$V$ of dimension $d=2n$. Let
$$
\textbf{c} : {\rm Cl}(V)\to \End_\Cbb(\Scal)
$$
be the spinor representation of the Clifford algebra of $V$. We use
the conventions  of \cite{B-G-V}[chapter 3] for the spinor
representation. In particular, as a vector space, the Clifford
algebra is identified with the real exterior algebra of $V$. The
orientation of $V$ gives a decomposition
$\Scal=\Scal^+\oplus\Scal^-$ which is is stable under the action of
the group ${\rm Spin}(V)\subset {\rm Cl}(V)$.

We consider the ${\rm Spin}(V)$-equivariant vector bundles $\Scal_V^\pm:=
V\times \Scal^\pm$ over $V$ : recall that the action of ${\rm Spin}(V)$ on
the base $V$ is through the twofold covering $\tau : {\rm
Spin}(V)\to {\rm SO}(V)$.

The Clifford module $\Scal$ is provided with an Hermitian inner
product such that $\textbf{c}(x)^*=-\textbf{c}(x)$, for $x\in V$.
We work with the equivariant morphism $\sigma_V: \Scal_V^+\to \Scal_V^-$
defined by: for $x\in V$,
$$
\sigma_V(x):= -i\textbf{c}(x): \Scal^+\longrightarrow \Scal^-.
$$

Then the odd linear map $v_\sigma(x):\Scal\to\Scal$ is equal to
$-i\clif(x)$. We choose on $\Scal_V^{\pm}$ the trivial connections
$\nabla^+=\nabla^-=d$. Thus the super-connection $A_t:=\A+it
v_\sigma$ is
$$
\A_t=  d + t\textbf{c}(x).
$$

The Lie algebra $\spingot(V)$ of ${\rm Spin}(V)$ is  identified
with the Lie sub-algebra $\wedge^2 V$ of ${\rm Cl}(V)$ and the
exponential map is the exponential inside the Clifford algebra. The
differential of the action of ${\rm Spin}(V)$ in $\Scal$ is
$Y\mapsto \textbf{c}(Y)$. We denote by $Y\in\spingot(V)\mapsto
Y^\tau\in\sogot(V)$ the differential of the homomorphism $\tau$ (it
is an isomorphism). We need the function $\jdemi_V:\sogot(V)\to
\Rbb$ defined by
$$
\jdemi_V(X)=
\det{}^{1/2}\left(\frac{\e^{X/2}-\e^{-X/2}}{X/2}\right).
$$
Then $\jdemi_V(X)$ is invertible near $X=0$ as $\jdemi_V(0)=1$.

For $Y\in \spingot(V)$, the moment $\mu(Y)=\Lcal(Y)-
[\A_t,\iota(VY)]$ is equal to $\bc(Y)$.
Hence the equivariant curvature of the super-connection $\A_t$ is
the  function $\F_t:\spingot(V)\to\Acal(V,\End_\Cbb(\Scal))$, given by
\begin{equation}\label{eq:curvature-spin}
    \F_t(Y)= -t^2\|x\|^2 + t\sum_k dx_k \bc_k +
\sum_{k<l}Y_{kl}\bc_k \bc_l.
\end{equation}
Here $Y=\sum_{k<l}Y_{kl}e_k\wedge e_l\in \spingot(V)$ and
$\bc(Y)=\sum_{k<l}Y_{kl}\bc_k \bc_l$ : we have
denoted by $\bc_k$ the odd endomorphism of $\Scal$ produced by
$e_k\in V\subset {\rm Cl}(V)$.

This formula is very similar to the form $f_t$ that we used to
construct the Thom form, see Subsection \ref{subsec:Thom-vector}. For
$Y\in\spingot(V)$, we write
\begin{eqnarray*}
f_t(Y^{\tau})&=&
 -t^2\|x\|^2 + t\sum_k dx_k e_k + \frac{1}{2}\sum_{k<l}Y^{\tau}_{kl}~e_k\wedge e_l\\
&=&-t^2\|x\|^2 + t\sum_k dx_k e_k +\sum_{k<l}Y_{kl}~e_k\wedge e_l.
\end{eqnarray*}

In order to compute the relative Chern character   of $\sigma$ we
follow the strategy of \cite{B-G-V}[Section 7.7]. However, our
convention for the Chern character $\ch(\A)=\str(\e^{\A^2})$ is
different from the one of \cite{B-G-V}, which decided that
$\ch(\A)=\str(\e^{-\A^2})$. Thus we carefully check signs in our
formulae.

We consider in parallel the closed equivariant forms on $V$
$$
\ch(\A_t)(Y)=\str\left(\e^{\F_t(Y)}\right)\quad , \quad
\mathrm{C}_{\wedge}^t(Y^\tau):=
\bere\left(\e^{f_t(Y^\tau)}\right),\quad Y\in\spingot(V).
$$

In the first case, the exponential is computed in the super-algebra
$\End(\Scal)\otimes \Acal(V)$  and in the second case, the
exponential is computed in the super-algebra $ \Acal(V)\otimes \wedge V.$

We also consider in parallel the equivariant forms
$$
\eta_{\textbf{c}}^t(Y)= -\str\left((\sum_k x_k
\bc_k)\e^{\F_t(Y)}\right)\quad , \quad \eta_{\wedge}^t(Y^\tau):=
-\bere\left((\sum_k x_k e_k)\e^{f_t(Y^\tau)}\right).
$$

\begin{prop}\label{prop:ch-wedge}
Let $Y\in \spingot(V)$. We have
$$
\ch(\A_t)(Y)=(-2i)^n ~
\jdemi_V(Y^\tau)~\mathrm{C}_{\wedge}^t(Y^\tau),\quad
\eta_{\textbf{c}}^t(Y)=(-2i)^n ~ \jdemi_V(Y^\tau) ~
\eta_{\wedge}^t(Y^\tau).
$$
\end{prop}

\begin{proof} The first identity is proved in \cite{B-G-V}[Section 7.7] with
other conventions. For clarity we perform the computations. We fix
$Y\in\spingot(V)\simeq\wedge^2 V$ and we take an \emph{oriented}
orthonormal base $(e_i)$ of $V$ such that $Y=\sum_{k=1}^n \lambda_k
e_{2k-1}\wedge~e_{2k}$. Hence $\jdemi_V(Y^\tau)=\Pi_k
\frac{\sin\lambda_k}{\lambda_k}$. For $1\leq k\leq n$, let
\begin{eqnarray*}
    B_k &=& t (dx_{2k-1} \bc_{2k-1} + dx_{2k} \bc_{2k}) + \lambda_k \bc_{2k-1}\bc_{2k} \\
    b_k &=& t (dx_{2k-1} e_{2k-1} + dx_{2k} e_{2k}) + \lambda_k e_{2k-1}\wedge e_{2k}.
  \end{eqnarray*}
Since $B_i B_j=B_j B_i$ and $b_i b_j=b_j b_i$, we have
$\e^{\F_t(Y)}= \e^{-t^2\|x\|^2}\Pi_k \e^{B_k}$ and
$\e^{f_t(Y^\tau)}= \e^{-t^2\|x\|^2}\Pi_k \e^{b_k}$.

\begin{lem}\label{lem: exp-B}
 We have
 \begin{eqnarray*}
   \e^{B_k} &=& \cos \lambda_k + t^2(\frac{\sin\lambda_k -
   \lambda_k \cos\lambda_k}{\lambda_k^2})dx_{2k-1} dx_{2k} +\\
    & & \frac{\sin\lambda_k}{\lambda_k}(\lambda_k -t^2 dx_{2k-1} dx_{2k}) \bc_{2k-1} \bc_{2k}
  +t \frac{\sin\lambda_k}{\lambda_k}(dx_{2k-1} \bc_{2k-1}+ dx_{2k} \bc_{2k})
 \end{eqnarray*}
 in $\Acal(V,\End_\Cbb(\Scal))$, and
$$
\e^{b_k}= 1 + (\lambda_k -t^2 dx_{2k-1} dx_{2k}) e_{2k-1}\wedge
e_{2k}
  +t (dx_{2k-1} e_{2k-1}+ dx_{2k} e_{2k})
$$
in $\Acal(V)\otimes\wedge V$.
\end{lem}
\begin{proof}The identity for $\e^{b_k}$ is obvious since $(b_k)^i=0$ when $i>2$. When
$\lambda_k=0$ the identity for $\e^{B_k}$ can be proved directly
since $(dx_{2k-1} \bc_{2k-1} + dx_{2k} \bc_{2k})^i=0$ for $i>2$.
When $\lambda_k\neq 0$ we can write
$$
B_k=\lambda_k(\bc_{2k-1}+ t\lambda_k^{-1}dx_{2k})(\bc_{2k}-
t\lambda_k^{-1}dx_{2k-1}) -t^2 \lambda_k^{-1}dx_{2k-1}dx_{2k}.
$$
If we let $\xi_1=\bc_{2k-1}+ t\lambda_k^{-1}dx_{2k}$ and
$\xi_2=\bc_{2k}- t\lambda_k^{-1}dx_{2k-1}$, then
$\xi_1^2=\xi_2^2=-1$, $\xi_1\xi_2+\xi_2\xi_1=0$ and the $\xi_i$
commute with $dx_{2k-1}dx_{2k}$. Thus, we see that
\begin{eqnarray*}
  \e^{B_k} &=& \e^{\lambda_k\xi_1\xi_2} \e^{-t^2 \lambda_k^{-1}dx_{2k-1}dx_{2k}} \\
          &=& \left(\cos\lambda_k +\sin\lambda\xi_1\xi_2\right)
          \left(1-t^2 \lambda_k^{-1}dx_{2k-1}dx_{2k}\right)\\
          &=&\Big(\cos\lambda_k +\sin\lambda(\bc_{2k-1}+ t\lambda_k^{-1}dx_{2k})
          (\bc_{2k}- t\lambda_k^{-1}dx_{2k-1})\Big)\times \\
          && \Big(1-t^2 \lambda_k^{-1}dx_{2k-1}dx_{2k}\Big).
\end{eqnarray*}

\end{proof}

\bigskip

Since $\str(\bc_1\cdots \bc_{2n})= (-2i)^n$ and $\str(
\bc_{i_1}\cdots \bc_{i_l})=0$ for $l<2n$, Lemma \ref{lem: exp-B}
gives that $\ch(\A_t)(Y)=(-2i)^n
\left(\jdemi_V\mathrm{C}_{\wedge}^t\right)(Y^\tau)$ where
$$
\mathrm{C}^t_\wedge(Y^\tau)=\e^{-t^2\|x\|^2}\Pi_{k=1}^n(\lambda_k
-t^2 dx_{2k-1} dx_{2k}).
$$
We found also that $\eta_{\textbf{c}} ^t(Y)=(-2i)^n
\left(\jdemi_V\eta_{\wedge}^t\right)(Y^\tau)$ with
$$
\eta_{\wedge}^t(Y^\tau)= -t \e^{-t^2\|x\|^2}\sum_{k=1}^n\Pi_{i\neq
k}(\lambda_i -t^2 dx_{2i-1} dx_{2i}) (x_{2k} dx_{2k-1} - x_{2k-1}
dx_{2k}).
$$

\end{proof}

\bigskip

The difference of the Chern character of the bundles
$\Scal^{\pm}_V$ with trivial connection $d$ is
$\ch(\Scal_V^+)(Y)-\ch(\Scal_V^+)(Y)=\ch(\A_0)(Y).$
By the preceding calculation, we obtain
$$
\ch(\A_0)(Y)
=(-2i)^n ~ \jdemi_V(Y^\tau)
\Pf\left(\hbox{$\frac{Y^\tau}{2}$}\right).
$$

To compute the relative Chern character of the morphism $\sigma_V$
we need the ${\rm Spin}(V)$-equivariant form
$$
\beta_{\textbf{c}}(Y)=\int_0^{\infty}\eta_{\textbf{c}}^t(Y)dt
$$
which is defined on $V\times \{0\}$. By Proposition
\ref{prop:ch-wedge}, we obtain

\begin{lemm}
The ${\rm Spin}(V)$-equivariant form $\beta_{\textbf{c}}(Y)$
and $\beta_{\wedge}(Y^{\tau})=\int_0^{\infty}\eta_{\wedge}^t(Y^{\tau})dt$ are related on
$V\times \{0\}$ by
$$
\beta_{\textbf{c}}(Y)=(-2i)^n \jdemi_V(Y^\tau)
\beta_{\wedge}(Y^\tau),\quad Y\in\spingot(V).
$$
\end{lemm}

We then obtain the following comparison between the Thom classes
and the Chern characters of the symbol $\sigma_V$.

\begin{theo}\label{prop:chern-V-spin}
$\bullet$ We have the following equality in
$\Hcal^\infty(\spingot(V),V, V\setminus \{0\})$
\begin{equation}\label{eq:chernr-V}
\chr(\sigma_V)(Y) =(2i\pi)^n~ \jdemi_V(Y^\tau)~ \tur(V)(Y^\tau).
\end{equation}

$\bullet$ We have the following equality in
$\Hcal^\infty_c(\spingot(V),V)$
\begin{equation}\label{eq:chernc-V}
\chc(\sigma_V)(Y) =(2i\pi)^n ~ \jdemi_V(Y^\tau)~ \tuc(V)(Y^\tau).
\end{equation}

$\bullet$ We have the following equality in
$\Hcal^\infty_{\dr}(\spingot(V),V)$
\begin{equation}\label{eq:cherndr-V}
\ch_Q(\sigma_V)(Y) =(2i\pi)^n ~ \jdemi_V(Y^\tau)~ \tumq(Y^\tau).
\end{equation}

\end{theo}

\subsection{The Bott class of a vector bundle}

Let $M$ be a manifold equipped with the action of a compact Lie group
$K$. Let $p:\Vcal\to M$ be an oriented $K$-vector bundle of even rank $2n$ over M.

\subsubsection{The Spin case}

We assume that $\Vcal$ has a $K$-equivariant spin structure. Thus
$\Vcal$ is associated to a $K$-equivariant principal bundle $P\to
M$ with structure group ${\tilde G}=\spin(V)$. We denote by
${\tilde \ggot}$ the Lie algebra $\spingot(V)$.

Let $\omega$ be a $K$-invariant connection one form on $P$, with
curvature form $\Omega$. For $X\in\kgot$, we denote by
$\mu(X)=-\iota(VX)\omega\, \in\f(P)\otimes{\tilde \ggot}$ the
moment of $X$. The equivariant curvature form is
$\Omega(X)=\Omega+\mu(X),\, X\in\kgot$. For any ${\tilde
G}$-manifold $Z$, we consider the Chern-Weil homomorphisms
$$
\phi_\omega^Z:\Acal^\infty({\tilde \ggot},Z) \longrightarrow
\Acal^\infty(\kgot, \Zcal)
$$
where $\Zcal=P\times_{\tilde G} Z$.

Let $\Scal_M= P\times_{\tilde G} \Scal$ be the corresponding
$K$-equivariant spinor bundle over $M$. Let $p^*\Scal_M\to
\Vcal$ be the pull-back of $\Scal_M$ to $\Vcal$ by the projection
$p:\Vcal\to M$. We have $p^*\Scal_M=P\times_{\tilde G}
\Scal_V$ where $\Scal_V=V\times \Scal$ is the trivial bundle over
$V$.

Let $\nabla^{\Scal_M}$ be the connection on $\Scal_M\to M$ induced by the connection
form $\omega$~: $\nabla^{\Scal_M} :=d + \clif(\omega)$ where $\clif$ is the
representation of ${\rm Spin}(V)$ on $\Scal$.

Let $\xx:\Vcal\to p^*\Vcal$ be the canonical section.
We consider now the $K$-equivariant morphism on $\Vcal$,
$\sigma_\Vcal:p^*\Scal^+_M\longrightarrow p^*\Scal^-_M$,
defined by
$$
\sigma_\Vcal=-i\, \clif(\xx).
$$
We consider the family of super-connections on $p^*\Scal_M$
defined by $\A^{\sigma_\Vcal}_t=p^*\nabla^{\Scal_M}+ tc(\xx)$.
In the previous section we worked with a family $\A_t$ of
super-connections on the trivial bundle $\Scal\times V\to V$. Let
$\ch(\sigma_\Vcal,\A,t)\in \Acal^\infty(\kgot,\Vcal)$ and
$\ch(\A_t)\in\Acal^\infty({\tilde \ggot}, V)$ be the corresponding
Chern forms. Let $\eta(\sigma_\Vcal,\A,t)\in
\Acal^\infty(\kgot,\Vcal)$ and $\eta^t\in\Acal^\infty({\tilde
\ggot}, V)$ be the corresponding transgression forms.

\begin{lem}\label{lem:chern-weil-chern}
We have the following equalities:
$$
\phi_\omega\left(\ch(\A_t)\right)=\ch(\sigma_\Vcal,\A,t)\quad \mathrm{in}
\quad  \Acal^\infty(\kgot,\Vcal),
$$
and
$$
\phi_\omega\left(\eta^t\right)=\eta(\sigma_\Vcal,\A,t)\quad \mathrm{in}
\quad  \Acal^\infty(\kgot,\Vcal\setminus M).
$$\end{lem}

\begin{proof}  See \cite{B-G-V}, Section7.7.
\end{proof}

\medskip

Let $\nabla^\Vcal$ be the connection on $\Vcal\to M$ induced by the
connection form $\omega$~: $\nabla^\Vcal =d + \tau(\omega)$ where
$\tau:\spin(V)\to \So(V)$ is the double cover. Let
$\F^\Vcal(X)=(\nabla^\Vcal-\iota(VX))^2+ \Lcal(X),\, X\in\kgot$ be
the equivariant curvature of $\nabla^\Vcal$.

\begin{defi}
We associate to the $K$-equivariant (real) vector bundle $\Vcal\to M$ the closed
$K$-equivariant form on $M$ defined by
$$
\jdemi(\nabla^\Vcal)(X):=\det{}^{1/2}
\left(\dfrac{\e^{\frac{\F^\Vcal(X)}{2}}-\e^{-\frac{\F^\Vcal(X)}{2}}}{\F^\Vcal(X)}\right),
\quad X\in\kgot.
$$
We denote $\jdemi(\Vcal)(X)$ its cohomology class in $\Hcal^{\infty}(\kgot,M)$.
It is an invertible class near $X=0$ and its inverse  $\jdemi(\Vcal)(X)^{-1}$
is the equivariant $\hat A$-genus of $\Vcal$.
\end{defi}

It is easy to see that the image of the invariant polynomial
$Y\mapsto \jdemi_V(Y^\tau)$ by the Chern-Weil homomorphism
$\phi_\omega$ is equal to $\jdemi(\nabla^\Vcal)$ (see
\cite{B-G-V}, Section7.7).

\medskip

We consider the sub-space
$\Acal^{\infty}_{\cf}(\kgot,\Vcal)\subset
\Acal^{\infty}(\kgot,\Vcal)$ of $K$-equivariant forms on $\Vcal$
which have a compact support in the fibers of $p:\Vcal\to M$.
Let $\Hcal^{\infty}_{\cf}(\kgot,\Vcal)$ be the corresponding
cohomology space. The Chern-Weil homomorphism
$\phi_\omega:\Acal^\infty({\tilde \ggot}, V)\to
\Acal^{\infty}(\kgot,\Vcal)$ maps the sub-space
$\Acal^\infty_c({\tilde \ggot}, V)$ into
$\Acal^{\infty}_{\cf}(\kgot,\Vcal)$.

Consider now the equivariant morphism $\sigma_\Vcal$ on $\Vcal$.
The support $\sigma_\Vcal$ is equal to $M$, hence its relative
Chern Character $\chr(\sigma_\Vcal)$ belongs to
$\Hcal^{\infty}(\kgot,\Vcal,\Vcal\setminus M)$.

If we take the image of the equalities of Theorem \ref{prop:chern-V-spin}
by the Chern-Weil homomorphism $\phi_\omega$
we obtain the following

\begin{prop}\label{prop:ch-spin-Thom}

We have the equalities:
$$
\chr(\sigma_\Vcal)= (2i\pi)^n p^*\left(\jdemi(\Vcal)\right)
\tur(\Vcal)\quad \mathrm{in} \quad
\Hcal^{\infty}(\kgot,\Vcal,\Vcal\setminus M).
$$
$$
\chcf(\sigma_\Vcal)= (2i\pi)^n p^*\left(\jdemi(\Vcal)\right)
\tucf(\Vcal)\quad \mathrm{in} \quad
\Hcal^{\infty}_{\cf}(\kgot,\Vcal).
$$
$$\ch_{Q}(\sigma_\Vcal)= (2i\pi)^n p^*\left(\jdemi(\Vcal)\right)
\tumq(\Vcal)\quad \mathrm{in} \quad
\Hcal^{\infty}_{\dr}(\kgot,\Vcal).
$$
\end{prop}

\subsubsection{The Spin$^\Cbb$ case}

We assume here that the vector bundle $p:\Vcal\to M$ has a $K$-equivariant $\spinc$
structure. Thus $\Vcal$ is associated to a $K$-equivariant principal
bundle $P^{\rm c}\to M$ with structure group $G^{\rm c}:=\spinc(V)$.

Let $U(1):=\{\e^{i\theta}\}$ be the circle group with  Lie algebra of $\ugot(1)\sim \Rbb$.
The group $\spinc(V)$ is the quotient $\spin(V)\times_{\Zbb_2}U(1)$,
where $\Zbb_2$ acts by $(-1,-1)$. There are two canonical group
homomorphisms
$$
\tau:\spinc(V)\to \So(V)\quad ,\quad \Det :\spinc(V)\to U(1)\
$$
such that  $\tau^{\rm c}=(\tau,\Det):\spinc(V)\to \So(V)\times U(1)$
is a double covering map.

\begin{defi}The $K$-equivariant line bundle
$\Lbb_\Vcal:=P^{\rm c}\times_{\Det}\Cbb$ over $M$ is called the determinant line
bundle associated to the $\spinc$ structure on $\Vcal$.
\end{defi}

Let $\nabla^{\Lbb_\Vcal}$ be an invariant connection on
$\Lbb_\Vcal$ adapted to an invariant Hermitian metric. Let
$\F^{\Lbb_\Vcal}(X), X\in \kgot$ be its equivariant curvature
$2$-form. Even if the line bundle $\Lbb_\Vcal$ does not admit a
square root, we define (formally) the Chern character of  the square
root as follows.

\begin{defi}The Chern character
$\ch(\Lbb^{1/2}_\Vcal)\in \Hcal^{\infty}(\kgot,M)$ is defined
by the equivariant form $\e^{\frac{1}{2}\F^{\Lbb_\Vcal}(X)}$.
\end{defi}

Since the spinor representation extends to $\spinc(V)$, the $\spinc$ structure
on $\Vcal$ induces a $K$-equivariant spinor bundle $\Scal^{\rm c}_M:=
P^{\rm c}\times_{G^{\rm c}} \Scal$ on $M$. Like in the $\spin$ case, one considers
the $K$-equivariant morphism
$\sigma^{\rm c}_\Vcal:p^*\Scal^{{\rm c},+}_M\longrightarrow p^*\Scal^{{\rm c},-}_M$,
defined
by $\sigma^{\rm c}_\Vcal=-i\, \clif(\xx)$
where $\xx:\Vcal\to p^*\Vcal$ is the canonical section.

\begin{prop}\label{prop:ch-spinc-Thom}
We have the equalities
$$
\chr(\sigma^{\rm c}_\Vcal)= (2i\pi)^n
p^*\left(\jdemi(\Vcal)\ch(\Lbb^{1/2}_\Vcal)\right)
\tur(\Vcal)\quad \mathrm{in} \quad
\Hcal^{\infty}(\kgot,\Vcal,\Vcal\setminus M).
$$
$$
\chcf(\sigma^{\rm c}_\Vcal)= (2i\pi)^n
p^*\left(\jdemi(\Vcal)\ch(\Lbb^{1/2}_\Vcal)\right)
\tucf(\Vcal)\quad \mathrm{in} \quad
\Hcal^{\infty}_{\cf}(\kgot,\Vcal).
$$
$$
\ch_Q(\sigma^{\rm c}_\Vcal)= (2i\pi)^n
p^*\left(\jdemi(\Vcal)\ch(\Lbb^{1/2}_\Vcal)\right)
\tumq(\Vcal)\quad \mathrm{in} \quad
\Hcal^{\infty}_{\dr}(\kgot,\Vcal).
$$
\end{prop}

\begin{proof} It is an easy matter to extend the proof of the $\spin$ case.
The Lie algebra of $\spingot^{\rm c}(V)$ of $\spinc(V)$ is  identified with
$\spingot(V)\times \Rbb$.

First one consider the the case of an oriented Euclidean vector
space $V$ of dimension $2n$ equipped with $\spinc(V)$-equivariant vector bundles
$\Scal^{{\rm c},\pm}_V:=V\times \Scal^{{\rm c},\pm}$ over $V$. Recall that the
action of $\spinc(V)$ on the base $V$ is through
$\tau^{\rm c} : \spinc(V)\to \so(V)$. Here $\Scal^{{\rm c},\pm}$ are the
spinor spaces $\Scal^{\pm}$ but with the (extended) action of $\spinc(V)$.

Then we consider the $\spinc(V)$-equivariant morphism
$\sigma^{\rm c}_V:\Scal^{{\rm c},+}_V\longrightarrow \Scal^{{\rm c},-}_V$,
defined by $\sigma^{\rm c}_\Vcal=-i\, \clif(\xx)$.
The $\spinc(V)$-equivariant curvature of the super-connection
$\A_t=  d + t\textbf{c}(x)$ is
\begin{equation}\label{eq:curvature-spinc}
    \F_t(Y,\theta)= \F_t(Y) +i\theta,\quad (Y,\theta)\in \spingot^{\rm c}(V),
\end{equation}
where $\F_t(Y)$ is the $\spin(V)$-equivariant curvature compute
in (\ref{eq:curvature-spin}). Then we get the easy extension of Proposition
\ref{prop:chern-V-spin}

\begin{lem}\label{lem:chern-V-spinc}
We have the following equality in $\Acal^\infty_c(\spingot^{\rm
c}(V),V)$
\begin{equation}\label{eq:chern-V-spinc}
\ch_c(\sigma^{\rm c}_V)(Y,\theta) =(2i\pi)^n \, \e^{i\theta}\,
\jdemi_V(Y^\tau) \tuc(V)(Y^\tau).
\end{equation}
\end{lem}

We come back to the situation of the $K$-equivariant vector bundle $\Vcal\to M$.
We consider a $K$-invariant connection one form $\omega$ on the
$\spinc(V)$-principal bundle $P^{\rm c}\to M$.

We prove Proposition \ref{prop:ch-spin-Thom} after taking the image
of (\ref{eq:chern-V-spinc}) by the Chern-Weil homomorphism
$\phi_\omega$. Note that
$\phi_\omega(\e^{i\theta})=\ch(\Lbb^{1/2}_\Vcal)$ in
$\Hcal^{\infty}(\kgot,M)$.

\end{proof}\bigskip

\subsubsection{The complex case}

In this section we treat the special case where the $\spinc$
structure comes from a complex structure.

We assume that $p:\Vcal\to M$ is a  $K$-equivariant complex vector
bundle equipped with compatible Hermitian inner product $\langle -,-\rangle$ and
connection $\nabla^\Vcal$. We consider the super-vector bundle
$\wedge_\Cbb \Vcal\to M$, where $\wedge_\Cbb \Vcal\to M$ means that we
consider $\Vcal$ as a complex vector bundle.

We consider now the $K$-equivariant morphism on $\Vcal$,
$\sigma^\Cbb_\Vcal:p^*(\wedge_\Cbb^+ \Vcal)\longrightarrow
p^*(\wedge_\Cbb^- \Vcal)$, defined by
$$
\sigma^\Cbb_\Vcal(v)=-i\left( \iota(v) - \varepsilon(v) \right)\quad
{\rm on} \quad \wedge_\Cbb^+ \Vcal_{p(v)},
$$
where $\iota(v)$ and $\varepsilon(v)$ are respectively the contraction
by $v$ ($\iota(v)(w)=\langle w,v\rangle$) and the wedge product by $v$.
Since $\left( \iota(v) - \varepsilon(v) \right)^2= -\|v\|^2$,
we know that the support of $\sigma^\Cbb_\Vcal$
is the zero section of $\Vcal$.

Let $\F^\Vcal(X), X\in\kgot$ be the equivariant curvature of $\nabla^\Vcal$.
\begin{defi}
The equivariant Todd form of $(\Vcal,\nabla^\Vcal)$ {\em is defined by
for $X$ small}, by $$
\todd(\nabla^\Vcal)(X):=\det{}_\Cbb
\left(\dfrac{\F^\Vcal(X)}{\e^{\F^\Vcal(X)}-1}\right).
$$
We denote $\todd(\Vcal)$ its cohomology class in $\Hcal^{\infty}(\kgot,M)$.
\end{defi}

Remark that the inverse $\det{}_\Cbb
\left(\frac{\e^{\F^\Vcal(X)}-1}{\F^\Vcal(X)}\right)$ of the equivariant Todd form
 is defined for any
$X\in\kgot$.

\begin{prop}

We have the equalities:
$$
\chr(\sigma^\Cbb_\Vcal)= (2i\pi)^n
p^*\left(\todd(\Vcal)^{-1}\right) \tur(\Vcal)\quad \mathrm{in}
\quad \Hcal^{\infty}(\kgot,\Vcal,\Vcal\setminus M),
$$
$$
\chcf(\sigma^\Cbb_\Vcal)= (2i\pi)^n
p^*\left(\todd(\Vcal)^{-1}\right) \tucf(\Vcal)\quad \mathrm{in}
\quad \Hcal^{\infty}_{\cf}(\kgot,\Vcal),
$$
$$
\ch_Q(\sigma^\Cbb_\Vcal)= (2i\pi)^n
p^*\left(\todd(\Vcal)^{-1}\right) \tumq(\Vcal)\quad \mathrm{in}
\quad \Hcal^{\infty}_{\dr}(\kgot,\Vcal).
$$
\end{prop}

\begin{proof} The complex structure on the bundle $\Vcal$ induces
canonically a $\spinc$ structure where the bundle of spinors is
$\wedge_\Cbb \Vcal$. The corresponding determinant line bundle is
$\Lbb_\Vcal:=\wedge_\Cbb^{\rm max} \Vcal$. Then one has just to
check that
$$
\jdemi(\Vcal)\ch(\Lbb^{1/2}_\Vcal)=  \todd(\Vcal)^{-1},
$$
and we conclude with Proposition \ref{prop:ch-spinc-Thom}.
\end{proof}\bigskip

\section{Appendix}

We give  proofs of the estimates used in this article. They are all
based on  Volterra's expansion formula: if $R$ and $S$ are elements
in a finite dimensional associative algebra, then
\begin{equation}\label{volterra}
\e^{(R+S)}=\e^{R}+\sum_{k=1}^{\infty} \int_{\Delta_k}\e^{s_1 R} S
\e^{s_2 R} S \cdots S \e^{s_{k}R} S \e^{s_{k+1}R}ds_1\cdots ds_{k}
\end{equation}
where $\Delta_{k}$ is the simplex $\{s_i\geq 0;
s_1+s_2+\cdots+s_{k}+s_{k+1}=1\}.$ We recall that the volume of
$\Delta_{k}$  for the measure $ds_1\cdots ds_{k} $ is
$\frac{1}{k!}$.

Now, let $\Acal=\oplus_{i=0}^q \Acal_i$  be a finite dimensional
graded commutative algebra with a norm $\|\cdot\|$ such that
$\|ab\|\leq \|a\|\|b\|$. We assume $\Acal_0=\Cbb$ and we denote by
$\Acal_+=\oplus_{i=1}^q \Acal_i$. Thus $\omega^{q+1}=0$ for any
$\omega \in \Acal_+.$ Let $V$ be a finite dimensional  Hermitian
vector space. Then $\End(V)\otimes \Acal$ is an algebra with a norm
still denoted by $\|\cdot\|$. If $S\in \End(V)$, we denote also by
$S$ the element $S\otimes 1$ in $\End(V)\otimes \Acal$.

\begin{rem} In the rest of this section we will denote $\cst(a,b,\cdots)$ some positive constant
which depends on the parameter $a,b,\cdots$.
\end{rem}

\subsection{First estimates}

We denote $\herm(V)\subset\End(V)$ the subspace formed by the
Hermitian endomorphisms. When $R\in\herm(V)$, we denote
$\sm(R)\in\Rbb$ the smallest eigenvalue of $R$ : we have
$$
\Big|\!\Big| \e^{-R}\Big|\!\Big|=\e^{-\sm(R)}.
$$

\begin{lem}\label{suffit}
Let $\Pcal(t)=\sum_{k=0}^q \frac{t^k}{k!}$. Then, for any $S\in
\End(V)\otimes \Acal$, $T\in \End(V)\otimes \Acal_+$, and
$R\in\herm(V)$, we have
$$
\|\e^{-R+S+T}\|\leq \e^{-\sm(R)}\e^{\|S\|}\Pcal(\|T\|).
$$
\end{lem}

\begin{proof} Let $c=\sm(R)$.
Then $\|\e^{-u R}\|= \e^{-uc}$ for all $u\geq 0$. Using Volterra's
expansion for the couple  $sR,s S$, we obtain $\|\e^{s(-R+S)}\|\leq
\e^{-sc}\e^{s\|S\|}.$ Indeed,
$\e^{s(-R+S)}=\e^{-sR}+\sum_{k=1}^{\infty} I_k$ with
$$
I_k= s^k\int_{\Delta_k}\e^{-s_1s R}S\cdots S\e^{-s_ks R}
 S \e^{-s_{k+1}sR}ds_1\cdots ds_k.
$$
The term $I_k$ is bounded in norm by $\frac{s^k}{k!}\|S\|^k\e^{-sc}$.
Summing in $k$, we obtain $\|\e^{-s(R+S)}\|\leq \e^{-sc}\e^{s\|S\|}$
for $s\geq 0$. We reapply Volterra's expansion to compute
$\e^{(-R+S)+T}$ as the sum
$$\e^{-R+S}+\sum_{k\geq 1}^q
\int_{\Delta_k}\e^{s_1(-R+S)}T\cdots T\e^{s_k(-R+S)}
T\e^{s_{k+1}(-R+S)}ds_1\cdots ds_k.
$$

Here the sum in $k$ is finite and stops at $k=q$.  The norm of the
$k^{th}$ term is bounded by $\frac{1}{k!}\e^{-c}\e^{\|S\|}\|T\|^k$.
Summing up in $k$, we  obtain our estimate.
\end{proof}\bigskip

\bigskip

For proving  Proposition \ref{estimates},   we need to consider
the following situation. Let $E$ be a (finite dimensional) vector
space. We consider the following smooth maps
\begin{itemize}
  \item $x\mapsto S(x)$ from $E$ to $\End(V)\otimes\Acal$.
  \item $(t,x)\mapsto t^2 R(x)$ from $\Rbb\times E$ to $\herm(V)$.
  \item $(t,x)\mapsto T(t,x)=T_0(x)+t T_1(x)$ from $\Rbb\times E$ to
  $\End(V)\otimes\Acal_+$.
\end{itemize}


\begin{prop}\label{prop-estimation-generale}
Let $D(\partial)$ be a constant coefficient differential operator in
$x\in E$ of degree $r$. Let $\Kcal$ be a compact subset of $E$.
There exists a constant $\cst>0$ (depending on $\Kcal,R(x)$, $S(x),T_0(x),T_1(x)$
 and $D(\partial)$)  such that
\begin{equation}\label{eq:maj-D-exp}
    \Big|\!\Big|D(\partial)\cdot \e^{-t^2R(x)+S(x)+T(t,x)}\Big|\!\Big|
\leq \cst\, (1+t)^{2r+q}\,  \e^{-t^2\sm(R(x))},
\end{equation}
for all $(x,t)\in\Kcal\times \Rbb^{\geq 0}$.

Here the integer $q$
is highest degree of the graded algebra $\Acal$.
\end{prop}

\begin{coro}\label{coro:integration-smooth}

Let $\Ucal$ be an open subset of $E$ such that $R(x)$ is positive
definite for any $x\in\Ucal$, that is $\sm(R(x))>0$ for all $x\in
\Ucal$. Then the integral
$$
\int_0^{\infty}\e^{-t^2R(x)+S(x)+T(t,x)}dt
$$
defines a smooth map from $\Ucal$ into $\End(V)\otimes\Acal$.
\end{coro}

\begin{proof} We fix a basis $v_1,\ldots,v_p$ of $E$. Let us denote
$\partial_i$ the partial derivative along the vector $v_i$. For
any sequence $I:=[i_1,\ldots,i_n]$ of integers
$i_k\in\{1,\ldots,p\}$, we denote $\partial_I$ the differential
operator of order $n=|I|$ defined by the product $\prod_{k=1}^n
\partial_{i_k}$.

For any smooth function $g:E\to\End(V)\otimes\Acal$ we define the
functions
$$
\Big|\!\Big| g \Big|\!\Big|_{n}(x):=\sup_{|I|\leq
n}\ \|\partial_I\cdot g(x)\|
$$
and the semi-norms $\| g\|_{\Kcal,n}:=\sup_{x\in\Kcal} \| g\|_{n}(x)$ attached to a compact subset
$\Kcal$ of $E$. We will use the trivial fact that $\| g\|_{n}(x)\leq \| g\|_{m}(x)$ when $n\leq m$.
Since any constant differential operator $D(\partial)$ is a finite
sum $\sum_I a_I\partial_I$, it is enough to proves
(\ref{eq:maj-D-exp}) for the $\partial_I$.

First, we analyze $\partial_I\cdot\left(\e^{-t^2R(x)}\right)$.
The Volterra expansion formula gives
\begin{equation}\label{eq-derive-1}
\partial_i\cdot\left(\e^{-t^2R(x)}\right)= -t^2
\int_{\Delta_1}\e^{-s_1 t^2 R(x)} \partial_i\cdot R(x)\, \e^{-s_2 t^2
R(x)}ds_1 ,
\end{equation}
and then $\|\partial_i\cdot \e^{-t^2R(x)}\| \leq
\| R\|_{1}(x)\, (1+t)^2\, \e^{-t^2\sm(R(x))}$ for
$(x,t)\in E\times \Rbb^{\geq 0}$.

With (\ref{eq-derive-1}), one can easily prove by induction on the
degree of $\partial_I$ that: if $| I  |=n$  then
\begin{equation}\label{eq:maj-D-exp-R}
    \Big|\!\Big|\partial_I\cdot \e^{-t^2R(x)}\Big|\!\Big|
\leq \cst(n)\, \Big(1+\|R\|_{n}(x)\Big)^n\, (1+t)^{2n}\,  \e^{-t^2\sm(R(x))}
\end{equation}
for $(x,t)\in E\times \Rbb^{\geq 0}$. Note that (\ref{eq:maj-D-exp-R}) is still true
when $I=\emptyset$ with $\cst(0)=1$.

Now we look at  $\partial_I\cdot\left(\e^{-t^2R(x)+S(x)}\right)$ for
 $| I  |=n$. The
Volterra expansion formula gives $\e^{-t^2R(x)+S(x)}= \e^{-t^2R(x)}+
\sum_{k=1}^{\infty} \mathcal{Z}_k(x)$ with
$$
\mathcal{Z}_k(x)=\int_{\Delta_k}\e^{-s_1(t^2R(x))}S(x) \e^{-s_2(t^2R(x))}S(x)
\cdots S(x) \e^{-s_{k+1}(t^2 R(x))}ds_1\cdots ds_k.
$$
The term $\partial_I\cdot \mathcal{Z}_k(x)$ is  equal to the sum,
indexed by the partitions (we allow some of the $I_j$ to
be empty.) $\Pcal:=\{I_1,I_2,\ldots, I_{2k+1}\}$ of $I$, of the
terms
\begin{equation}\label{eq-Z-I}
\mathcal{Z}_k(\Pcal)(x):=
\end{equation}
$$
\int_{\Delta_k}\left(\partial_{I_1}\!\cdot\! \e^{-s_1(t^2R(x))}\right)
\left(\partial_{I_2}\!\cdot\! S(x) \right)\cdots
\left(\partial_{I_{2k}}\!\cdot \!S(x) \right)
\left(\partial_{I_{2k+1}}\!\cdot \!\e^{-s_{k+1}(t^2 R(x))}\right)
ds_1\cdots ds_k
$$
which are, thanks to (\ref{eq:maj-D-exp-R}), smaller in norm than
\begin{equation}\label{eq:maj-Z-k-P}
\cst(\Pcal)\,\Big(1+\|R\|_{n^+_\Pcal}(x)\Big)^{n^+_\Pcal}\,
\frac{\left(\|S\|_{n^-_\Pcal}(x)\right)^k}{k!}\,
(1+t)^{2n^+_\Pcal}\,  \e^{-t^2\sm(R(x))}.
\end{equation}
The integer $n^+_\Pcal,n_\Pcal^-$ are respectively  equal to the
sums $| I_1 | + | I_3| + \cdots +| I_{2k+1} |$, $| I_2| + | I_4| +
\cdots +| I_{2k} |$, and then $n^+_\Pcal +n_\Pcal^-=n$. The
constant $\cst(\Pcal)$ is equal to the products $\cst(| I_1
|)\cst(| I_3|)\cdots\cst(| I_{2k+1}|)$.  Since the sum $\sum_\Pcal
\cst(\Pcal)$ is bounded by a constant $\cst'(n)$, we find that
\begin{equation}\label{eq:maj-D-exp-R-S}
    \Big|\!\Big|\partial_I\cdot \e^{-t^2R(x)+S(x)}\Big|\!\Big|
\leq \cst'(n)\,\Big(1+\|R\|_{n}(x)\Big)^n \e^{\|S\|_{n}(x)} (1+t)^{2n}\,  \e^{-t^2\sm(R(x))}
\end{equation}
for $(x,t)\in E\times \Rbb^{\geq 0}$. Note that (\ref{eq:maj-D-exp-R-S}) is still true
when $I=\emptyset$ with $\cst'(0)=1$.

Finally we look at  $\partial_I\cdot\left(\e^{-t^2R(x)+S(x)+T(t,x)}\right)$ for
 $| I  |=n$. The
Volterra expansion formula gives $\e^{-t^2R(x)+S(x)+T(t,x)}=
\e^{-t^2R(x)+S(x)}+\sum_{k=1}^{q} \mathcal{W}_k(x)$ with
$$
\mathcal{W}_k(x)=\int_{\Delta_k}\e^{s_1(-t^2R(x)+S(x))}T(t,x)
\cdots T(t,x) \e^{s_{k+1}(-t^2 R(x)+S(x))}ds_1\cdots ds_k.
$$
Note that the term $ \mathcal{W}_k(x)$ vanishes for $k>q$. If we use
(\ref{eq:maj-D-exp-R-S}), we get for $(x,t)\in E\times \Rbb^{\geq 0}$ :
\begin{eqnarray*}
\lefteqn{\|\partial_I \cdot  \mathcal{W}_k(x)\|\leq
\cst''(n)\Big(\|T_0\|_{n}(x)+\|T_1\|_{n}(x)\Big)^k \times} \\
& &\Big(1+\|R\|_{n}(x)\Big)^n
\frac{(1+t)^{2n+k}}{k!}\,\e^{\|S\|_{n}(x)}\,  \e^{-t^2\sm(R(x))}.\nonumber
\end{eqnarray*}
Finally we get  for $(x,t)\in E\times \Rbb^{\geq 0}$ :
\begin{eqnarray}\label{eq:maj-D-exp-R-S-T}
\lefteqn{\quad \Big|\!\Big|\partial_I\cdot \e^{-t^2R(x)+S(x)+ T(t,x)}\Big|\!\Big|
\leq \cst''(n) \Big(1+\|R\|_{n}(x)\Big)^n\times }\\
& &\, {\rm P}\Big(\|T_0\|_{n}(x)+\|T_1\|_{n}(x)\Big)
\,\e^{\|S\|_{n}(x)}\,  (1+t)^{2n+q} \, \e^{-t^2\sm(R(x))}\nonumber
\end{eqnarray}
where ${\rm P}$ is the polynomial ${\rm P}(z)=\sum_{k=0}^q
\frac{z^k}{k!}$.

So (\ref{eq:maj-D-exp}) is proved with
$$
\cst= \cst''(n) \sup_{x\in \Kcal}
\left\{\Big(1+\|R\|_{n}(x)\Big)^n{\rm P}\Big(\|T_0\|_{n}(x)+\|T_1\|_{n}(x)\Big)
\,\e^{\|S\|_{n}(x)}\right\}.
$$

\end{proof}\bigskip

\subsection{Second estimates}\label{subsec:appendix2}

Consider now the case where $E=W\times \kgot$ : the variable $x\in E$ will be replaced
by $(y,X)\in W\times \kgot$. {\bf We suppose that the maps $R$ and $T$ are constant
relatively to the parameter $X\in \kgot$.}

Let $\Kcal=\Kcal'\times\Kcal''$ be a compact
subset of $W\times\kgot$. Let $D(\partial)$ be a constant coefficient differential operator in
$(y,X)\in W\times \kgot$ of degree $r$ : let $r_W$ be its degree relatively to the variable $y\in W$.

\begin{prop}\label{prop-estimation-generale-X}
There exists a constant $\cst>0$, depending on $\Kcal,R(y),S(y,X)$,\break
$T_0(y),T_1(y)$ and $D(\partial)$,  such that
\begin{equation}\label{eq:maj-D-exp-bis}
    \Big|\!\Big|D(\partial)\cdot \e^{-t^2R(y)+S(y,X)+T(t,y)}\Big|\!\Big|
\leq \cst\, (1+t)^{2r_W+q}\,  \e^{-t^2\sm(R(y))},
\end{equation}
for all $(y,X,t)\in\Kcal'\times\Kcal''\times \Rbb^{\geq 0}$.
Here $q$
is the highest degree of the graded algebra $\Acal$.
\end{prop}

\begin{proof} We follow the proof of Proposition \ref{prop-estimation-generale}. We have just to explain
why we can replace in (\ref{eq:maj-D-exp}) the factor $(1+t)^{2 r}$ by $(1+t)^{2 r_W}$.

We choose some basis $v_1,\ldots,v_{\pi_1}$ of $W$ and
 $X_1,\ldots,X_{\pi_2}$ of $\kgot$.  Let us denote
$\partial_i^1,\partial_j^2$ the partial derivatives along the vector $v_i$ and
$X_j$.  For any sequence
$$
I:=\underbrace{\{i_1,\ldots,i_{n_1}\}}_{I(1)}\cup \underbrace{\{
j_1,\ldots,j_{n_2}\}}_{I(2)}
$$
of integers where $i_k\in\{1,\ldots,\pi_1\}$ and
$j_k\in\{1,\ldots,\pi_2\}$, we denote $\partial_I$ the differential
operator of order $|I|=n_1+n_2$ defined by the product
$\prod_{k=1}^n \partial_{i_k}^1\prod_{l=1}^m \partial^2_{j_k}$.

We first notice that $\partial_I\cdot \e^{-t^2R(y)}=0$ if
$I(2)\neq\emptyset$. Now we look at  $\partial_I\cdot\left(\e^{-t^2R(y)+S(y,X)}\right)$ for
 $I=I(1)\cup I(2)$.  The term
 $\mathcal{Z}_k(\Pcal)$  of (\ref{eq-Z-I}) vanishes when there exists
 a sub-sequence $I_{2l+1}$ with $I_{2l+1}(2)\neq\emptyset$.
 In the other cases, the integer $n^+_\Pcal=| I_1 | + | I_3| + \cdots +| I_{2k+1}|$ appearing in
 (\ref{eq:maj-Z-k-P}) is smaller than
 $| I(1)  |=n_1$.  So the inequalities (\ref{eq:maj-D-exp-R-S}) and (\ref{eq:maj-D-exp-R-S-T})
hold with the factor
 $(1+t)^{2n}$ replaced by  $(1+t)^{2n_1}$.

\end{proof}

\bigskip

The preceding estimates hold if we work in the algebra
$\End(\Ecal)\otimes \Acal$, where $\Ecal$ is a super-vector space
and $\Acal$ a super-commutative algebra.


\end{document}